\documentclass[12pt]{report}
\usepackage{amsmath}
\usepackage[cp866]{inputenc}
\usepackage[english,russian]{babel}
\usepackage{amsfonts,amssymb}
\usepackage{latexsym}
\usepackage{personal}

\usepackage[english,russian]{babel} 
\usepackage{amssymb,amsmath} 

\newtheorem{thm}{Теорема}[chapter]
\newtheorem{lem}{Лемма}[chapter]
\newtheorem{rem}{Замечание}[chapter]
\newtheorem{exa}{Пример}[chapter]
\newtheorem{pro}{Предложение}[chapter]
\newtheorem{cor}{Следствие}[chapter]
\newtheorem{opr}{Определение}[chapter]

\begin{document} \large
\setcounter{equation}{0}
\thispagestyle{empty}
\begin{center}
ВОРОНЕЖСКИЙ ГОСУДАРСТВЕННЫЙ УНИВЕРСИТЕТ\\

\vskip10mm

\hskip90mm на правах рукописи

\vskip15mm \Large МАКАРЕНКОВ \ ОЛЕГ \ ЮРЬЕВИЧ

\vskip10mm {\huge \bf{Методы теории топологической степени в
задачах
И.~\hskip-1.8mmГ.~\hskip-1.8mmМалкина~-~В.~\hskip-1.8mmК.~\hskip-1.8mmМельникова
для периодически возмущенных систем}}

\vskip17mm \rm

01.01.02 --- дифференциальные уравнения

\vskip7mm

Диссертация на соискание ученой степени\\
кандидата физико-математических наук
\end{center}
\vskip10mm

\hskip60mm \hfill Научный руководитель

\hskip60mm \hfill доктор физико-математических наук,

\hskip60mm \hfill профессор  М.~И.~Каменский

\vskip40mm
\begin{center}
Воронеж -- 2006
\end{center}


\setcounter{secnumdepth}{4} \setcounter{tocdepth}{2}

\tableofcontents

\pagebreak

\chapter*{Введение}
\addcontentsline{toc}{chapter}{Введение}

 Топологическая степень $
  d_{\mathbb{R}^n}(F,U)
$ векторного поля $F:\mathbb{R}^n\to\mathbb{R}^n$ по отношению к
открытому множеству $U\subset\mathbb{R}^n$ в cлучае односвязного
множества $U,$ ограниченного положительно ориентированной
жордановой кривой $q,$ и $n=2$ введена А.~Пуанкаре  и известна под
названием индекса кривой $q$ по отношению к полю $F$ (см.
\cite{poin}, Гл.~3). А.~Пуанкаре использовал полученную
характеристику для анализа существования, числа и типа особых
точек двумерных автономных систем. К нему же восходит основная
теорема теории топологической степени: {\it если
$d_{\mathbb{R}^2}(F,U)\not=0,$ то в $U$ имеется особая точка поля
$F$} и свойство аддитивности топологической степени (первое
основное свойство степени), именно, {\it если $U=\overline{U_1\cup
U_2}$ и $U_1\cap U_2=0,$ где $U_1,U_2\subset\mathbb{R}^2$ --
открытые множества, ограниченные положительно ориентированными
жордановыми кривыми, то
$d_{\mathbb{R}^2}(F,U)=d_{\mathbb{R}^2}(F,U_1)+d_{\mathbb{R}^2}(F,U_2)$}
(см. \cite{poin}, с.~38). Также А.~Пуанкаре доказал, что если
множество $U$ содержит простую особую точку (простой нуль) поля
$F$ и достаточно мало, то $\left|d_{\mathbb{R}^2}(F,U)\right|=1$
(в зависимости от того $d_{\mathbb{R}^2}(F,U)=1$ или
$d_{\mathbb{R}^2}(F,U)=-1$ А.~Пуанкаре делал выводы о типе особой
точки), {\it если же в этом малом множестве нет нулей поля $F,$ то
$d_{\mathbb{R}^2}(F,U)=0$ }(второе основное свойство степени) (см.
\cite{poin}, с.~39). Для случая произвольных открытого
ограниченного множества $U\in\mathbb{R}^n$ и $n\in\mathbb{N}$
конструкция топологической степени получена Л.~Брауером
\cite{brouwer}, кто также сформулировал третье основное свойство
топологической степени (принцип продолжения Брауера) о том, что
{\it степень $d_{\mathbb{R}^n}(F,U)$ остается постоянной, если
область $U$ и отображение $F$ непрерывно меняются так, что в образ
границы $\partial U$ этой области нигде не попадает нуль} (см.
\cite{brouwer}, свойство с.~105). Наконец, Ж.~Лерэ и Ю.~Шаудер
рассмотрели случай, когда $F$ является разностью тождественного и
компактного отображений, заданных в банаховом пространстве.
Доказывая возможность аппроксимации этой ситуации некоторой
конечномерной и используя в последней степень Брауера, Ж.~Лерэ и
Ю.~Шаудер обосновали определение степени в банаховом
(бесконечномерном) пространстве (см. \cite{lersch}, \S 1).

В диссертационной работе изучаются возможности применения теории
топологической степени  к задачам И.~Г.~Малкина и В.~К.~Мельникова
о существовании $T$-периодических решений в системе обыкновенных
дифференциальных уравнений
\begin{equation}\label{ps_I}\addtocounter{form}{1}\tag{\arabic{form}}
  \dot x=f(t,x)+\varepsilon g(t,x,\varepsilon),
\end{equation}
где $f:\mathbb{R}\times\mathbb{R}^n\to\mathbb{R}^n,$
$g:\mathbb{R}\times\mathbb{R}^n\times[0,1]\to\mathbb{R}^n$ --
$T$-периодические по первой переменной непрерывные функции и
$\varepsilon>0$ -- малый параметр. К системам вида (\ref{ps_I})
приводится большое число уравнений, описывающих  разнообразные
нелинейные процессы, в частности, уравнения Ван дер Поля,
Дуффинга, "синус Гордона"\ в отсутствии демпфирования, плоского
маятника, "хищник-жертва"\ при учете периодического изменения
климата. Одной из наиболее важных рассматриваемых при этом задач
 является задача о существовании в системе
(\ref{ps_I}) $T$-периодических решений. Аналитические методы
решения поставленной задачи, как правило, предполагают, что правые
части системы (\ref{ps_I}) некоторое число раз непрерывно
дифференцируемы, а также, что известно семейство
$\left\{\widetilde{x}_\lambda\right\}_{\lambda\in\Lambda}$
$T$-периодических решений порождающей системы
\begin{equation}\label{np_I}\addtocounter{form}{1}\tag{\arabic{form}}
  \dot x=f(t,x).
\end{equation}
Одним из основных аналитических методов является основанный на
теореме о неявной функции метод малого параметра Пуанкаре (см.
Б.~П.~Демидович \cite{dem}, Гл.~III, \S~24, М.~Розо \cite{rozo},
Гл.~9, \S~1), развитием которого для различных ситуаций занимались
Л.~С.~Понтрягин \cite{pont1}, А.~А.~Андронов-А.~Витт \cite{and},
Н.~Г.~Булгаков \cite{bul}, Н.~М.~Крылов - Н.~Н.~Боголюбов -
Ю.~А.~Митропольский \cite{bog}, А.~М.~Кац \cite{kac}, И.~Г.~Малкин
\cite{mal} , В.~К.~Мельников \cite{mel} и другие. В работах всех
указанных авторов строится соответствующая задаче бифуркационная
функция $M,$ и предъявляется условие о существовании у этой
функции простого нуля $\lambda_0\in\Lambda,$ то есть такого числа,
что $M(\lambda_0)=0$ и $M'(\lambda_0)\not=0.$ Преимуществом
геометрических методов обычно является то, что они работают в
случае, когда возмущение $g$ всего лишь непрерывно, а также не
требуют нахождения простых нулей бифуркационных функций. Вместо
этого предполагается известным поведение решений системы
(\ref{ps_I}) с начальными условиями, принадлежащими границе
$\partial U$ такого открытого ограниченного множества
$U\subset\mathbb{R}^n,$ для которого указанное поведение легко
устанавливается. Одним из основных геометрических методов
доказательства существования $T$-периодических решений является
принцип неподвижной точки. Наиболее удобное его применение связано
с вычислением топологической степени $d_{\mathbb{R}^n}(I-P,U)$
некоторого вспомогательного оператора
$P:\mathbb{R}^n\to\mathbb{R}^n,$ неподвижные точки которого
совпадают с начальными условиями $T$-периодических решений системы
(\ref{ps_I}), относительно множества $U$ и с проверкой отличия
этой топологической степени от нуля. В качестве вспомогательного
оператора используется оператор Пуанкаре
$\mathcal{P}_\varepsilon:\mathbb{R}^n\to\mathbb{R}^n,$ ставящий в
соответствие каждой точке $\xi$ значение единственного решения $x$
системы (\ref{ps_I}) с начальным условием $x(0)=\xi$ в момент
времени $T.$

Первая формула для вычисления топологической степени
$d_{\mathbb{R}^n}(I-\mathcal{P}_\varepsilon,U)$ для систем типа
(\ref{ps_I}) получена М.~А.~Красносельским и А.~И.~Перовым (см.
\cite{kraper} и \cite{perov}) и связана с развитием результата
И.~Берштейна - А.~Халаная \cite{berhal}. Он основан на
предположении о том, что множество $U$ удовлетворяет условию
невозвращаемости, то есть из границы $\partial U$ этого множества
исходят (в нулевой момент времени) только такие решения, которые
не пересекают $\partial U$ при всех $t\in(0,T].$ В этом случае
М.~А.~Красносельским и А.~И.~Перовым установлена формула
$d_{\mathbb{R}^n}(I-\mathcal{P}_\varepsilon,U)=d_{\mathbb{R}^n}(-f(0,\cdot),U),$
позволившая легко считать
$d_{\mathbb{R}^n}(I-\mathcal{P}_\varepsilon,U)$ и доказывать
существование $T$-периодических решений для (\ref{ps_I}) во многих
задачах, где метод малого параметра Пуанкаре ответа не дает,
включая все те, где функция $g$ всего лишь непрерывна. Модификация
формулы Красносельского-Перова для так называемых $m$-систем
получена Э.~Мухамадиевым \cite{muh2}, при этом в левой части
рассматриваемой формулы вместо ограниченного множества $U$ берется
некоторое бесконечно большое множество. Последним принципиальным
результатом в этом направлении является работа А.~Капетто,
Ж.~Мавена и Ф.~Занолина \cite{camaza}, где установлено, что если
система (\ref{np_I}) автономна, то для справедливости формулы
Красносельского-Перова достаточно требовать, чтобы $\partial U$ не
содержало начальных условий $T$-периодических решений системы
(\ref{np_I}).

Вторая формула для вычисления топологической степени
$d_{\mathbb{R}^n}(I-\mathcal{P}_\varepsilon,U)$ получена
Ж.~Мавеном \cite{mawphd1} и предполагает, что $f=0.$ Ж.~Мавен
установил, что если соответствующий оператор усреднения $\Phi^T$
Крылова-Боголюбова-Митропольского невырожден на $\partial U,$ то
$d_{\mathbb{R}^n}(I-\mathcal{P}_\varepsilon,U)=d_{\mathbb{R}^n}(-\Phi^T,U),$
не смотря на то, что $d_{\mathbb{R}^n}(I-\mathcal{P}_0,U)$ для
рассматриваемой системы не определено. Полученная формула
позволила доказать существование $T$-периодических решений во
многих таких системах, где условия аналитических методов
А.~А.~Андронова, Н.~Г.~Булгакова, Н.~М.~Крылова, Н.~Н.~Боголюбова
и Ю.~А.~Митропольского не выполнены. Варианты указанной формулы
для различных случаев, в которых условия Ж.~Мавена не выполнены,
получены М.~И.~Каменским \cite{kam} на основе теоремы
Красносельского-Крейна \cite{krakre} (см. также \cite{dem}, Гл.~V,
\S~3) о предельном переходе под знаком интеграла. Развитие
последней теоремы для систем с наследственностью сделано в работе
В.~В.~Стрыгина \cite{str}, что позволило обосновать формулу Мавена
и для таких систем.

Н.~А.~Бобылев и М.~А.~Красносельский заметили \cite{bob}, что ни
одна из рассмотренных формул не дает геометрического метода
решения задач И.~Г.~Малкина \cite{mal} и В.~К.~Мельникова
\cite{mel} о существовании для системы (\ref{ps_I})
$T$-периодических решений с начальными условиями, принадлежащими
окрестности $T$-периодического цикла $\widetilde{x}$ порождающей
системы (\ref{np_I}) в случае, когда последняя автономна
(И.~Г.~Малкиным рассматривался случай изолированного цикла
$\widetilde{x},$ а В.~К.~Мельниковым случай, когда цикл
$\widetilde{x}$ вложен в некоторое семейство циклов порождающей
системы). Указанное замечание обусловлено тем, что выбирая
множество $U$ лежащим в окрестности цикла $\widetilde{x}$ и
удовлетворяющим условиям невозвращаемости, как правило, имеем
равенство $d_{\mathbb{R}^n}(-f(0,\cdot),U)=0.$ Использование же
формулы Мавена возможно только при дополнительном предположении о
том, что система (\ref{ps_I}) приводится к такой $T$-периодической
системе типа (\ref{ps_I}), в которой $f=0$ (см. К.~Шнайдер
\cite{sch}). Последнее возможно в единственном случае, когда
система (\ref{np_I}), линеаризованная на $\widetilde{x},$ имеет
только $T$-периодические решения, что естественным образом
выполнено лишь для линейных систем (\ref{np_I}). Возникает
естественная проблема: разработать формулы вычисления
топологической степени
$d_{\mathbb{R}^n}(I-\mathcal{P}_\varepsilon,U)$ для более широких
классов множеств $U$ и порождающих систем (2), которые позволили
бы получить геометрические методы решения задач И.~Г.~Малкина и
В.~К.~Мельникова  с одной стороны и превращались бы в формулы
Красносельского-Перова и Мавена в рассмотренных ими ситуациях с
другой стороны. Возможный вариант решения сформулированной
проблемы предлагается в настоящей диссертационной работе.

Актуальность разработки указанных геометрических аналогов связана
еще и с тем, что целый ряд полученных в последнее время
математических моделей приводит к системам (\ref{ps_I}), в которых
функция $g$ не дифференцируема, например, асимметрический
осциллятор Е.~Н.~Дансера \cite{dancer}, модель  колебаний
подвесных мостов А.~С.~Лазера-П.~Дж.~Маккенна \cite{laz} и другие.

 Диссертация состоит из трех глав. В первой главе
рассматривается случай, когда система (\ref{ps_I}) имеет вид
\begin{equation}\label{ps_O}\addtocounter{form}{1}\tag{\arabic{form}}
  \dot x=f(x)+\varepsilon g(t,x,\varepsilon),
\end{equation}
и  предполагается, что граница $\partial U$ множества
$U\subset\mathbb{R}^n$ содержит конечное число начальных условий
$T$-периодических решений автономной порождающей системы
\begin{equation}\label{np_O}\addtocounter{form}{1}\tag{\arabic{form}}
  \dot x=f(x).
\end{equation}
Сначала разрабатывается формула для вычисления топологической
степени $d_{\mathbb{R}^n}(I-\mathcal{P}_\varepsilon,U)$ оператора
Пуанкаре $\mathcal{P}_\varepsilon$ системы (\ref{ps_O}), затем
даются приложения этой формулы к задаче о существовании в системе
(\ref{ps_O}) $T$-периодических решений с принадлежащими множеству
$U$ начальными условиями. Хотя при этом предполагается, что
оператор Пуанкаре $\mathcal{P}_\varepsilon$ для рассматриваемой
системы определен при всех достаточно малых $\varepsilon>0$ (то
есть выполнено условие $(A_{\mathcal{P}})$ единственности и
продолжимости на всю ось решений возмущенной системы с любым
начальным условием), в главе даются аналоги полученных теорем для
случая, когда указанное предположение не выполнено. В этом
последнем случае вместо оператора Пуанкаре
$\mathcal{P}_\varepsilon$ используется интегральный оператор
$$
(Q_\varepsilon x)(t)=x(T)+\int\limits_0^t
f(x(\tau))d\tau+\varepsilon\int\limits_0^t
g(\tau,x(\tau),\varepsilon)d\tau, \quad t\in[0,T],
$$
и вместо топологической степени Брауера
$d_{\mathbb{R}^n}(I-\mathcal{P}_\varepsilon,U)$ -- топологическая
степень Лерэ-Шаудера
$d_{C([0,T],\mathbb{R}^n)}(I-Q_\varepsilon,W_U),$ где множество
$W_U\in C([0,T],\mathbb{R}^n)$ выбирается таким образом, чтобы
множества $U$ и $W_U$ имели так называемую общую сердцевину (см.
\cite{krazab}, Гл.~3, \S~24) по отношению к $T$-периодическим
решениям системы (\ref{ps_O}).

Основным ограничением, используемым в главе~1, является
предположение о том, что каждый $T$-периодический цикл
$\widetilde{x}$ системы (\ref{np_O}) с начальным условием из
$\partial U$ является простым, то есть алгебраическая кратность
мультипликатора $+1$ системы
\begin{equation}\label{ls_O}\addtocounter{form}{1}\tag{\arabic{form}}
  \dot y=f'(\widetilde{x}(t))y
\end{equation}
равна $1,$ что соответствует требованиям работы И.~Г.~Малкина
\cite{mal}. В диссертации показано, что вклад каждого такого цикла
в величину топологической степени
$d_{\mathbb{R}^n}(I-\mathcal{P}_\varepsilon,U)$ может быть
посчитан (теорема~1.4) при помощи соответствующих бифуркационных
функций Малкина
$$
M_{\widetilde{x}}(\theta) = {\rm sign}\left<\dot
{\widetilde{x}}(0),\widetilde{z}(0)\right>\int\limits_0^T\left<\widetilde{z}(\tau),
g(\tau-\theta,\widetilde{x}(\tau),0)\right>d\tau,
$$
где $\widetilde{z}$ -- произвольное нетривиальное
$T$-периодическое решение системы $
  \dot z=-(f'_x(\widetilde{x}(t)))^*z.
$ В случае, когда $\partial U$ не содержит начальных условий
$T$-периодических циклов системы (\ref{np_O}), установленная в
теореме~1.4 формула~(\ref{kmnformU}) для  вычисления
$d_{\mathbb{R}^n}(I-\mathcal{P}_\varepsilon,U)$  совпадает с
формулой Красносельского-Перова. Однако, в покрываемом
теоремой~1.4 классе множеств $U$ уже имеются такие, использование
которых в формуле~(\ref{kmnformU}) позволяет получить
геометрический вариант решения задачи И.~Г.~Малкина \cite{mal} о
существовании $T$-периодических решений в системах (\ref{ps_O})
вблизи цикла $\widetilde{x}$ (теорема~1.6), которое, согласно
замечанию Бобылева-Красносельского, не может быть получено на
основании формулы Красносельского-Перова.

Во  второй главе рассматриваются системы общего вида (\ref{ps_I})
в предположении, что $\mathcal{P}_0(\xi)=\xi$ для любого
$\xi\in\partial U.$ Оказывается (теорема~2.2), выполнения
указанного предположения для справедливости формулы Мавена
достаточно (из условий Мавена следует, что
$\mathcal{P}_0(\xi)=\xi$ для любого $\xi\in\mathbb{R}^n$), если
только так называемый обобщенный оператор усреднения
$\Phi^s:\mathbb{R}^n\to\mathbb{R}^n$ невырожден на $\partial U$
при всех $s\in[0,T].$ Оператор $\Phi^s$ впервые указан в работе
М.~И.~Каменского-О.~Ю.~Макаренкова-П.~Нистри \cite{kmndan} и
совпадает при $s=T$ с классическим оператором усреднения
Крылова-Боголюбова-Митропольского, входящим в формулу Мавена.
Распространение формулы Мавена на такой значительно более широкий
класс множеств $U$ позволило получить новые теоремы о
существовании для системы (\ref{ps_I}) $T$-периодических решений
вблизи $\partial U$ (теоремы~2.4, 2.5 и 2.6). При этом, в
теоремах~2.5 и 2.6 рассматривается случай системы (\ref{ps_O}),
заданной в пространстве $\mathbb{R}^2,$ и в качестве множества $U$
берется внутренность $T$-периодического цикла $\widetilde{x}$
системы (\ref{np_O}). Указанный выбор множества $U$ вместе с
доказанной в главе~2 формулой  для $\Phi^s(\widetilde {x}(t)),$
дающей разложение вектора $\Phi^s(\widetilde {x}(t))$ по векторам
$\dot{\widetilde{ x}}(t)$ и $\widehat {y}(t)$ (лемма~2.4), где
$\widehat {y}$ -- линейно независимое с $\dot{\widetilde{ x}}$
решение системы (\ref{ls_O}), позволил получить геометрический
метод решения задачи В.~К.~Мельникова (теорема~2.6). Одним из
преимуществ полученного метода, по сравнению с методом Мельникова,
является то, что он дает существование для возмущенной системы
(\ref{ps_O}) двух $T$-периодических решений, лежащих по разные
стороны от порождающего цикла $\widetilde{x}.$ Работа
предложенного метода проиллюстрирована на примерах уравнения
Дуффинга (пример~2.1), системы Гринспана-Холмса (пример~2.2) и
одной его модификации, в которой порождающий цикл $\widetilde{x}$
вырожден в том смысле, что все решения системы (\ref{ls_O})
являются $T$-периодическими (пример~2.3). При этом для вычисления
степени $d_{\mathbb{R}^2}(-\Phi^T,U)$ используется разработанный в
этой же главе метод, связанный с некоторыми предположениями типа
четности поля $\Phi^T,$ используемыми в теоремах Борсука-Улама
\cite{borsuk}.

В  третьей главе изучаются свойства $T$-периодических решений
возмущенных систем (\ref{ps_I}) и (\ref{ps_O}), связанные со
скоростью их сходимости при $\varepsilon\to 0.$

Пусть $\{\varepsilon_k\}_{k\in\mathbb{N}}$ -- сходящаяся к нулю
последовательность значений параметра системы (\ref{ps_I}) и
$\{x_k\}_{k\in\mathbb{N}}$ -- соответствующая последовательность
$T$-периодических решений этой системы  такая, что
\begin{equation}\label{as1_3_}\addtocounter{form}{1}\tag{\arabic{form}}
  x_k(t)\to\widetilde{x}(t)\quad\mbox{при}\ \ k\to \infty,
\end{equation}
где $\widetilde{x}$ -- $T$-периодическое решение порождающей
системы (\ref{np_I}). Обозначим через $\Omega(\cdot,t_0,\xi)$
решение $x$ порождающей системы (\ref{np_I}) такое, что
$x(t_0)=\xi.$ Доказанная в главе~3 альтернатива (теорема~3.1)
утверждает, что
 либо
начальные условия $T$-периодических решений системы (\ref{ps_I})
сходятся  к начальному условию $\widetilde{x}(0)$ порождающего
решения $\widetilde{x}$ вдоль плоскости
$\left\{l\in\mathbb{R}^n:\left(\Omega'_{(3)}(T,0,\widetilde{x}(0))-I\right)l=0\right\},$
либо сходимость имеет скорость $\varepsilon>0.$ При этом, в
последнем случае описание поведения решений $x_k$ при $k\to\infty$
может быть уточнено на основании обобщенного оператора усреднения
$\Phi^s.$

Если функция $g$ непрерывно дифференцируема, и свойство
(\ref{as1_3_}) получено применением теорем И.~Г.~Малкина
\cite{mal}, то сходимость в (\ref{as1_3_}) со скоростью
$\varepsilon_k$  уже гарантирована, и теорема~3.1 ничего нового не
дает. Однако, если функция $g$ всего лишь непрерывна, или свойство
(\ref{as1_3_}) получено иными методами, например, методом
Мельникова \cite{mel} или при помощи теорем глав~1 и 2, то теоремы
о скорости сходимости в (\ref{as1_3_}) в литературе отсутствуют, и
теорема~3.1 частично заполняет этот пробел. Полный ответ об
асимптотике расстояния между траекториями решений $x_k$ и
$\widetilde{x}$ в случае, когда функция $g$ непрерывна, дает
теорема~3.2 обсуждаемой главы, но в последней теореме
дополнительно предполагается, что система (\ref{ps_I}) имеет вид
(\ref{ps_O}), и $\widetilde{x}$ является простым циклом. При этом
одно из следствий теоремы~3.1 (следствие~3.7) дает условия, при
которых расстояния между траекториями решений $x_k$ и
$\widetilde{x}$ стремится к нулю со скоростью большей, чем
$\varepsilon_k.$

Каждая глава завершается сопоставлением полученных утверждений с
имеющимися в литературе.

Результаты диссертационной работы докладывались на следующих
семинарах: академика Д.~В.~Аносова и профессора Ю.~С.~Ильяшенко
(Математический институт им. В.~А.~Стеклова РАН, 2006), профессора
Ж.~Мавена (университет г.~Леувен-ла-Нуов, Бельгия, 2005),
профессора  П.~Нистри (университет г.~Сиены, Италия, 2005),
профессора А.~И.~Перова (ВГУ, Воронеж, 2005), профессора Н.~Хирано
(университет г.~Йокогамы, Япония, 2004),
 НОЦ "Волновые
процессы в неоднородных и нелинейных средах" (ВГУ, Воронеж, 2004),
а также на следующих международных конференциях: Barcelona
Conference in Planar Vector Fields (Барселона, Испания, 2006),
"Современные проблемы прикладной математики и математического
моделирования" (Воронеж, 2005), "Trends in Differential Equations
and Dynamical Systems" (Реджио Эмилья, Италия, 2005), "12th
International Workshop on Nonlinear Dynamics of Electronic
Systems" (Евора, Португалия, 2004), "International Symposium on
Dynamical Systems Theory and Its Applications to Biology and
Environmental Sciences" (Хамаматсу, Япония, 2004).

Исследования, включенные в настоящую диссертацию, поддержаны
 грантом РФФИ №~05-01-00100,
 а также грантом для молодых участников проекта VZ--010 "Волновые процессы
 в неоднородных и нелинейных средах"$\;$ Минобразования РФ и CRDF (США).

Основные результаты диссертации опубликованы в работах
\cite{kmndan}, \cite{makope}-\cite{mak2005}, \cite{kmnnon},
\cite{jap}-\cite{barcelona}. Из совместных работ
\cite{kmndan,kmnnon} в диссертацию вошли только принадлежащие
Макаренкову~О.~Ю. результаты.

Автор выражает глубокую благодарность своему научному руководителю
профессору Каменскому Михаилу Игоревичу за постановку задачи,
обсуждение результатов и организацию работы над диссертацией.

\pagebreak

\chapter{Возмущения систем, у которых пересечение множества
начальных условий $T$-периодических решений и границы некоторого
открытого множества $U\subset\mathbb{R}^n$ конечно}\label{gl2}

\setcounter{subsection}{0}

В  настоящей главе исследуется существование $T$-периодических
решений в системах вида
\begin{equation}\label{ps_2}
  \dot x=f(x)+\varepsilon g(t,x,\varepsilon),
\end{equation}
начальные условия которых принадлежат такому открытому
ограниченному множеству $U\subset\mathbb{R}^n,$ граница которого
содержит конечное число начальных условий $T$-периодических
решений порождающей системы
\begin{equation}\label{np_2}
  \dot x=f(x).
\end{equation}

На протяжении всей главы, если другое не оговорено дополнительно,
предполагается, что $f:\mathbb{R}^n\to\mathbb{R}^n$ -- непрерывно
дифференцируемая и $g:\mathbb{R}\times\mathbb{R}^n\times[0,1]$ --
непрерывная функции. Разработке методов решения сформулированной
задачи и их приложениям предпошлем некоторые определения и
свойства, которые будут многократно использоваться на протяжении
этой и последующих глав.

\section{Предварительные сведения}\label{predv}

\paragraph{Основные обозначения.}
  Через $0_{m\times k}$ обозначается нулевая $m\times k$-матрица и
  через
$(a_1,...,a_k),$ где $a_i\in\mathbb{R}^m,$ -- $m\times k$-матрица,
столбцы которой суть векторы $a_1,...,a_k.$ Введенное обозначение
в настоящей диссертации не приводит к путанице с обозначением
интервала $(\theta_1,\theta_2),$ где
$\theta_1,\theta_2\in\mathbb{R}.$ Всюду $\left<\cdot,\cdot\right>$
-- обычное скалярное произведение в $\mathbb{R}^n.$ Если
$\xi\in\mathbb{R}^n$ и $K\subset\mathbb{R}^n$ -- компактное
множество, то $\rho(\xi,K)$ -- расстояние от $\xi$ до $K,$ то есть
$\rho(\xi,K)=\min_{\zeta\in K}\left\|\xi-\zeta\right\|.$
$C([0,T],\mathbb{R}^n)$ - это пространство всех непрерывных
функций, действующих из $[0,T]$ в $\mathbb{R}^n$ с нормой
$\|x\|=\max\limits_{t\in[0,T]}\|x(t)\|.$ Функция $o(\varepsilon),$
используемая в настоящей диссертационной работе, обладает всегда
тем свойством, что $o(\varepsilon)/\varepsilon\to 0$ при
$\varepsilon\to 0.$ При этом, функция $o$ может зависеть и от
других переменных,  такие переменные  в записи $o$ опускаются,
если указанная сходимость к нулю имеет место равномерно по ним.
Через $B_\delta(V)$ обозначается $\delta$-окрестность множества
$V,$ а через $\partial V$ -- его граница в норме содержащего $V$
пространства. Если $f:V\to V_1$ -- некоторая функция, то
$f(V)=\cup_{\xi\in V}f(\xi)$ и $f_{(i)}'$ -- производная функции
$f$ по $i$-й переменной. Тот факт, что равенство $f(\xi)=\xi$
справедливо при всех $\xi\in V$ записывается как $f(\xi)=\xi,\
\xi\in V.$ Запись $\{\xi_1,\xi_2\}$ означает множество, состоящее
из двух элементов $\xi_1$ и $\xi_2.$

\paragraph{Основные определения и свойства.}

Рассмотрим систему обыкновенных дифференциальных уравнений
\begin{equation}\label{Np_0}
\dot x=F(t,x),
\end{equation}
где $F:\mathbb{R}\times\mathbb{R}^n\to\mathbb{R}^n$ -- непрерывная
функция.

Некоторые из приводимых ниже определений  намечены во введении,
сейчас они формулируются со всей строгостью (см.
М.~А.~Красносельский \cite{kraope}).

\begin{opr}\label{shift_ope}
Если решение $x_{t_0,\xi}$ системы (\ref{Np_0}) с начальным
условием $x_{t_0,\xi}(t_0)=\xi$ существует, единственно и
продолжимо на отрезок $[0,T]$ при любых $t_0\in[0,T]$ и
$\xi\in\mathbb{R}^n,$ то оператор $\Omega,$ определяемый для
$(t,t_0,\xi)\in [0,T]\times[0,T]\times\mathbb{R}^n$ как
$$
  \Omega(t,t_0,\xi)=x_{t_0,\xi}(t),
$$
называется оператором сдвига по траекториям системы (\ref{Np_0}).
\end{opr}

\begin{opr}\label{}
Пусть правая часть системы (\ref{Np_0}) $T$-периодична по первой
переменной. Если решение $x$ системы (\ref{Np_0}) с начальным
условием $x(t_0)=\xi$ существует, единственно и продолжимо на
отрезок $[0,T]$ при любых $t_0\in[0,T]$ и $\xi\in\mathbb{R}^n,$ то
оператор
$$
  \mathcal{P}(\xi)=\Omega(T,0,\xi),\quad\xi\in\mathbb{R}^n
$$
называется оператором Пуанкаре $T$-периодической системы
(\ref{Np_0}).
\end{opr}

\begin{opr}\label{}  Если система
(\ref{Np_0}) является $T$-периодической, то интегральным
оператором, соответствующим задаче о $T$-периодических решениях
для (\ref{Np_0}), называется оператор $Q:C([0,T],\mathbb{R}^n)\to
C([0,T],\mathbb{R}^n),$ задаваемый как
$$
  (Qx)(t)=x(T)+\int\limits_0^t F(\tau,x(\tau))d\tau.
$$
\end{opr}

\begin{opr}\label{}
$T$-периодические решения $T$-периодической системы (\ref{Np_0})
будем отождествлять с неподвижными точками интегрального
оператора, соответствующего задаче о $T$-периодических решениях
для (\ref{Np_0}).
\end{opr}

Напомним, что траекторией решения $x$ системы (\ref{Np_0})
называется образ отображения $x:\mathbb{R}\to\mathbb{R}^n$ (см.
В.~И.~Арнольд \cite{ARN}, с.~11).

Систематически будет использоваться нижеследующая
лемма~\ref{zamena}, связанная с возмущенной системой
\begin{equation}\label{ps_0}
\dot x=f(t,x)+\varepsilon g(t,x,\varepsilon),
\end{equation}
где $f:\mathbb{R}\times\mathbb{R}^n\to \mathbb{R}^n$ -- непрерывно
дифференцируемая и $g:\mathbb{R}\times\mathbb{R}^n\times[0,1]\to
\mathbb{R}^n$ -- непрерывная функции, $\varepsilon>0$ -- параметр.
Предположим, что $f(t+T,\xi)=f(t,\xi)$ и
$g(t+T,\xi,\varepsilon)=g(t,\xi,\varepsilon)$ для всех
$t\in\mathbb{R},$ $\xi\in\mathbb{R}^n$ и $\varepsilon\in[0,1].$
Пусть $\Omega$ -- оператор сдвига по траекториям системы
(\ref{ps_0}) при $\varepsilon=0.$

\begin{lem}\label{zamena}
Функция $x\in C([0,T],\mathbb{R}^n)$ является $T$-периодическим
решением системы (\ref{ps_0}) тогда и только тогда, когда функция
\begin{equation}\label{zp_}
  \nu(t)=\Omega(0,t,x(t)),\ t \in [0,T],
\end{equation}
является решением системы
$$
\dot\nu=\varepsilon
\Omega'_{(3)}(0,t,\Omega(t,0,\nu))g(t,\Omega(t,0,\nu),\varepsilon),
$$
удовлетворяющим условию $\nu(0)=\Omega(T,0,\nu(T)).$
\end{lem}

\noindent {\bf Доказательство. } Произведем в системе (\ref{ps_0})
замену переменных
\begin{equation}\label{zpn}
  x(t)=\Omega(t,0,\nu(t)).
\end{equation}
Формула  (\ref{zpn}) каждому $\nu \in C([0,T],\mathbb{R}^n)$
ставит в соответствие $x \in C([0,T],\mathbb{R}^n) $ гомеоморфно,
и обратное отображение дается формулой (\ref{zp_}). Следовательно,
функция $x$ является решением системы (\ref{ps_0}) тогда и только
тогда, когда функция $\nu,$ введенная по закону (\ref{zp_}),
удовлетворяет следующему равенству
\begin{eqnarray}\label{f9}
  \Omega'_{(1)}(t,0,\nu(t))+\Omega'_{(3)}(t,0,\nu(t))\dot \nu(t)=\varepsilon g(t,\Omega(t,0,\nu(t)),\varepsilon)+ \nonumber \\
  +f(t,\Omega(t,0,\nu(t))).
\end{eqnarray}
По определению функции $\Omega$ имеем
\begin{equation}\label{f11}
  \Omega'_{(1)}(t,0,\nu(t))=f(t,\Omega(t,0,\nu(t))).
\end{equation}
Пользуясь соотношением (\ref{f11}), система (\ref{f9}) может быть
переписана в виде
\begin{equation}\label{1_}
  \dot \nu(t)=\varepsilon
  \Omega'_{(3)}(0,t,\Omega(t,0,\nu(t)))g(t,\Omega(t,0,\nu(t)),\varepsilon).
\end{equation}
Рассмотрим произвольное $T$-периодическое решение $x$ системы
(\ref{ps_0}). Имеем
$$
  \nu(0)=\Omega(0,0,x(0))=x(0)=x(T)=\Omega(T,0,\nu(T)).
$$

Лемма доказана.

\section{Связь функций Малкина и топологической степени
 оператора, соответствующего задаче о
$T$-периодических решениях с начальными условиями в $U$}

Пусть $\Omega$ -- оператор сдвига по траекториям системы
(\ref{np_2}). Основное предположение настоящей главы следующее:

$(A_0)$ множество
\[
  \mathfrak{S}^U=\bigcup_{\xi\in\partial
  U:\Omega(T,0,\xi)=\xi}\left\{x\in C([0,T],\mathbb{R}^n):x(t)=\Omega(t,0,\xi),\ t\in[0,T]\right\}
\]
конечно, и для каждого $\widetilde{x}\in \mathfrak{S}^U$
алгебраическая кратность мультипликатора $+1$ системы
\begin{equation}\label{ls_2}
  \dot y=f'_x(\widetilde{x}(t))y
\end{equation}
 равна $1.$ Однако ряд теорем на пути к основному результату
 предполагает более слабые условия.

Предположение $(A_0)$ характерно для случая, когда $U$ является
окрестностью в $\mathbb{R}^n$ некоторой точки
$\widetilde{x}(\theta)$ изолированного цикла $\widetilde{x}$
системы (\ref{np_2}), подробно изученного И.~Г.~Малкиным в
\cite{mal}. В настоящей главе будет, в частности, установлен ряд
обобщений результата И.~Г.~Малкина.

\begin{opr}\label{simplecycle}
Циклы $\widetilde{x}$ порождающей системы (\ref{np_2}) такие, что
алгебраическая кратность мультипликатора $+1$ линейной системы
(\ref{ls_2}) равна $1,$ будем называть простыми.
\end{opr}

Пусть $Q_\varepsilon:C([0,T],\mathbb{R}^n)\to
C([0,T],\mathbb{R}^n)$ -- интегральный оператор, соответствующий
задаче о $T$-периодических решениях для системы (\ref{ps_2}), и
\[
W_U=\left\{\widehat {x}\in
C([0,T],\mathbb{R}^n):\Omega(0,t,\widehat{ x}(t))\in U,\  t\in
[0,T]\right\}.
\]

\noindent Каждому простому циклу $\widetilde{x}$ может быть
поставлена в соответствие бифуркационная функция Малкина (см.
\cite{mal}, формула~3.13)
\begin{equation}\label{fun_mal}
M_{\widetilde{x}}(\theta) = {\rm sign}\left<\dot
{\widetilde{x}}(0),\widetilde{z}(0)\right>\int\limits_0^T\left<\widetilde{z}(\tau),
g(\tau-\theta,\widetilde{x}(\tau),0)\right>d\tau,
\end{equation}
где $\widetilde{z}$ -- произвольное нетривиальное
$T$-периодическое решение системы
\begin{equation}\label{ss_2}
  \dot z=-(f'_x(\widetilde{x}(t)))^*z.
\end{equation}
Всюду ниже через $\beta(\widetilde{x})$ обозначается число, равное
сумме кратностей больших $+1$ мультипликаторов системы
(\ref{ls_2}).

Для достижения основного результата главы теорем~1.3-1.4 нам
понадобится ряд вспомогательных утверждений, которые, однако,
могут иметь самостоятельный интерес для теории топологической
степени. Первое из них следующее.

\begin{thm}\label{MathNach1}
Пусть $\widetilde{x}$ -- простой $T$-периодический цикл системы
(\ref{np_2}). Пусть $0\le\theta_1<\theta_2\le
\theta_1+\frac{T}{p},$ где $p\in\mathbb{N}$ и $\frac{T}{p}$ --
наименьший период цикла $\widetilde{x}.$ Предположим, что
$M_{\widetilde{x}}(\theta_1)\not=0$ и
$M_{\widetilde{x}}(\theta_2)\not=0.$ Тогда для заданного
$\alpha>0$ существует $\delta_0>0$ и семейство открытых множеств
$\{V_\delta\}_{\delta\in(0,\delta_0]},$ удовлетворяющих свойствам:

1) $\widetilde{x}((\theta_1,\theta_2))\subset V_\delta\subset
B_{\delta}(\widetilde{x}((\theta_1,\theta_2))),$

2) $ \partial V_\delta\cap
\widetilde{x}([\theta_1,\theta_2])=\{\widetilde{x}(\theta_1),\widetilde{x}(\theta_2)\},$

\noindent такие, что для любого $\delta\in(0,\delta_0]$ и любого
$\varepsilon\in(0,\delta^{1+\alpha}]$ степень
$d(I-Q_\varepsilon,W_{V_\delta})$ определена и может быть найдена
по следующей формуле
\begin{equation}\label{STAR}
d(I-Q_\varepsilon,W_{V_\delta})=-(-1)^{\beta(\widetilde{x})}d_{\mathbb{R}}
(M_{\widetilde{x}},(\theta_1,\theta_2)).
\end{equation}
\end{thm}

Введем некоторые дополнительные понятия и утверждения необходимые
для доказательства теоремы. Пусть $\widetilde{x}$ -- простой цикл
системы  (\ref{np_2}), тогда существует (см. \cite{dem}, \S~20,
лемма~1) фундаментальная матрица $Y(t)$ системы (\ref{ls_2}) вида
\begin{equation}\label{Yt}
  Y(t)=\Phi(t)\left(\begin{array}{cc} {\rm e}^{\Lambda t} & 0_{{n-1}\times 1} \\
  0_{1\times{n-1}}
  & 1 \end{array}\right),
\end{equation}
где $\Phi$ -- $T$-периодическая матрица Флоке и $\Lambda$ --
постоянная $(n-1)\times(n-1)$-матрица с собственными значениями
отличными от $0.$ Для $k$-й компоненты вектора
$\zeta\in\mathbb{R}^n$ в дальнейшем используется обозначение
$\zeta^k$ или $[\zeta]^k.$ Для любого $\delta>0$ определим
множество $C_\delta\subset\mathbb{R}^n$ следующим образом
$$C_\delta=\left\{\zeta\in\mathbb{R}^n:\|P_{n-1}\zeta\|<\delta,\
\zeta^n\in\left(-\frac{\theta_2-\theta_1}{2},\frac{\theta_2-\theta_1}
{2}\right)\right\},$$ где
$$P_{n-1}\zeta=\left(\begin{array}{c}\zeta^1\\ \vdots \\ \zeta^{n-1} \\
0
 \end{array}\right).
$$
Пусть $\Gamma:B_\Delta(C_\delta)\to \Gamma(B_\Delta(C_\delta)),$
$\Delta>0,$ задается формулой
$$
 \Gamma(\zeta)= \frac{Y(\zeta^n+\overline{\theta})}{\|Y\|_{[0,T]}}
 P_{n-1}\zeta+\widetilde{x}(\zeta^n+\overline{\theta}),
$$
где $$ \overline{\theta}=\frac{\theta_1+\theta_2}{2}\quad\mbox{\rm
и}\quad \|Y\|_{[0,T]}=\max_{\theta\in[0,T]}\|Y(\theta)\|.$$

Справедливы следующие предварительные свойства. \vskip0.2truecm

\begin{lem}\label{scal_prod}
Для всех $\theta\in[0,T]$ и всех $\zeta\in\mathbb{R}^n$ имеет
место соотношение
$\left<Y(\theta)P_{n-1}\zeta,\widetilde{z}(\theta)\right>=0.$
Обратно, если $\left<\xi,\widetilde{z}(\theta)\right>=0,$ то
существует $\zeta\in\mathbb{R}^n$ такое, что
$\left<Y(\theta)P_{n-1}\zeta,\widetilde{z}(\theta)\right>=0.$
\end{lem}

\noindent {\bf Доказательство.} Пусть $\zeta\in\mathbb{R}^n.$
Рассмотрим следующий вектор
\[
  \widehat{\zeta}=\left(\begin{array}{cc}\left(I-{\rm e}^{\Lambda
  T}\right)^{-1} & 0 \\ 0 & 0 \end{array}\right)\zeta.
\]
По лемме Перрона  (см. \cite{perron} или \cite{dem}, Гл.~III,
\S~12) имеем
\[
  \left<Y(\theta+T)P_{n-1}\widehat{\zeta},\widetilde{z}(\theta)\right>=\left<
  Y(\theta)P_{n-1}\widehat{\zeta},\widetilde{z}(\theta)\right>\quad
  \mbox{\rm для\ всех\ }\theta\in[0,T].
\]
Следовательно,
\[
0=
\left<\left(Y(\theta)-Y(\theta+T)\right)P_{n-1}\widehat{\zeta},\widetilde{z}(\theta)\right>=
\]
\[
= \left<\Phi(\theta)\left(\begin{array}{cc} {\rm
e}^{\Lambda\theta}\left(I-{\rm e}^{\Lambda T}\right) & 0 \\ 0 & 0
\end{array}\right)P_{n-1}\widehat{\zeta},\widetilde{z}(\theta)\right>=
\]
\[
= \left<\Phi(\theta)\left(\begin{array}{cc} {\rm
e}^{\Lambda\theta} & 0 \\ 0 & 0
\end{array}\right)P_{n-1}\zeta,\widetilde{z}(\theta)\right>=\left<Y(\theta)P_{n-1}
\zeta,\widetilde{z}(\theta)\right>
\]
для любого $\theta\in[0,T],$ что является первым утверждением
леммы~\ref{scal_prod}. Чтобы доказать второе утверждение леммы,
положим
$$
  L_\xi=\left\{\xi\in\mathbb{R}^n:\left<\xi,\widetilde{z}(\theta)\right>=0\right\},\quad
  L_\zeta=\bigcup_{\zeta\in\mathbb{R}^n}Y(\theta)P_{n-1}\zeta.
$$
Заметим, что $L_{\xi}$ и $L_{\zeta}$ являются линейными
подпространствами пространства $\mathbb{R}^n,$ и ${\rm
dim}L_\xi=n-1.$ Так как $Y(\theta)P_{n-1}$ является линейным
невырожденным отображением, действующим из  $P_{n-1}\mathbb{R}^n$
в $Y(\theta)P_{n-1}\mathbb{R}^n,$ то  ${\rm dim}L_\zeta={\rm
dim}P_{n-1}\mathbb{R}^n=n-1.$ Но, согласно первому утверждению
леммы, $L_\xi\supset L_\zeta$ и значит, можно заключить, что
$L_\xi=L_\zeta.$

Лемма доказана.

\begin{lem}\label{lem_homeom}
Для любого $\Delta\in(0,\Delta_0]$ и любого
$\delta\in(0,\delta_0]$ отображение $\Gamma$ является
гомеоморфизмом $B_\Delta(C_\delta)$ на
$\Gamma(B_\Delta(C_\delta))$ при условии, что  $\Delta_0>0$ и
$\delta_0>0$ достаточно малы. Более того, множество
$\Gamma(B_\Delta(C_\delta))$ открыто в $\mathbb{R}^n,$ и
$\Gamma^{-1}$ непрерывно дифференцируемо на множестве
$\Gamma(B_\Delta(C_\delta)).$
\end{lem}

\noindent {\bf Доказательство.} Очевидно, что $\Gamma$ непрерывно.
Покажем, что  $\Gamma:B_\Delta(C_\delta)\to
\Gamma(B_\Delta(C_\delta))$ инъективно для достаточно малых
$\Delta>0$ и $\delta>0.$ Для этого предположим противное, тогда
существуют
$\{a_k\}_{k\in\mathbb{N}},\{b_k\}_{k\in\mathbb{N}}\subset\mathbb{R}^n,$
$a_k\not= b_k,$ $a_k\to a_0,$ $b_k\to b_0$ при $k\to\infty,$
\begin{equation}\label{Pform}
P_{n-1}a_0=P_{n-1}b_0=0,
\end{equation}
такие, что
\begin{eqnarray}\label{contr}
\frac{Y(a_k^n)}{\|Y\|_{[0,T]}}(P_{n-1}a_k)+\widetilde{x}(a_k^n)=\frac{Y(b_k^n)}
{\|Y\|_{[0,T]}}(P_{n-1}b_k)+\widetilde{x}(b_k^n).
\end{eqnarray}
Без ограничения общности можем считать, что либо $a^n_k=b^n_k$ для
любого  $k\in\mathbb{N},$ либо $a^n_k\not=b^n_k$ для любого
$k\in\mathbb{N}.$ Предположим, что  $a^n_k=b^n_k$ для любого
$k\in\mathbb{N},$ следовательно,
$$
  Y(a_k^n)(P_{n-1}a_k-P_{n-1}b_k)=0\quad\mbox{\rm для\ любого\
  }k\in\mathbb{N},
$$
и, таким образом,
$$ P_{n-1}a_k=P_{n-1}b_k\quad\mbox{\rm для\ любого\
  }k\in\mathbb{N},
$$
противореча тому свойству, что $a_k\not= b_k$ для любого
$k\in\mathbb{N}.$ Рассмотрим теперь случай когда $a^n_k\not=b^n_k$
для любого $k\in\mathbb{N}.$ Так как из (\ref{contr}) имеем, что
$\widetilde{x}(a_0^n)=\widetilde{x}(b_0^n)$ и  (по предположению
теоремы~\ref{MathNach1}) $|a_0^n-b_0^n|<\frac{T}{p},$ где
$\frac{T}{p}$ -- наименьший период цикла $\widetilde{x},$ то
$a_0^n=b_0^n=:\theta_0.$ Используя лемму~\ref{scal_prod}, из
(\ref{contr}) имеем
$$
  \left<\widetilde{x}(a_k^n)-\widetilde{x}(b_k^n),\widetilde{z}(a_k^n)\right>=\left<\frac{Y(b_k^n)}
  {\|Y\|_{[0,T]}}P_{n-1}b_k^n,\widetilde{z}(a_k^n)\right>=$$
$$\qquad=
  \left<\frac{Y(b_k^n)-Y(a_k^n)}{\|Y\|_{[0,T]}}P_{n-1}b_k^n,\widetilde{z}(a_k^n)\right>,
$$
или, после деления на $a_k^n-b_k^n,$
$$
  \left<\frac{\widetilde{x}(a_k^n)-\widetilde{x}(b_k^n)}{a_k^n-b_k^n},\widetilde{z}(a_k^n)\right>=-
  \frac{1}{\|Y\|_{[0,T]}}\left<\frac{Y(a_k^n)-Y(b_k^n)}{a_k^n-b_k^n}
  P_{n-1}b_k^n,\widetilde{z}(a_k^n)\right>.
$$
Переходя к пределу при  $k\to\infty$ в предыдущем равенстве  и
учитывая, что  $P_{n-1}b_k^n\to 0$ при $k\to\infty,$ получаем
$$
\left<\dot
{\widetilde{x}}(\theta_0),\widetilde{z}(\theta_0)\right>=0,
$$
что является противоречием (см. \cite{malb}, Гл.~III,
формула~12.9). Следовательно, существуют $\Delta_0>0$ и $\delta>0$
такие, что отображение $\Gamma:B_\Delta(C_\delta)\to
\Gamma(B_\Delta(C_\delta))$ инъективно для $\Delta\in(0,\Delta_0]$
и $\delta\in(0,\delta_0].$ Покажем, что $\Delta_0>$ и $\delta_0>0$
могут быть выбраны таким образом, что
\begin{equation}\label{need2}
\Gamma(B_\Delta(C_\delta)) \mbox{\rm\ \ открыто\ в\ }
\mathbb{R}^n\mbox{\rm\  для\ любых\ } \Delta\in(0,\Delta_0]
\mbox{\rm\ и\  }\delta\in(0,\delta_0].
\end{equation}
Заметим, что для любого $\zeta\in\mathbb{R}^n,$ удовлетворяющего
условию $P_{n-1}\zeta=0,$ имеем
$$\Gamma'(\zeta)=\frac{1}{\|Y\|_{[0,T]}}\Phi\left(\zeta^n+\overline
{\theta}\right)\left(\begin{array}{cc} {\rm
e}^{\Lambda(\zeta^n+\overline{\theta}}) & 0 \\ 0 & 0
\end{array}\right)+
\left(0\ \ldots\ 0\ \ \dot
{\widetilde{x}}(\zeta^n+\overline{\theta})\right)
$$
и, таким образом, для любого $\zeta\in\mathbb{R}^n$ такого, что
$P_{n-1}\zeta=0$ производная  $\Gamma'(\zeta)$ обратима.
Следовательно, без ограничения общности можем считать, что
$\Delta_0>0$ и $\delta_0>0$ достаточно малы так, что линейное
отображение $\Gamma'(\zeta)$ обратимо для любого $\zeta\in
B_\Delta(C_\delta)$ с $\Delta\in(0,\Delta_0]$ и
$\delta\in(0,\delta_0].$ По теореме об обратной функции (см.
\cite{rud}, теорема~9.17) имеем, что $\Gamma$ локально обратимо на
$B_\Delta(C_\delta)$ при $\Delta\in(0,\Delta_0]$ и
$\delta\in(0,\delta_0],$ то есть оно переводит любую достаточно
малую окрестность (в $\mathbb{R}^n$) элемента $\zeta$ в открытое
множество пространства $\mathbb{R}^n,$ что, в свою очередь,
означает (\ref{need2}). Более того, из теоремы об обратном
отображении следует, что  $\Gamma^{-1}$ непрерывно дифференцируемо
в $\Gamma(B_\Delta(C_\delta)).$

Лемма доказана.

\noindent {\bf Доказательство теоремы  \ref{MathNach1}.}  Прежде
всего заметим, что если $x$ -- решение уравнения  $x=Q_\varepsilon
x,$ то в силу леммы~\ref{zamena}  $\nu(t)=\Omega(0,t,x(t))$ --
решение уравнения $\nu=G_\varepsilon \nu,$  где
$G_\varepsilon:C([0,T],\mathbb{R}^n)\to C([0,T],\mathbb{R}^n)$
определяется как
$$
  \hskip-3cm(G_\varepsilon
  \nu)(t)=\Omega(T,0,\nu(T))+
$$
$$\qquad\hskip3cm+\varepsilon\int\limits_0^t
  \left(\Omega'_{(3)}(\tau,0,\nu(\tau))\right)^{-1}g(\tau,\Omega(\tau,0,\nu(\tau)),\varepsilon)d\tau.
$$
Более того, так как для любого открытого ограниченного множества
$V\subset\mathbb{R}^n$ гомеоморфизм $(S x)(t)=\Omega(0,t,x(t))$
отображает каждую окрестность множества $W_V$ в окрестность
множества
$$
  \widehat{W}_V=\left\{u\in C([0,T],\mathbb{R}^n):u(t)\in V,\ \mbox{\rm для\ любого\
  }t\in[0,T]\right\},
$$
то  (см. \cite{krazab}, теорема~26.4) имеем, что
$$
  d(I-Q_\varepsilon,W_{\Gamma(C_\delta)})=d(I-G_\varepsilon,
  \widehat{W}_{\Gamma(C_\delta)})
$$
в случае, если $d(I-G_\varepsilon,\widehat{W}_{\Gamma(C_\delta)})$
определен. Чтобы доказать, что
$d(I-G_\varepsilon,\widehat{W}_{\Gamma(C_\delta)})$ определен,
рассмотрим векторное поле
\[
A_\varepsilon(\xi)=x'_{(2)}\left(T-\varepsilon
\widetilde{M}\left(\left[\Gamma^{-1}(\xi)\right]^n\right),
\widetilde{x}\left(\left[\Gamma^{-1}
  (\xi)\right]^n+\overline{\theta}\right)
  \right)\circ \]
  \[\circ\left(\xi-\widetilde{x}\left(\left[\Gamma^{-1}(\xi)\right]^n+
  \overline{\theta}\right)\right)+
\]
\[
   +
  \widetilde{x}\left(\left[\Gamma^{-1}(\xi)\right]^n+\overline{\theta}-
  \varepsilon \widetilde{M}\left(   \left[\Gamma^{-1}(\xi)\right]^n
  \right)\right),\qquad \xi\in\Gamma(B_\Delta(C_\delta)),
\]
где $\Gamma, \Delta, \delta>0$ -- те, о которых говорится в
лемме~\ref{lem_homeom}, и $\widetilde{M}:\mathbb{R}\to\mathbb{R}$
определено как
$$
\widetilde{M}(t)=\left\{   \begin{array}{cl} \ \ |t|, & \mbox{\rm
\ \ если\ }
M_{\widetilde{x}}(\theta_1)<0\mbox{\rm\ и\ }M_{\widetilde{x}}(\theta_2)<0, \\
-|t|, & \mbox{\rm \ \ если\ }
M_{\widetilde{x}}(\theta_1)>0\mbox{\rm\ и\
}M_{\widetilde{x}}(\theta_2)>0, \\
-d_{\mathbb{R}}\left(M_{\widetilde{x}},\left(\theta_1,\theta_2\right)\right)
\cdot t, & \mbox{\rm \ \ в остальных случаях.}
\end{array} \right.
$$
Теперь докажем, что существует $\delta_0>0$ такое, что для любого
$\delta\in(0,\delta_0]$ и любого
$\varepsilon\in(0,\delta^{1+\alpha}]$ топологические степени
$d(I-G_\varepsilon,\widehat{W}_{\Gamma(C_{\delta})})$ и
$d_{\mathbb{R}^n}(I-A_\varepsilon,\Gamma(C_{\delta}))$ определены
и
\begin{equation}\label{co2}
  d(I-G_\varepsilon,\widehat{W}_{\Gamma(C_{\delta})})=d_{\mathbb{R}^n}
  (I-A_\varepsilon,\Gamma(C_{\delta})).
\end{equation}
Для этого определим вспомогательное векторное поле
$\widehat{A}_\varepsilon:C([0,T],\mathbb{R}^n)\to
C([0,T],\mathbb{R}^n),$ полагая $(\widehat{A}_\varepsilon
\nu)(t)={A}_\varepsilon (\nu(T))$ для любого  $t\in[0,T]$ и любого
$\nu\in C([0,T],\mathbb{R}^n).$ Так как
$\widehat{W}_{\Gamma(C_{\delta})}\cap\mathbb{R}^n=\Gamma(C_{\delta}),$
то по теореме о сужении (см. \cite{krazab}, теорема~27.1) степень
$d_{\mathbb{R}^n}(I-A_\varepsilon,\Gamma(C_{\delta}))$ определена,
если только определена степень
$d(I-\widehat{A}_\varepsilon,\widehat{W}_{\Gamma(C_{\delta})}),$
более того, $d_{\mathbb{R}^n}(I-A_\varepsilon,\Gamma(C_{\delta}))=
d(I-\widehat{A}_\varepsilon,\widehat{W}_{\Gamma(C_{\delta})}).$
Следовательно, желательно показать, что существует $\delta_0>0$
такое, что для любого $\delta\in(0,\delta_0]$ и любого
$\varepsilon\in(0,\delta^{1+\alpha}]$ топологические степени
$d(I-G_\varepsilon,\widehat{W}_{\Gamma(C_{\delta})})$ и
$d(I-\widehat{A}_\varepsilon,\widehat{W}_{\Gamma(C_{\delta})})$
определены и
\begin{equation}\label{co2p}
d(I-G_\varepsilon,\widehat{W}_{\Gamma(C_{\delta})})=
d(I-\widehat{A}_\varepsilon,\widehat{W}_{\Gamma(C_{\delta})}).
\end{equation}

\noindent Чтобы доказать (\ref{co2p}), определим
$F_\varepsilon:C([0,T],\mathbb{R}^n)\to C([0,T],\mathbb{R}^n)$ как
$$({F_\varepsilon}\nu)(t)=\int\limits_0^t
  \left(\Omega'_{(3)}(\tau,0,\nu(\tau))\right)^{-1}g(\tau,\Omega(\tau,0,\nu(\tau)),
  \varepsilon)d\tau, \quad t\in[0,T],$$
и рассмотрим линейную деформацию
\[
  D_\varepsilon\big(\lambda,\nu)(t)=\lambda\big(\nu(t)-\Omega(T,0,\nu(T))-\varepsilon
  ({F_\varepsilon}\nu)(t)\big)+\]
  \[\hskip5cm+(1-\lambda)\left(\nu(t)-
  \left(\widehat{A}_\varepsilon \nu\right)(t)\right),
\]
где $\lambda\in[0,1],\  u\in\partial
\widehat{W}_{\Gamma(C_\delta)},\
  \delta\in(0,\delta_0)$. Эквивалентно
\[
D_\varepsilon(\lambda,\nu)(t)=\lambda\left(\nu(t)-\Omega(T,0,\nu(T))\right)+
(1-\lambda)\nu(t)-
\]
\[
   -(1-\lambda)x'_{(2)}\left(T-\varepsilon \widetilde{M}\left(\left[\Gamma^{-1}(\nu(T))\right]^n\right),R_{\widetilde{x}}(\nu(T))\right)
   (\nu(T)-R_{\widetilde{x}}(\nu(T)))-
\]
\[
   - \lambda\varepsilon({F_\varepsilon}\nu)(t)-(1-\lambda)\,\widetilde{x}\left(\left[\Gamma^{-1}(\nu(T))
   \right]^n+
  \overline{\theta}-\varepsilon \widetilde{M}\left(   \left[\Gamma^{-1}(\nu(T))\right]^n
  \right)\right),
\]
где $\lambda\in[0,1],\ \nu\in\partial
\widehat{W}_{\Gamma(C_\delta)},\ \delta\in(0,\delta_0)$ и
$$R_{\widetilde{x}}(\xi)=\widetilde{x}\left(\left[\Gamma^{-1}(\xi)\right]^n+\overline{\theta}\right).$$
Покажем, что для всех достаточно малых  $\delta\in(0,\delta_0],$
$\varepsilon\in(0,\delta^{1+\alpha}]$ и каждых $\lambda\in[0,1],$
$\nu\in\partial \widehat{W}_{\Gamma(C_{\delta})}.$  Предположим
противное, значит существуют
$\{\delta_k\}_{k\in\mathbb{N}}\subset\mathbb{R}_+,$ $\delta_k\to
0$ при $k\to\infty,$ $\{\varepsilon_k\}_{k\in\mathbb{N}},$
$\varepsilon_k\in(0,\delta_k^{1+\alpha}),$
$\{\nu_k\}_{k\in\mathbb{N}},$ $\nu_k\in\partial
\widehat{W}_{\Gamma(C_{\delta_k})},$
 $\{\lambda_k\}_{k\in\mathbb{N}}\subset[0,1]$ такие, что
\begin{equation}
\begin{aligned}
0=&\lambda_k\left(\nu_k(t)-\Omega(T,0,\nu_k(T))\right)+(1-\lambda_k)\nu_k(t)-
\\ & -(1-\lambda_k)x'_{(2)}\left(T-\varepsilon_k
\widetilde{M}\left(\left[\Gamma^{-1}(\nu_k(T))\right]^n\right),R_{\widetilde{x}}(\nu_k(T))\right)\circ
\\ &
\circ(\nu_k(T)-R_{\widetilde{x}}(\nu_k(T)))-\lambda_k\varepsilon_k({F_{\varepsilon_k}}\nu_k)(t)-
\\ & -(1-\lambda_k)
\widetilde{x}\left(\left[\Gamma^{-1}(\nu_k(T))\right]^n+\overline{\theta}-
\varepsilon_k \widetilde{M}\left(
\left[\Gamma^{-1}(\nu_k(T))\right]^n \right)\right).\label{fir}
\end{aligned}
\end{equation}
Из (\ref{fir}) имеем
$$
\begin{aligned}
  \nu_k(t)=&\lambda_k \Omega(T,0,\nu_k(T))+\nonumber\\&+(1-\lambda_k)x'_{(2)}
  \left(T-\varepsilon_k
  \widetilde{M}\left(\left[\Gamma^{-1}(\nu_k(T))\right]^n\right),R_{\widetilde{x}}(\nu_k(T))\right)\circ
  \nonumber\\ &\circ(\nu_k(T)-R_{\widetilde{x}}(\nu_k(T)))+\lambda_k\varepsilon_k({F_{\varepsilon_k}}\nu_k)(t)+\nonumber\\
  & +(1-\lambda_k)
  \widetilde{x}\left(\left[\Gamma^{-1}(\nu_k(T))\right]^n+\overline{\theta}-
  \varepsilon_k \widetilde{M}\left(   \left[\Gamma^{-1}(\nu_k(T))\right]^n
  \right)\right)\nonumber
\end{aligned}
$$
и, следовательно,
\begin{equation}\label{cons}
\dot \nu_k(t)=\lambda_k\varepsilon_k
\left(\Omega'_{(3)}(t,0,\nu_k(t))\right)^{-1}g(t,\Omega(t,0,\nu_k(t)),\varepsilon_k).
\end{equation}
Из  (\ref{cons}) следует, что без ограничения общности можно
предположить существование  $\xi_0\in\mathbb{R}^n$ такого, что
\[
  \nu_k(t)\to \xi_0\mbox{\ \rm при\ }k\to\infty
\]
равномерно по отношению к  $t\in[0,T].$ Так как
$\nu_k(0)\in\Gamma(C_{\delta_k})\in
B_{\delta_k}(\widetilde{x}([\theta_1,\theta_2])),$ то $\xi_0\in
\widetilde{x}([\theta_1,\theta_2]).$ Теперь, чтобы получить
противоречие, возьмем $t=T$ и перепишем  (\ref{fir}) в виде
\[
     0=\lambda_k\left(\nu_k(T)-\Omega(T,0,\nu_k(T))\right)+(1-\lambda_k)\nu_k(T)-
\]
\[
    -(1-\lambda_k)\Omega'_{(3)}\left(T-\varepsilon_k
    \widetilde{M}\left(\left[\Gamma^{-1}(\nu_k(T))\right]^n\right),0,
   R_{\widetilde{x}}(\nu_k(T))\right)\circ
\]
\[
\circ(\nu_k(T)-R_{\widetilde{x}}(\nu_k(T)))-\lambda_k\varepsilon_k({F_{\varepsilon_k}}\nu_k)(T)-
\]
\[
    -(1-\lambda_k)\widetilde{x}
  \left(\left[\Gamma^{-1}(\nu_k(T))\right]^n+\overline{\theta}-
  \varepsilon_k \widetilde{M}\left(   \left[\Gamma^{-1}(\nu_k(T))\right]^n
  \right)\right)=
\]
\[
  =\lambda_k\left(\nu_k(T)-\Omega(T,0,\nu_k(T))\right)+
\]
\[+(1-\lambda_k)
 \left(I-\Omega'_{(3)}\left(T-\varepsilon_k
 \widetilde{M}\left(\left[\Gamma^{-1}(\nu_k(T))\right]^n\right),0,
 R_{\widetilde{x}}(\nu_k(T))\right)\right)\circ
\]
\[
  \circ (\nu_k(T)-R_{\widetilde{x}}(\nu_k(T)))-
  \lambda_k\varepsilon_k({F_{\varepsilon_k}}\nu_k)(T)+(1-\lambda_k)
  R_{\widetilde{x}}(\nu_k(T))-
\]
\[
   -(1-\lambda_k)\widetilde{x}\left(\left[\Gamma^{-1}(\nu_k(T))\right]^n+
  \overline{\theta}-\varepsilon_k \widetilde{M}\left(   \left[\Gamma^{-1}(\nu_k(T))\right]^n
  \right)\right)=
\]
\[
   = \lambda_k\left(\nu_k(T)-\Omega(T,0,\nu_k(T))\right)+(1-\lambda_k)\circ\nonumber
\]
\[
   \circ
  \left(I-\Omega'_{(3)}\left(T-\varepsilon_k
  \widetilde{M}\left(\left[\Gamma^{-1}(\nu_k(T))\right]^n\right),0,
  R_{\widetilde{x}}(\nu_k(T))\right)\right)\circ
\]
\[
   \circ(\nu_k(T)-R_{\widetilde{x}}(\nu_k(T)))-\lambda_k\varepsilon_k
   ({F_{\varepsilon_k}}\nu_k)(T)+
\]
\[
   +\varepsilon_k(1-\lambda_k)
  \dot {\widetilde{x}}\left(\left[\Gamma^{-1}(\nu_k(T))\right]^n+\overline{\theta}\right)
  \widetilde{M}\left(   \left[\Gamma^{-1}(\nu_k(T) )\right]^n
  \right)+o(\varepsilon_k).
\]
Таким образом, замечая, что
\[
\Omega(T,0,\xi)-\xi
=\Omega(T,0,\xi)-R_{\widetilde{x}}(\xi)+R_{\widetilde{x}}(\xi)-\xi=
\]
\[
=
\Omega(T,0,R_{\widetilde{x}}(\xi)+(\xi-R_{\widetilde{x}}(\xi)))-R_{\widetilde{x}}(\xi)+
R_{\widetilde{x}}(\xi)-\xi=
\]
\[
=\Omega'_{(3)}(T,0,R_{\widetilde{x}}(\xi))(\xi-R_{\widetilde{x}}(\xi))-(\xi-R_{\widetilde{x}}(\xi))+
o(\xi-R_{\widetilde{x}}(\xi))=
\]
\[
=\left(\Omega'_{(3)}(T,0,R_{\widetilde{x}}(\xi))-I\right)(\xi-R_{\widetilde{x}}(\xi))
+o(\xi-R_{\widetilde{x}}(\xi)),
\]
имеем
\begin{equation}\label{cc}
\begin{aligned}
  &
  \lambda_k\left(I-\Omega'_{(3)}(T,0,R_{\widetilde{x}}(\nu_k(T)))\right)(\nu_k(T)-R_{\widetilde{x}}(\nu_k(T)))-\\
& -\lambda_k o(\nu_k(T)-R_{\widetilde{x}}(\nu_k(T)))
+(1-\lambda_k)\circ
\\
  &\circ
  \left(I-\Omega'_{(3)}\left(T-\varepsilon_k
  \widetilde{M}\left(\left[\Gamma^{-1}(\nu_k(T))\right]^n\right),0,
  R_{\widetilde{x}}(\nu_k(T))\right)\right)\circ\\
  &\circ(\nu_k(T)-R_{\widetilde{x}}(\nu_k(T)))-\lambda_k\varepsilon_k
  ({F_{\varepsilon_k}}\nu_k)(T)+o(\varepsilon_k)+
\\
  & +\varepsilon_k(1-\lambda_k)
  \dot {\widetilde{x}}\left(\left[\Gamma^{-1}(\nu_k(T))\right]^n+\overline{\theta}\right)
  \widetilde{M}\left(   \left[\Gamma^{-1}(\nu_k(T) )\right]^n
  \right)=0.
\end{aligned}
\end{equation}
Можно считать, что последовательности
$\{\lambda_k\}_{k\in\mathbb{N}}$ и
$\left\{\dfrac{\nu_k(T)-R_{\widetilde{x}}(\nu_k(T))}{\|\nu_k(T)-R_{\widetilde{x}}(\nu_k(T))\|}
\right\}_{k\in\mathbb{N}}$ сходятся, пусть
$\lambda_0=\lim_{k\to\infty}\lambda_k$ и
$l_0=\lim_{k\to\infty}\dfrac{\nu_k(T)-R_{\widetilde{x}}(\nu_k(T))}{\|\nu_k(T)-R_{\widetilde{x}}(\nu_k(T))\|}.$
Так как $\nu_k\in\partial\widehat{W}_{\Gamma(C_{\delta_k})},$ то
существует $t_k\in[0,T]$ такое, что  $\nu_k(t_k)\in\partial
\Gamma(C_{\delta_k}).$ Положим $\zeta_k=\Gamma^{-1}(\nu_k(t_k)),$
без ограничения общности можем предполагать, что либо
\begin{equation}\label{case2}
\zeta_k^n+\overline{\theta}\in(\theta_1,\theta_2)\quad\mbox{\rm
для\ любого} \; k\in\mathbb{N},
\end{equation}
либо
\begin{equation}\label{case1}
  \zeta_k^n+\overline{\theta}\in\{\theta_1\}\cup\{\theta_2\}\quad\mbox{\rm для\ любого}
\;k\in\mathbb{N}.
\end{equation}

\vskip0.2truecm \noindent Покажем, что  случай (\ref{case2})
невозможен. По лемме~\ref{lem_homeom}, $\Gamma$ является
гомеоморфизмом множества $C_{\delta_k}$ на $\Gamma(C_{\delta_k})$
для достаточно малых  $\Delta>0,$ и
$\nu_k(t_k)\in\partial\Gamma(C_{\delta_k}),$ значит
\begin{equation}\label{tm}
\zeta_k=\Gamma^{-1}(\nu_k(t_k))\in\partial C_{\delta_k}.
\end{equation}
Следовательно, (\ref{case2}) и (\ref{tm}) влекут
\begin{equation}\label{imp}
  \|P_{n-1}\zeta_k\|=\delta_k\quad\mbox{\rm для\ всех\
  }k\in\mathbb{N}.
\end{equation}
Так как
\[
  \|P_{n-1}\zeta_k\|=\|Y^{-1}(\theta)Y(\theta)P_{n-1}\zeta_k\|\le\|Y^{-1}
  (\theta)\| \|Y(\theta)P_{n-1}\zeta_k\|,
\]
то существует  $c>0$ такое, что
\[
  \|Y(\theta)P_{n-1}\zeta_k\|\ge c\|P_{n-1}\zeta_k\|=c\delta_k
\]
для любого  $\theta\in [0,T]$, и, таким образом, имеем
$$
  \hskip-6cm\|\nu_k(t_k)-R_{\widetilde{x}}(\nu_k(t_k))\|=
$$
\begin{equation}\label{he1}
 \hskip2cm=\|\Gamma(\zeta_k)-\widetilde{x}(\zeta_k^n+
 \overline{\theta})\|=
 \left\|
 Y\left(\zeta_k^n+\overline{\theta}\right)P_{n-1}\zeta_k
 \right\|\ge c\delta_k
\end{equation}
для любого  $k\in\mathbb{N}$. С другой стороны из (\ref{cons})
заключаем, что существует  $c_1>0$ такое, что
\begin{equation}\label{he2}
  \|\nu_k(T)-\nu_k(t_k)\|\le c_1\varepsilon_k\quad \mbox{\rm для\ любого \
  }k\in\mathbb{N}.
\end{equation}
Наконец, из леммы~\ref{lem_homeom} следует, что функция
$\widetilde{x}\left(\left[\Gamma^{-1}(\cdot)\right]^n+\overline{\theta}\right)$
непрерывно дифференцируема и, учитывая (\ref{he2}), существует
$c_2>0$ такое, что
\begin{equation}
\begin{aligned}
  &\|R_{\widetilde{x}}(\nu_k(T))-R_{\widetilde{x}}(\nu_k(t_k))\|=\\
  &=\left\|
  \widetilde{x}\left(\left[\Gamma^{-1}(\nu_k(T))\right]^n+\overline{\theta}\right)-\widetilde{x}\left
  (\left[\Gamma^{-1}(\nu_k(t_k))\right]^n+\overline{\theta}\right)\right\|
  \le
\\ & \le  c_2\|\nu_k(T)-\nu_k(t_k)\|\le c_1 c_2\varepsilon_k\quad \mbox{\rm для\ всех \
  }k\in\mathbb{N}.\label{he3}
\end{aligned}
\end{equation}
Теперь возможно оценить $\|\nu_k(T)-R_{\widetilde{x}}(\nu_k(T))\|$
снизу. Имеем
\begin{equation}
\begin{aligned}
& \|\nu_k(T)-R_{\widetilde{x}}(\nu_k(T))\|=
\\&=\|\nu_k(t_k)-R_{\widetilde{x}}(\nu_k(t_k))+\\
&+\nu_k(T)-
\nu_k(t_k)-(R_{\widetilde{x}}(\nu_k(T))-R_{\widetilde{x}}
(\nu_k(t_k)))\|\ge
\\
& \ge\left|\|\nu_k(t_k)-R_{\widetilde{x}}(\nu_k(t_k))\|-\right.
\\
& \left.\|\nu_k(T)-\nu_k(t_k)-(R_{\widetilde{x}}(\nu_k(T))-
R_{\widetilde{x}}(\nu_k(t_k)))\|\right|\label{he5}.
\end{aligned}
\end{equation}
Так как  $\varepsilon_k\in(0,\delta_k^{1+\alpha}),$ то существует
 $k_0\in\mathbb{N}$ такое, что  $c_1\varepsilon_k+c_1
c_2\varepsilon_k<c\delta_k$ для всех $k\ge k_0.$ Следовательно, из
(\ref{he2}) и (\ref{he3}) имеем
\begin{equation}\label{he4}
\|\nu_k(T)-\nu_k(t_k)-(R_{\widetilde{x}}(\nu_k(T))-R_{\widetilde{x}}(\nu_k(t_k)))\|\le
 c_1\varepsilon_k+c_1
c_2\varepsilon_k< c\delta_k,
\end{equation}
для всех $k\ge k_0$. Используя (\ref{he1}) и (\ref{he4}), можем
переписать (\ref{he5}) как
\begin{equation}
\begin{aligned}
& \|\nu_k(T)-R_{\widetilde{x}}(\nu_k(T))\|\ge
\|\nu_k(t_k)-R_{\widetilde{x}}(\nu_k(t_k))\|-
\\& -\|\nu_k(T)-\nu_k(t_k)-(R_{\widetilde{x}}(\nu_k(T))-
R_{\widetilde{x}}(\nu_k(t_k)))\|\label{he6}
\end{aligned}
\end{equation}
и, таким образом,
\[
  \|\nu_k(T)-R_{\widetilde{x}}(\nu_k(T))\|\ge c\delta_k - c_1\varepsilon_k-c_1
c_2\varepsilon_k\quad\mbox{\rm для\ любого\ }k\ge k_0.
\]
Используя это неравенство для любого  $k\ge k_0$ получаем
\begin{equation}\label{pro}
\begin{aligned}
  & \frac{\varepsilon_k}{\|\nu_k(T)-R_{\widetilde{x}}(\nu_k(T))\|}\le
  \frac{\varepsilon_k}{c\delta_k - c_1\varepsilon_k-c_1
  c_2\varepsilon_k}\le
\\ \le & \frac{\delta_k^{1+\alpha}}{c\delta_k- c_1
\delta_k^{1+\alpha}-c_1 c_2
\delta_k^{1+\alpha}}=\frac{\delta_k^\alpha}{c-c_1
\delta_k^\alpha-c_1 c_2  \delta_k^\alpha}.
\end{aligned}
\end{equation}
Используя  (\ref{pro}) и переходя к пределу при $k\to\infty$ в
(\ref{cc}) деленном на $\|\nu_k(T)-R_{\widetilde{x}}(\nu_k(T))\|,$
получаем
\begin{equation}\label{res}
  \left(I-\Omega'_{(3)}(T,0,\widetilde{x}(\zeta_0^n+\overline{\theta}))\right)l_0=0.
\end{equation}
Чтобы доказать, что  (\ref{res}) не верно, установим, что
\begin{equation}\label{lus}
  \left<(I-\Omega'_{(3)}(T,0,\xi_0))l_0,\widetilde{z}\left(\left[\Gamma^{-1}(\xi_0)\right]^n+
  \overline{\theta}\right)\right>=0.
\end{equation}
Действительно,
\[
  \left<\frac{\nu_k(T)-R_{\widetilde{x}}(\nu_k(T))}{\|\nu_k(T)-R_{\widetilde{x}}(\nu_k(T))\|},\widetilde{z}
  \left(\left[\Gamma^{-1}(\nu_k(T))\right]^n+\overline{\theta}\right)\right>=
\]
\[= \frac{1}{\|\nu_k(T)-R_{\widetilde{x}}(\nu_k(T))\|}
   \left<\Gamma(\Gamma^{-1}(\nu_k(T)))-\widetilde{x}\left(\left[\Gamma^{-1}(\nu_k(T))\right]^n+
   \overline{\theta}\right),\right.
\]
\[
   \hskip3cm\left.\widetilde{z}\left(\left[\Gamma^{-1}(\nu_k(T))\right]^n+
   \overline{\theta}\right)\right>=
\]
\[
  =\frac{1}{\|\nu_k(T)-R_{\widetilde{x}}(\nu_k(T))\|}\left<Y\left(\left[\Gamma^{-1}
  (\nu_k(T))\right]^n+\overline{\theta}\right)P_{n-1}\Gamma^{-1}(\nu_k(T)),\right.
\]
\[\hskip3cm\left.\widetilde{z}
  \left(\left[\Gamma^{-1}(\nu_k(T))\right]^n+\overline{\theta}\right)\right>,
\]
и, таким образом, пользуясь леммой~\ref{scal_prod}, можем
заключить, что
\begin{equation}\label{willu}
  \left<\frac{\nu_k(T)-R_{\widetilde{x}}(\nu_k(T))}{\|\nu_k(T)-R_{\widetilde{x}}(\nu_k(T))\|},\widetilde{z}
  \left(\left[\Gamma^{-1}(\nu_k(T))\right]^n+\overline{\theta}\right)\right>=
  0,\qquad k\in\mathbb{N}.
\end{equation}
По определению вектора $l_0$ из (\ref{willu}) получаем
\begin{equation}\label{willu1}
  \left<l_0,\widetilde{z}\left(\zeta^n_0+\overline{\theta}\right)\right>=0.
\end{equation}
Но  $\|l_0\|=1,$ и на основании из леммы~\ref{scal_prod} имеем,
что существует $l_*\not=0$ такое, что
\begin{equation}\label{lstar}
l_0=Y\left(\zeta^n_0+\overline{\theta}\right)P_{n-1}l_*\quad\mbox{\rm
и}\quad P_{n-1}l_*=l_*
\end{equation}
и, замечая, что  (см., например, \cite{kraope}, теорема~2.1),
\begin{equation}\label{xylink}
  \Omega'_{(3)}(t,0,\widetilde{x}(\tau))=Y(t+\tau)Y^{-1}(\tau)\quad \mbox{\rm для\ всех\
  }t,\tau\in\mathbb{R},
\end{equation}
имеем
\begin{equation}\label{abc}
\begin{aligned}
   &
  \left(I-\Omega'_{(3)}(T,0,\widetilde{x}(\zeta_0^n+\overline{\theta}))\right)l_0= \left(I-Y\left(T+\zeta_0^n+\overline{\theta}\right)
   Y^{-1}\left(\zeta_0^n+\overline{\theta}\right)\right)l_0=
\\
   & =\left(Y\left(\zeta_0^n+\overline{\theta}\right)-Y\left(T+\zeta_0^n+\overline{\theta}\right)
   \right)P_{n-1}l_*=
\\
   &=\Phi\left(\zeta_0^n+\overline{\theta}\right)\left(\left(
    \begin{array}{ll} {\rm
    e}^{\Lambda(\zeta_0^n+\overline{\theta})}& 0 \\
    0 & 1\end{array}\right)-\left(\begin{array}{ll} {\rm
    e}^{\Lambda(T+\zeta_0^n+\overline{\theta})}& 0 \\
    0 & 1\end{array}\right)\right)P_{n-1}l_*=\\
\\ &=\Phi\left(\zeta_0^n+\overline{\theta}\right)
    \left(
    \begin{array}{ll} {\rm
    e}^{\Lambda(\zeta_0^n+\overline{\theta})} (I-{\rm e}^{\Lambda T})& 0 \\
    0 & 0\end{array}\right)P_{n-1}l_*
\end{aligned}
\end{equation}
в противоречии с (\ref{res}).

\vskip0.2truecm Покажем теперь, что случай (\ref{case1}) также
приводит к противоречию. Если, переходя при необходимости к
подпоследовательности,
$\dfrac{\varepsilon_k}{\|\nu_k(T)-R_{\widetilde{x}}(\nu_k(T))\|}
\to 0,$ то возможно действовать как и прежде для получения ложного
факта (\ref{res}). Рассмотрим теперь случай, когда
$\dfrac{\varepsilon_k}{\|\nu_k(T)-R_{\widetilde{x}}(\nu_k(T))\|}
\to l$, с $l>0$ или $l=+\infty.$ Из (\ref{cc}) заключаем, что
\begin{equation}\label{caca}
\begin{aligned}
& \frac{\varepsilon_k}{\|\nu_k(T)-R_{\widetilde{x}}(\nu_k(T))\|
}\left<\Xi_k(\widetilde{x})(T),\widetilde{z}\left(\left[\Gamma^{-1}(\nu_k(T))\right]^n+\overline{\theta}\right)\right>=
\\
&
=\left<\Upsilon_k(\widetilde{x})(T),\widetilde{z}\left(\left[\Gamma^{-1}(\nu_k(T))\right]^n+
\overline{\theta}\right)\right>,
\end{aligned}
\end{equation}
 где
\[
 \Xi_k(\widetilde{x})(T):=\lambda_k
({F_{\varepsilon_k}}\nu_k)(T)-(1-\lambda_k)\dot {\widetilde{x}}
  \left(\left[\Gamma^{-1}(\nu_k(T))\right]^n+\overline{\theta}\right)\cdot
\]
\[\hskip3cm\cdot
  \widetilde{M}\left(   \left[\Gamma^{-1}(\nu_k(T) )\right]^n
  \right)+\frac{o(\varepsilon_k)}{\varepsilon_k},
\]
\[
\Upsilon_k(\widetilde{x})(T):=
\lambda_k\left(I-\Omega'_{(3)}(T,0,R_{\widetilde{x}}(\nu_k(T)))\right)\frac{\nu_k(T)-R_{\widetilde{x}}(\nu_k(T)))}
{\|\nu_k(T)-R_{\widetilde{x}}(\nu_k(T))\|} -
\]
\[
\ \ \ -\lambda_k
\frac{o(\nu_k(T)-R_{\widetilde{x}}(\nu_k(T)))}{\|\nu_k(T)-R_{\widetilde{x}}(\nu_k(T))\|}+
\]
\[\ \ \  +(1-\lambda_k)\left(I-\Omega'_{(3)}\left(T-\varepsilon_k
\widetilde{M}\left(\left[\Gamma^{-1}(\nu_k(T))\right]^n\right),0,R_{\widetilde{x}}(\nu_k(T))\right)\right)\circ
\]
\[\circ\frac{\nu_k(T)-R_{\widetilde{x}}(\nu_k(T))}
{\|\nu_k(T)-R_{\widetilde{x}}(\nu_k(T))\|}.
\]
Используя представление (\ref{lstar}), формулу (\ref{abc}) и
лемму~\ref{scal_prod}, получаем
\[
\left<\left(I-\Omega'_{(3)}\left(T,0,\widetilde{x}\left(\zeta_0^n+\overline{\theta}\right)\right)\right)l_0,
    \widetilde{z}\left(\zeta_0^n+\overline{\theta}\right)\right>=
\]
\[=
    \left<Y\left(\zeta_0^n+\overline{\theta}\right)\left(I-{\rm e}^{\Lambda T}\right)
    P_{n-1}l_*,
    \widetilde{z}\left(\zeta_0^n+\overline{\theta}\right)\right>=
\]
\[
    =\left<Y\left(\zeta_0^n+\overline{\theta}\right)
    P_{n-1}\left(I-{\rm e}^{\Lambda T}\right)l_*,
    \widetilde{z}\left(\zeta_0^n+\overline{\theta}\right)\right>=0.
\]
Следовательно,
$$
\left<\Upsilon_k(\widetilde{x})(T),\widetilde{z}\left
  (\left[\Gamma^{-1}(\nu_k(T))\right]^n+\overline{\theta}\right)\right>\
  \to 0\quad\mbox{\rm при\ }k\to\infty,
$$
и из  (\ref{caca}) заключаем, что
\[
\left<\Xi_k(\widetilde{x})(T),\widetilde{z}\left(\left[\Gamma^{-1}(\nu_k(T))\right]^n+
\overline{\theta}\right)\right> \to 0\quad\mbox{\rm при\
}k\to\infty
\]
откуда, в свою очередь, следует
\begin{equation}\label{fa}
  \hskip-0.05cm\left<
\lambda_0\widehat{
F}\left(\widetilde{x}\left(\zeta^n_0+\overline{\theta}\right)\right)-(1-\lambda_0)\dot{
\widetilde{x}}\left(\zeta^n_0+\overline{\theta}\right)
  \widetilde{M}\left(   \zeta^n_0
  \right),\widetilde{z}\left(\zeta^n_0+\overline{\theta}\right)\right>=0,
\end{equation}
где
\[
\widehat {F}(\xi)=\int\limits_0^T
  \left(\Omega'_{(3)}(\tau,0,\xi)\right)^{-1}g(\tau,\Omega(\tau,0,\xi),0)d\tau.
\]
По лемме Перрона (см. \cite{perron} или \cite{dem}, Гл.~III,
\S~12)
\[
\left< \dot
{\widetilde{x}}\left(\zeta^n_0+\overline{\theta}\right)
  \widetilde{M}\left(   \zeta^n_0
  \right),\widetilde{z}\left(\zeta^n_0+\overline{\theta}\right)\right>=
\left<\dot
{\widetilde{x}}(0),\widetilde{z}(0)\right>\widetilde{M}\left(
\zeta^n_0
  \right)
\]
и, таким образом,  (\ref{fa}) может быть переписано как
$$
  \lambda_0\,{\rm sign}\left<\dot {\widetilde{x}}(0),\widetilde{z}(0)\right>\left<
\widehat
{F}\left(\widetilde{x}(\zeta^n_0+\overline{\theta})\right),
\widetilde{z}(\zeta^n_0+\overline{\theta})\right>-
$$
\begin{equation}\label{cbr}
-(1-\lambda_0)\left|\left<\dot
{\widetilde{x}}(0),\widetilde{z}(0)\right>\right|M\left( \zeta^n_0
  \right) =0.
\end{equation}
Покажем, что
\begin{equation}\label{gf}
{\rm sign}\left<\dot
{\widetilde{x}}(0),\widetilde{z}(0)\right>\left< \widehat{
F}\left(\widetilde{x}\left(\theta\right)\right),
\widetilde{z}\left(\theta\right)\right>=M_{\widetilde{x}}(\theta),\qquad\theta\in[0,T].
\end{equation}
Обозначим через $Z(t)$ и $Z_0(t)$ фундаментальные матрицы
сопряженной системы (\ref{ss_2}) такие, что $Z(0)=I$ и
$Z_0(t)=(Z_{n-1}(t),\widetilde{z}(t)),$ где $Z_{n-1}(t)$ является
$n\times(n-1)$-матрицей, чьи столбцы являются не
$T$-периодическими линейно-независимыми  собственными функциями
системы (\ref{ss_2}). Так как
\[
  \left(x'_{(3)}(\tau,0,\widetilde{x}(\theta))\right)^{-1}=Y(\theta)Y^{-1}(\tau+\theta)=
\]
\[=
  \left(Z^{-1}(\theta)\right)^* Z^*(\tau+\theta)=\left(Z^{-1}_0(\theta)
  \right)^* Z^*_0(\tau+\theta),
\]
 (см., напр., \cite{dem}, Гл.~III, \S~12), и
$\widetilde{z}(\theta)=\left(Z_{n-1}(\theta),\widetilde{z}(\theta)\right)\left(\begin{array}{c} 0\\ \vdots\\ 0 \\
1\end{array}\right)$ то имеем
\[
\left< \widehat {F}\left(\widetilde{x}\left(\theta\right)\right),
\widetilde{z}\left(\theta\right)\right>=\left<\left(Z^{-1}_0(\theta)\right)
^*\int\limits_0^T
Z^*_0(\tau+\theta)\,g(\tau,\widetilde{x}(\tau+\theta),0)d\tau,
\widetilde{z}\left(\theta\right)\right>=
\]
\[
=\left<\int\limits_\theta^{T+\theta} \left(\begin{array}{c}
Z_{n-1}^*(\tau)\\ \widetilde{z}(\tau)\end{array}\right)
g(\tau-\theta,\widetilde{x}(\tau),0)d\tau,
\left(\begin{array}{c} 0\\ \vdots\\ 0 \\
1\end{array}\right)\right>=
\]
\[
=\int\limits_\theta^{T+\theta}
 \left<\widetilde{z}(\tau)
g(\tau-\theta,\widetilde{x}(\tau),0)\right>d\tau=M_{\widetilde{x}}(\theta)
\]
и, таким образом,  (\ref{gf}) выполнено. Учитывая (\ref{gf}),
можем  переписать   (\ref{cbr}) как
\[
  \lambda_0 M_{\widetilde{x}}\left(\zeta^n_0+\overline{\theta}\right)
-(1-\lambda_0)\left|\left<\dot
{\widetilde{x}}(0),\widetilde{z}(0)\right>\right|\widetilde{M}\left(
\zeta^n_0
  \right) =0,
\]
где либо $\zeta^n_0+\overline{\theta}=\theta_1,$ либо
$\zeta^n_0+\overline{\theta}=\theta_2.$ Последнее, в свою очередь,
может быть переписано как
\begin{equation}\label{rev}
  \lambda_0 M_{\widetilde{x}}\left(\theta_i\right)
-(1-\lambda_0)\left|\left<\dot
{\widetilde{x}}(0),\widetilde{z}(0)\right>\right|\widetilde{M}\left(
(-1)^i|\zeta^n_0|
  \right) =0,
\end{equation}
где либо $i=1,$ либо $i=2.$ Если
$d_{\mathbb{R}}(M_{\widetilde{x}},(\theta_1,\theta_2))=0,$ то
(см., напр., \cite{krazab}, \S 3.2 по поводу определения
топологической степени в $\mathbb{R}$) для любого $i=1,2$ и любого
$a\ge 0$ имеем
$$
\widetilde{M}\left((-1)^i a\right)=-a{\rm
sign}(M_{\widetilde{x}}(\theta_1))=-a{\rm
sign}(M_{\widetilde{x}}(\theta_2)),
$$ и, таким образом, если
$d_{\mathbb{R}}(M_{\widetilde{x}},(\theta_1,\theta_2))=0,$ то
(\ref{rev}) может быть переписано как
\begin{equation}\label{rew1}
  \lambda_0 M_{\widetilde{x}}\left(\theta_i\right)
+(1-\lambda_0)\left|\left<\dot
{\widetilde{x}}(0),\widetilde{z}(0)\right>\right|\cdot|\zeta_0^n|\cdot{\rm
sign} \left(M_{\widetilde{x}}(\theta_i)\right) =0,
\end{equation}
где либо $i=1,$ либо $i=2.$ Если
$d_{\mathbb{R}}(M_{\widetilde{x}},(\theta_1,\theta_2))\not=0,$ то
для любых $i=1,2$ и любых $a\ge 0$ имеем
$$
\widetilde{M} \left((-1)^i
a\right)=d_{\mathbb{R}}(M_{\widetilde{x}},(\theta_1,\theta_2))\cdot(-1)^{i+1}
a=
$$
$$
=(-1)^{i+1} {\rm sign}(M_{\widetilde{x}}(\theta_i))(-1)^i a=-a{\rm
sign}(M_{\widetilde{x}}(\theta_i))$$ и, таким образом, (\ref{rev})
может быть вновь переписано как  (\ref{rew1}). Но (\ref{rew1})
противоречит либо предположению $M_{\widetilde{x}}(\theta_1)\not=
0$ (в случае $i=1$), либо предположению
$M_{\widetilde{x}}(\theta_2)\not= 0$ (в случае $i=2$).

\vskip0.2truecm \noindent Следовательно, ни (\ref{case1}), ни
(\ref{case2}) не могут осуществиться и, таким образом, существует
$\delta_0>0$ такое, что для любых $\delta\in(0,\delta_0]$ и любых
 $\varepsilon\in(0,\delta^{1+\alpha}]$ имеем, что
$D_\varepsilon(\lambda,\nu)\not=0$ для всех $\lambda\in[0,1]$ и
всех $\nu\in\partial\widehat{W}_{\Gamma(C_{\delta})}.$ Значит, для
любых $\delta\in(0,\delta_0]$ и
$\varepsilon\in(0,\delta^{1+\alpha}]$ топологические степени
$d(I-G_\varepsilon,\widehat{W}_{\Gamma(C_{\delta})})$ и
$d(I-\widehat{A}_\varepsilon,\widehat{W}_{\Gamma(C_{\delta})})$
определены, и (\ref{co2p}) выполнено. Как уже замечено,
(\ref{co2p}) влечет (\ref{co2}), значит, чтобы закончить
доказательство формулы (\ref{STAR}), остается показать, что
$d(I-A_\varepsilon,\Gamma(C_\delta))=(-1)^{\beta(\widetilde{x})}
d_{\mathbb{R}}\left(M_{\widetilde{x}},\left(\theta_1,\theta_2\right)\right)
$ для всех $\delta\in(0,\delta_0]$ и
$\varepsilon\in(0,\delta^{1+\alpha}].$ Пусть
$\delta\in(0,\delta_0]$ и $\varepsilon\in(0,\delta^{1+\alpha}].$
Так как $\Gamma$ является гомеоморфизмом множества
$B_\Delta(C_{\delta})$ на $\Gamma(B_\Delta(C_{\delta})),$ то (см.,
напр., \cite{krazab}, теорема~26.4) получаем
\[
  d_{\mathbb{R}^n}(I-A_\varepsilon,\Gamma(C_\delta))=d_{\mathbb{R}^n}
  (I-\Gamma^{-1}A_\varepsilon\Gamma,C_\delta).
\]
Пусть $\zeta\in C_\delta.$ Учитывая, что (\ref{xylink}) и действуя
как в (\ref{abc}), получаем
\[
  \zeta-(\Gamma^{-1}A_\varepsilon\Gamma)(\zeta)=\zeta-
  (\Gamma^{-1}A_\varepsilon)\left(\frac{Y(\zeta^n+\overline{\theta})}
  {\|Y\|_{M_T}}
  P_{n-1}\zeta+\widetilde{x}(\zeta^n+\overline{\theta})\right)=
\]
\[
  =\zeta-\Gamma^{-1}\left(\Omega'_{(3)}\left(T-\varepsilon \widetilde{M}(\zeta^n),0,\widetilde{x}\left(\zeta^n+
  \overline{\theta}\right)\right)\frac{Y(\zeta^n+\overline{\theta})}
  {\|Y\|_{M_T}}P_{n-1}\zeta+\right.
\]
\[\left.+
  \widetilde{x}(\zeta^n+\overline{\theta}-\varepsilon
    \widetilde{M}(\zeta^n))\right)=
\]
\[=\zeta-\Gamma^{-1}\Bigg(\frac{Y\left(\zeta^n+\overline{\theta}-\varepsilon
\widetilde{M}(\zeta^n)\right)}
  {\|Y\|_{M_T}}P_{n-1}\left(\begin{array}{cc}
  {\rm e}^{\Lambda T} & 0 \\
  0 &   0\end{array}\right)\zeta +
\]
\[
+\widetilde{x}\left(\zeta^n+\overline{\theta}-
  \varepsilon \widetilde{M}(\zeta^n)\right)\Bigg)=
\]
\[
  =\zeta-\left(\begin{array}{c}
  {\rm e}^{\Lambda T}\left(
  {\zeta|}_{\mathbb{R}^{n-1}}\right) \\
    \zeta^n-\varepsilon \widetilde{M}(\zeta^n)\end{array}\right)
\] и, таким образом,
\[
  d_{\mathbb{R}^n}(I-\Gamma^{-1}A_\varepsilon\Gamma,C_\delta)=
  d_{\mathbb{R}^n}\left(\left(I-{\rm e}^{\Lambda T}\right)\times
  \varepsilon \widetilde{M},C_\delta\right),
\]
где  $\left(I-{\rm e}^{\Lambda T}\right)\times
  \varepsilon \widetilde{M}=\left(I-{\rm e}^{\Lambda T},
  \varepsilon \widetilde{M}\right).$ По свойству топологической степени произведения векторных полей
   (см., напр., \cite{krazab}, теорема~7.4) имеем
\[
  d_{\mathbb{R}^n}\left(\left(I-{\rm e}^{\Lambda T}\right)\times
   \varepsilon \widetilde{M},C_\delta\right)=
\]
\[=
  d_{\mathbb{R}^n}\left(I-{\rm e}^{\Lambda
  T},B_{\delta}(0)\right)\cdot
  d_{\mathbb{R}}\left(\varepsilon  \widetilde{M},\left(-\frac{\theta_2-\theta_1}{2},
  \frac{\theta_2-\theta_1}{2}\right)\right),
\]
где  $d_{\mathbb{R}^n}\left(I-{\rm e}^{\Lambda
T},B_{\delta}(0)\right)=(-1)^{\beta(\widetilde{x})}$ согласно
(\cite{krazab}, теорема~6.1) и
\[ d_{\mathbb{R}}\left(\varepsilon
  \widetilde{M},\left(-\frac{\theta_2-\theta_1}{2},\frac{\theta_2-
  \theta_1}{2}\right)\right)=
-d_{\mathbb{R}}\left(M_{\widetilde{x}},\left(\theta_1,\theta_2\right)\right)
\]
прямым подсчетом. Таким образом, имеем
\[
  d_{\mathbb{R}^n}(I-\Gamma^{-1}A_\varepsilon\Gamma,C_\delta)=
  -(-1)^{\beta(\widetilde{x})}d_{\mathbb{R}}\left(M_{\widetilde{x}},\left(\theta_1,
  \theta_2\right)\right).
\]
Подводя итог, заключаем, что существует  $\delta_0>0$ такое, что
для любого $\delta\in(0,\delta_0]$ и любого
$\varepsilon\in(0,\delta^{1+\alpha}]$ топологическая степень
$d(I-Q_\varepsilon,W_{\Gamma(C_\delta)})$ определена и может быть
подсчитана по формуле
\[
d(I-Q_\varepsilon,W_{\Gamma(C_\delta)})=-(-1)^{\beta(\widetilde{x})}
d_{\mathbb{R}}\left(M_{\widetilde{x}},\left(\theta_1,\theta_2\right)\right).
\]
Чтобы завершить доказательство, остается показать, что
$V_\delta:=\Gamma(C_\delta)$ удовлетворяет свойствам  1) и 2). Для
этой цели,  положим $\xi\in\Gamma(C_\delta),$ значит
\[
  \xi=\frac{Y(\zeta^n+\overline{\theta})}{\|Y\|_{M_T}}
  P_{n-1}\zeta+\widetilde{x}(\zeta^n+\overline{\theta}).
\]
для некоторого $\zeta\in\mathbb{R}^n,$ удовлетворяющего
$\|P_{n-1}\zeta\|\le\delta$ и
$\left[\Gamma^{-1}(\xi)\right]^n+\overline{\theta}\in[\theta_1,\theta_2].$
Следовательно,
\[
\left\|\xi-\widetilde{x}\left(\left[\Gamma^{-1}(\xi)\right]^n+\overline{\theta}
\right)\right\|=
\left\|\frac{Y(\zeta^n+\overline{\theta})}{\|Y\|_{M_T}}P_{n-1}\zeta\right\|
\le\|P_{n-1}\zeta\|\le\delta
\]
и, значит, свойство 1) выполнено. По определению множества
$C_\delta$ имеем, что для любого $\delta\in(0,\delta_0)$ точки
$\left(0,...,0,-\dfrac{\theta_2-\theta_1}{2}\right)$ и
$\left(0,...,0,\dfrac{\theta_2-\theta_1}{2}\right)$ пренадлежат
границе множества $C_\delta.$ Следовательно, точки
$\widetilde{x}(\theta_1)$ и $\widetilde{x}(\theta_2)$ пренадлежат
границе множества $\Gamma(C_\delta).$ С другой стороны, если
$\xi=\widetilde{x}(\theta),$ где $\theta\in(\theta_1,\theta_2),$
то
\begin{equation}\label{starr}
\Gamma^{-1}(\xi)=\left(0,...,0,\theta-\overline{\theta}\right)\subset
C_\delta.
\end{equation}
Следовательно, $\xi\in\Gamma(C_\delta)$ и свойство  2) также
удовлетворено.

Доказательство теоремы \ref{MathNach1} завершено.

Далее нам необходима следующая лемма, принадлежащая И.~Г.~Малкину
(\cite{mal}, формула~3.13 или \cite{malb}, теорема с.~387).

\begin{lem}\label{malkin_lem} Если
$$
M_x(0)\not=0\quad\mbox{\it для\ некоторого\ }x\in\mathfrak{S}_W,
$$
то
$$
Q_\varepsilon x\not=x \mbox{ для \ любого\ достаточно\ малого\
}\varepsilon>0.
$$
\end{lem}

Для каждого $x\in\mathfrak{S}_W$ положим
\[
\Theta_W(x)=\left\{\theta_0\in(0,T):S_{\theta_0}\,x\in\partial W,\
S_\theta\,
  x\in W \mbox{\rm\ для\ любого\ }\theta\in(0,\theta_0)\right\},
\]
\[
(S_\theta\,x)(t) = x(t+\theta) \quad{\rm и}
\]

В дальнейшем постоянные функции пространства
$C([0,T],\mathbb{R}^n)$ отождествляются с соответствующими
элементами пространства $\mathbb{R}^n.$

\begin{thm}\label{MathNach2}
Предположим, что множество  $$\mathfrak{S}_W=\left\{x\in\partial
W:Q_0 x=x\right\}$$ конечно и содержит только простые
$T$-периодические циклы системы (\ref{np_2}). Предположим, что
$M_x(0)\not=0$ для всех $x\in\mathfrak{S}_W.$ Тогда для каждого
достаточно малого $\varepsilon>0$ топологическая степень
$d(I-Q_\varepsilon,W)$ определена, и справедлива следующая формула
$$
d(I-Q_\varepsilon,W)=(-1)^n
d_{\mathbb{R}^n}(f,W\cap\mathbb{R}^n)-
$$
\begin{equation}\label{kmnform}
-\sum_{x\in\mathfrak{S}_W:~ \Theta_W(x)\not=
\emptyset}(-1)^{\beta(x)}
d_{\mathbb{R}}\left(M_x,\left(0,\min\{\Theta_W(x)\}\right)\right).
\end{equation}
\end{thm}

\noindent {\bf Доказательство.} Для каждого  $x\in\mathfrak{S}_W,$
удовлетворяющего  $\Theta_W(x)\not=\emptyset,$ обозначим через
$\delta_0(x)$ и $\{V_\delta(x)\}_{\delta\in(0,\delta_0(x))}$ те
числа и множества, о которых говорится в теореме~\ref{MathNach1},
где $\widetilde{x}:=x,$ $\theta_1:=0$ и
$\theta_2:=\min\{\Theta_W(x)\}.$ Пусть
$\delta_1=\min_{x\in\mathfrak{S}_W:\Theta_W(x)\not=\emptyset}\delta_0(x)>0.$
Так как $M_x(0)\not=0$ для любого  $x\in\mathfrak{S}_W,$ то по
лемме Малкина (см. лемму~\ref{malkin_lem}) существуют
$\delta_*\in(0,\delta_1)$ и $\varepsilon_*>0$ такие, что
\begin{equation}\label{mark0}
   Q_\varepsilon \widetilde {x}\not=\widetilde {x}\quad\mbox{\rm для\ }
\widetilde{x}\in \overline{B_{\delta_*}
  \left(x\right)}\mbox{\rm \ при \ любых\ }
x\in\mathfrak{S}_W\mbox{\ \rm и\ }
  \varepsilon\in(0,\varepsilon_*).
\end{equation}
Пользуясь определением множества  $\mathfrak{S}_W,$ из
(\ref{mark0}) для любых $x\in\mathfrak{S}_W$ и
  $\varepsilon\in(0,\varepsilon_*)$
имеем
\begin{equation}\label{mark}
   Q_\varepsilon \widetilde{x}\not=\widetilde{x}\quad\mbox{\rm для\ любого\ }
\widetilde{x}\in \overline{B_{\delta_*}
  \left(x\right)}\cup\overline{B_{\delta_*}
  \left(S_{\min\{\Theta_W(x)\}}x\right)}.
\end{equation}
Пусть $\delta_{**}\in(0,\delta_*)$  достаточно мало так, что
\[
\left(B_{\delta_*} (x)\cup B_{\delta_*}(S_{\min\{\Theta_W(x)\}}
x)\cup
  W_{V_{\delta_{**}}(x)}\right)\backslash \overline{W}\subset B_{\delta_*}
  (x)\cup B_{\delta_*}(S_{\min\{\Theta_W(x)\}}
  x)
\]
для любых  $x\in \mathfrak{S}_W.$ Тогда, учитывая (\ref{mark}),
имеем
\[
  Q_\varepsilon \widetilde{x}\not=\widetilde{x}\quad{\rm for\ }
\widetilde{x}\in \left(B_{\delta_*} (x)\cup
B_{\delta_*}(S_{\min\{\Theta_W(x)\}} x)\cup
  W_{V_{\delta_{**}}(x)}\right)\backslash \overline{W}
\]
при  $x\in \mathfrak{S}_W$ и $\varepsilon\in(0,\varepsilon_*)$.
Следовательно, применяя формулу теоремы~\ref{MathNach1}, для
любого $x\in\mathfrak{S}_W$ такого, что
$\Theta_W(x)\not=\emptyset,$ и любого
$\varepsilon\in\left(0,\min\{\delta^{1+\alpha},\varepsilon_*\}\right)$
имеем
\begin{equation}
\begin{aligned}
&d\left(I-Q_\varepsilon,\left(B_{\delta_*} (x)\cup
B_{\delta_*}(S_{\min\{\Theta_W(x)\}} x)\cup
  W_{V_{\delta_{**}}(x)}\right)\cap W\right)=\\
&=d\left(I-Q_{\varepsilon},B_{\delta_*}
  (x)\cup B_{\delta_*}(S_{\min\{\Theta_W(x)\}} x)
  \cup W_{V_{\delta_{**}}(x)}\right)=\\
& = d(I-Q_\varepsilon,W_{V_{\delta_{**}}(x)})=-(-1)^{\beta(x)}
d_{\mathbb{R}}\left(M_x,\left(0,\min\{\Theta_W(x)\}\right)\right).\label{fo1}
\end{aligned}
\end{equation}
Пусть
\[
  \hskip-3cm\mathfrak{S}_W^0=\left\{x\in\mathfrak{S}_W:\mbox{\rm существует\
  }\delta_0>0\mbox{\rm\ такое,}\right.
\]
\[
 \left.\hskip3cm\mbox{\ что\ }S_\delta(x)\not\in\partial W\mbox{\rm\
  для\ всех\ }\delta\in(-\delta_0,0)\cup (0,\delta_0)\right\}.
\]
Из (\ref{mark0}) заключаем, что при каждых $ x\in\mathfrak{S}_W^0$
и $\varepsilon\in(0,\varepsilon_*)$ выполнено равенство
\begin{equation}\label{fo1rev}
d\left(I-Q_\varepsilon,B_{\delta_*}(x) \cap W\right)=0.
\end{equation}
Так как любая точка $x\in\mathfrak{S}_W$ является предельным
циклом системы (\ref{np_2}), и по предположению число таких точек
конечно, без ограничения общности можно считать, что $\delta_*>0$
настолько мало, что
\begin{equation}\label{bege}
\begin{array}{l}
 Q_0(\widehat{ x})\not=\widehat{ x} \quad \mbox{\rm для\
любых\ }\widehat{ x}\in C([0,T],\mathbb{R}^n)\mbox{\rm\ таких,}\\
\mbox{что\ }\widehat {x}(0)\in B_{\delta_*}(x([0,T]))\backslash
x([0,T]). \end{array}
\end{equation}
Следовательно, граница множества $W\backslash E_{\delta_*},$ где
\[
E_{\delta_*}:=\left(\bigcup_{x\in\mathfrak{S}_W:\Theta_W(x)\not=\emptyset}
\left(B_{\delta_*}
  (x)\cup B_{\delta_*}(S_{\Theta_W(x)} x)\cup
  W_{V_{\delta_{**}}(x)}\right)\cap W\right)\bigcup
\]
\[\bigcup\left(\bigcup_{x\in\mathfrak{S}_W^0}
  B_{\delta_*}(x)\cap W\right),
\]
не содержит $T$-периодических решений системы (\ref{np_2}). Но
результат Капетто-Мавена-Занолина (\cite{camaza}, следствие~2)
утверждает, что если оператор $Q_0$ не имеет неподвижных точек на
границе какого-нибудь открытого ограниченного множества $A\subset
C([0,T],\mathbb{R}^n),$ то
$$
d(I-Q_0,A)=(-1)^n d_{\mathbb{R}^n}(f,A\cap\mathbb{R}^n),
$$
поэтому
\begin{equation}\label{ma}
d(I-Q_0,W\backslash E_{\delta_*})=(-1)^n d_{\mathbb{R}^n}\left(f,
\left(W\backslash E_{\delta_*}\right)\cap\mathbb{R}^n\right).
\end{equation}
На основании (\ref{bege}) функция $f$ невырождена на множестве
$E_{\delta_*}\cap\mathbb{R}^n,$ и из  (\ref{ma}) заключаем, что
\begin{equation}\label{fo2}
  d(I-Q_0,W\backslash E_{\delta_*})=
  (-1)^n d_{\mathbb{R}^n}(f,W\cap\mathbb{R}^n).
\end{equation}
Суммируя (\ref{fo1}), (\ref{fo1rev}) и (\ref{fo2}), получаем
утверждение теоремы.

Теорема доказана.

\begin{rem}
Из формулы  (\ref{kmnform}) следует, что точки множества
$\mathfrak{S}_W$ такие, что  $S_\theta x\not\in\mathfrak{S}_W$ для
всех $\theta\in(0,T)$ не влияют на значение степени
$d(I-Q_\varepsilon,W),$ где $\varepsilon>0$ достаточно мало.
\end{rem}

Положим $X=\{x\in C([0,T],\mathbb{R}^n):x(0)=x(T)\}$ и обозначим
как $L:{\rm dom}L\subset X\to L^1([0,T],\mathbb{R}^n)$ линейный
оператор, определенный как $(Lx)(\cdot)=\dot x(\cdot)$ с ${\rm
dom} L=\{x\in X: x(\cdot) - \mbox{\rm\ абсолютно\ непрерывна}\}.$
Тогда  $L$ -- оператор Фредгольма нулевого индекса (см.
\cite{camaza}, пункт~II.1). Пусть $N_\varepsilon:X\to
L^1([0,T],\mathbb{R}^n)$ -- оператор Немыцкого, задаваемый как
$(N_\varepsilon x)(\cdot)=f(x(\cdot))+\varepsilon g
(\cdot,x(\cdot),\varepsilon).$ Таким образом, существование
$T$-периодических решений для системы (\ref{ps_2}) эквивалентно
разрешимости уравнения
\begin{equation}\label{coequ}
Lx=N_\varepsilon x,\quad x\in{\rm dom}L.
\end{equation}

В следующем утверждении предлагается формула для индекса
совпадения $D_L(L-N_\varepsilon,W\cap X)$ операторов $L$ и
$N_\varepsilon,$  (см. Ж.~Мавен \cite{mawbvp}, p.~19), подобная
формуле теоремы~\ref{MathNach2}.

\begin{cor}\label{maw_gen}
Предположим, что все условия теоремы~\ref{MathNach2} выполнены.
Тогда, для всех достаточно малых  $\varepsilon>0$ индекс
совпадения $D_L(L-N_\varepsilon,W\cap X)$ определен, и справедлива
следующая формула
\begin{equation}
\begin{aligned}
& D_L(L-N_\varepsilon,W\cap X)=\; (-1)^n
d_{\mathbb{R}^n}(f,W\cap\mathbb{R}^n)\;-\\
& -\sum_{x\in\mathfrak{S}_W: ~\Theta_W(x) \not=
\emptyset}(-1)^{\beta(x)}
d_{\mathbb{R}}\left(M_x,\left(0,\min\{\Theta_W(x)\}\right)\right).\label{kmnformD}
\end{aligned}
\end{equation}
\end{cor}

\noindent {\bf Доказательство.} Так как  $d(I-Q_\varepsilon,W)$
определен для достаточно малых $\varepsilon>0,$ то
$D_L(L-N_\varepsilon,W\cap X)$ также определен при этих значениях
$\varepsilon>0,$ см. (\cite{mawbvp}, Гл.~2, \S~2). Чтобы доказать
(\ref{kmnformD}), используем принцип родственности, разработанный
в (\cite{mawbvp}, Гл.~3). Во-первых, заметим, что нули оператора
$R_\varepsilon:C([0,T],\mathbb{R}^n)\to C([0,T],\mathbb{R}^n),$
определенного как
\[
(R_\varepsilon x)(t)=\; x(t)-x(0)-\int\limits_0^T
  \left(f(x(\tau))+\varepsilon g(\tau,x(\tau),\varepsilon)\right)d\tau
-
\]
\[
-\int\limits_0^t\left(f(x(\tau))+\varepsilon g
(\tau,x(\tau),\varepsilon)
  \right)d\tau
+t\int\limits_0^T
  \left(f(x(\tau))+\varepsilon
  g(\tau,x(\tau),\varepsilon)\right)d\tau,
\]
совпадают с неподвижными точками оператора  $Q_\varepsilon,$ и
значит $d(R_\varepsilon,W)$ также определен при достаточно малых
$\varepsilon>0.$ На основании принципа родственности для
топологический степеней эквивалентных интегральных операторов (см.
Ж.~Мавен \cite{mawbvp},
 теорема~III.1 при $a=1,$ $b=0$ и теорема~III.4)
\[
  d(R_\varepsilon,W)=d(I-Q_\varepsilon,W).
\]
Далее, используя для определения $D_L(L-N_\varepsilon,W\cap X)$
методы, разработанные в (\cite{mawbvp}, Гл.~III, \S~4) и принцип
родственности, связывающий индекс совпадения операторов $L$ и
$N_\varepsilon$ с топологической степенью эквивалентного
интегрального оператора (см. \cite{mawbvp}, теорема~III.7),
получаем равенство
\[
D_L(L-N_\varepsilon,W\cap X)=d(R_\varepsilon,W),
\]
которое завершает доказательство.

Следствие доказано.

Для перехода от условий теоремы~1.2 к основному предположению
$(A_0)$ (см. начало пункта) необходимы нижеследующие утверждения.

\begin{rem}\label{remark1} Имеют место соотношения
\[
  \mathfrak{S}_{W_U}=\mathfrak{S}^U
\]
и
\[
  \Theta_{W_U}(x)=\left\{\theta_0\in(0,T):x(\theta_0)\in\partial U,\
  x(\theta)\in U\mbox{\ \rm для\ всех\ }\theta\in(0,\theta_0)\right\}.
\]
\end{rem}

\begin{lem}\label{W_U_3} Пусть выполнено условие $(A_0).$ Тогда
$$
  d_{\mathbb{R}^n}(f,W_U\cap\mathbb{R}^n)=d_{\mathbb{R}^n}(f,U).
$$
\end{lem}

\noindent{\bf Доказательство.} Доказательство леммы производится
на основании принципа продолжения Лерэ-Шаудера. Пусть
\[
  U_\lambda=\left\{\xi\in\mathbb{R}^n: \Omega(0,\lambda t,\xi)  \in U\mbox{\rm\ для\
  любого\ }t\in[0,T]\right\},\quad\lambda\in[0,1],
\]
покажем, что
\begin{equation}\label{impo}
  0\not\in \widetilde{M}(\partial U_\lambda)\quad\mbox{\rm для\ любого\
  }\lambda\in[0,1].
\end{equation}
Предположим противное, тогда существует $\lambda_0\in[0,1]$ такое,
что $\xi_0\in\partial U_{\lambda_0}$ и $\widetilde{M}(\xi_0)=0.$
Заметим, что  $\Omega(0,\lambda_0 t,\xi_0)\in\overline{U}$ для
любого $t\in[0,T].$ Следовательно, существует $t_0\in[0,T]$ такое,
что  $\Omega(0,\lambda_0 t_0,\xi_0)\in\partial U,$ и на основании
факта $\widetilde{M}(\xi_0)=0,$ получаем, что  $\Omega(0,\lambda_0
t,\xi_0)$ постоянно по отношению к  $t\in[0,T].$ Значит, мы имеем
$\Omega(0,\lambda_0 t_0,\xi_0)=\Omega(0,0,\xi_0)=\xi_0$ и
получаем, что $\xi_0\in\partial U,$ противореча предположению о
том, что $\partial U$ содержит только начальные условия предельных
циклов системы (\ref{np_2}). Используя принцип продолжения
Лерэ-Шаудера (см. \cite{lersch}, Гл.~3, \S~16, "Фундаментальная
теорема" или \cite{bro}, теорема~10.7) из (\ref{impo}) заключаем,
что
\[
  d_{\mathbb{R}^n}(f,U_0)=d_{\mathbb{R}^n}(f,U_1).
\]
С другой стороны,  $U_0=U$ и $U_1=W_U\cap\mathbb{R}^n.$

Лемма доказана.

Положим $\Theta^U(\widetilde{x})=\Theta_{W_U},$ то есть

\centerline{$\Theta^U(\widetilde{x})=\left\{\theta_0\in(0,T):\widetilde{x}(\theta_0)\in\partial
U,\ \widetilde{x}(\theta)\in U, \ \theta\in(0,\theta_0)\right\}.$}

Напомним, что согласно условию $(A_0),$
\[
  \mathfrak{S}^U=\bigcup_{\xi\in\partial
  U:\Omega(T,0,\xi)=\xi}\left\{x\in C([0,T],\mathbb{R}^n):x(t)=\Omega(t,0,\xi),\
  t\in[0,T]\right\}.
\]

На основании замечания~1.2 и леммы~\ref{W_U_3} получаем следующее
следствие из теоремы~1.2.

\begin{thm}\label{main_3} Пусть выполнено условие $(A_0).$
Предположим, что $M_x(0)\not=0$ для всех $x\in\mathfrak{S}^U.$
Тогда для каждого достаточно малого $\varepsilon>0$ топологическая
степень $d(I-Q_\varepsilon,W_U)$ определена, и справедлива
следующая формула
$$
d(I-Q_\varepsilon,W_U)=(-1)^n d_{\mathbb{R}^n}(f,U)-
$$
\begin{equation}\label{kmnform_}
-\sum_{x\in\mathfrak{S}^U: ~\Theta^U(x)\not=
\emptyset}(-1)^{\beta(x)}
d_{\mathbb{R}}\left(M_x,\left(0,\min\{\Theta^U(x)\}\right)\right),
\end{equation}
\end{thm}

\noindent Предположим теперь, что

$(A_{\mathcal{P}})$ решение $x$ системы (\ref{ps_2}) с начальным
условием $x(t_0)=\xi$ существует, единственно и продолжимо на
отрезок $[0,T]$ при любых $t_0\in[0,T],$ $\xi\in\mathbb{R}^n$ и
$\varepsilon>0.$

При условии $(A_{\mathcal{P}})$ для системы (\ref{ps_2}) при любых
$\varepsilon>0$  определен оператор $\mathcal{P}_\varepsilon$
Пуанкаре, соответствующий задаче о $T$-периодических решениях для
(\ref{ps_2}) (см. определение~1.2). В этом случае можно
сформулировать аналогичное теореме~1.3 утверждение о
топологической степени оператора $I-\mathcal{P}_\varepsilon$ на
$U.$

Действительно, справедлив следующий основной результат.

\begin{thm}\label{cor2} Предположим, что  выполнены условия $(A_0)$ и $(A_{\mathcal{P}}).$ Если
$M_{\widetilde{x}}(0)\not=0$  для любого $\widetilde{x}\in
\mathfrak{S}^U,$ то для всех достаточно малых $\varepsilon>0$
топологическая степень
 $d_{\mathbb{R}^n}(I-\mathcal{P}_\varepsilon,U)$
определена и может быть посчитана по формуле
$$
d_{\mathbb{R}^n}(I-\mathcal{P}_\varepsilon,U)=(-1)^n
d_{\mathbb{R}^n}(f,U)-
$$
\begin{equation}\label{kmnformU}
-\sum_{\widetilde{x}\in\mathfrak{S}^U:~\Theta^U(\widetilde{x})\not=
\emptyset}(-1)^{\beta(\widetilde{x})}
d_{\mathbb{R}}\left(M_{\widetilde{x}},\left(0,\min\{\Theta^U(\widetilde{x})\}\right)\right).
\end{equation}
\end{thm}

\noindent{\bf Доказательство.} Из теоремы~\ref{MathNach2},
учитывая замечание~\ref{remark1}, заключаем, что существует
$\varepsilon_0>0$ такое, что для любого
$\varepsilon\in(0,\varepsilon_0]$ степень
 $d(I-Q_\varepsilon,W_U)$ определена, и
$$
d(I-Q_\varepsilon,W_U)=(-1)^n
d_{\mathbb{R}^n}(f,W_U\cap\mathbb{R}^n)-
$$
\begin{equation}\label{kmnformU1}
-\sum_{x\in\mathfrak{S}^U: ~\Theta^U(x)\not=
\emptyset}(-1)^{\beta(x)}
d_{\mathbb{R}}\left(M_x,\left(0,\min\{\Theta^U(x)\}\right)\right).
\end{equation}
Следовательно, чтобы доказать следствие, достаточно показать, что
\begin{equation}\label{show1}
d(I-Q_\varepsilon,W_U)=d_{\mathbb{R}^n}(I-\mathcal{P}_\varepsilon,U)\quad\mbox{\rm
для\ любого\ }\varepsilon\in(0,\varepsilon_0],
\end{equation}
и
\begin{equation}\label{show2}
  d_{\mathbb{R}^n}(f,W_U\cap\mathbb{R}^n)=d_{\mathbb{R}^n}(f,U).
\end{equation}
Справедливость (\ref{show2}) следует из леммы~\ref{W_U_3}. Чтобы
доказать (\ref{show1}), определим $W_U^\varepsilon\subset
C([0,T],\mathbb{R}^n)$ как
\[
  W_U^\varepsilon=\left\{\widehat {x}\in C([0,T],\mathbb{R}^n):
  \Omega_\varepsilon(0,t,\widehat {x}(t))\in U,\ \mbox{\rm для\ любого\
  }t\in[0,T]\right\},
\]
где $\Omega_\varepsilon$ -- оператор сдвига по траекториям
возмущенной системы (\ref{ps_2}). Утверждается, что существует
$\widehat{\varepsilon}_0\in(0,\varepsilon_0]$ такое, что
\begin{equation}\label{stepI}
  Q_\varepsilon x\not=x\quad\mbox{\rm для\ любых\
  }x\in  \left(W_U\backslash W_U^\varepsilon\right)\cup
  \left(W_U^\varepsilon \backslash W_U \right) \mbox{\ \rm и\
  }\varepsilon\in(0,\widehat{\varepsilon}_0].
\end{equation}
Предположим противное, тогда существуют последовательности
$\{\varepsilon_k\}_{k\in\mathbb{N}}\subset(0,\varepsilon_0],$
$\varepsilon_k\to 0$ при $k\to\infty,$
$\{x_k\}_{k\in\mathbb{N}}\subset C([0,T],\mathbb{R}^n)$ такие, что
\begin{equation}\label{propI}x_k\in\left(W_U\backslash
W_U^{\varepsilon_n}\right)\cup\left(W_U^{\varepsilon_n} \backslash
W_U \right),\end{equation} и
\begin{equation}\label{qwq}
x_k\to \widetilde{x}\mbox{\rm\ при\ } k\to\infty,\mbox{\rm\ \ где\
\ }Q_{\varepsilon_k}x_k=x_k.
\end{equation}
Легко видеть, что  (\ref{propI}) означает
$\widetilde{x}\in\partial W_U.$ Этот факт вместе с  (\ref{qwq}) и
предположением $M_{\widetilde{x}}(0)\not=0$ приводит к
противоречию с утверждением леммы Малкина (\ref{malkin_lem}).
Следовательно,  утверждение (\ref{stepI}) справедливо и значит
\[
  d(I-Q_\varepsilon,W_U)=d(I-Q_\varepsilon,W_U^\varepsilon)\quad\mbox{\rm
  для\ любого\ }\varepsilon\in(0,\widehat{\varepsilon}_0].
\]
Так как для любого $\varepsilon\ge 0$ множества $U$ и
$W_U^\varepsilon$ имеют общую серцивину по отношению к задаче о
$T$-периодических решениях для (\ref{ps_2}) (см. \cite{krazab}, \S
28.5), то  (см. \cite{krazab}, теорема~28.5) имеем
\[
  d(I-Q_\varepsilon,W_U^\varepsilon)=d_{\mathbb{R}^n}
  (I-\mathcal{P}_\varepsilon,U)\quad\mbox{\rm
  для\ любого\ }\varepsilon\ge 0
\]
и, таким образом,  (\ref{show1}) доказано.

Теорема доказана.

\begin{rem} Из формулы (\ref{kmnformU}) следует, что если цикл  $x\in\mathfrak{S}^U$
имеет с границей $\partial U$ единственную общую точку, то этот
цикл не влияет на $d_{\mathbb{R}^n}(I-\mathcal{P}_\varepsilon,U),$
где $\varepsilon>0$ достаточно мало.
\end{rem}

\section{Теоремы о продолжении $T$-периодических решений из
$\overline{U}$ по параметру}

На основании подходящего выбора множества  $W\subset
C([0,T],\mathbb{R}^n),$ участвующего в формуле теоремы~1.2, в этом
пункте формулируются некоторые результаты о существовании
$T$-периодических решений для системы (\ref{ps_2}), где
$\varepsilon\in(0,\varepsilon_0],$ сходящихся при $\varepsilon\to
0$ к лежащим в $\overline{U}$ $T$-периодическим решениям системы
(\ref{np_2}). Такие результаты называются результатами о
продолжении $T$-периодических решений системы (\ref{ps_2}) при
увеличении параметра $\varepsilon>0$ от нуля до $\varepsilon_0>0$
(см. \cite{camaza}).

\begin{thm}\label{thm3} Предположим, что все непостоянные  $T$-периодические решения системы
 (\ref{np_2}) являются простыми циклами этой системы. Тогда, для любого открытого ограниченного множества
$W\subset C([0,T],\mathbb{R}^n),$ содержащего все постоянные
решения системы (\ref{np_2}) и удовлетворяющего условиям
\[
\mathfrak{S}_W\mbox{\it\ конечно},\quad M_x(0)\not=0\quad\mbox{\it
для\ любого\ }x\in \mathfrak{S}_W
\]
и
\[
(-1)^n d_{\mathbb{R}^n}(f,W\cap\mathbb{R}^n)-
\]
\[-\sum\limits_{x\in\mathfrak{S}_W:
~\Theta_W(x)\not=\emptyset}(-1)^{\beta(x)}
d_{\mathbb{R}}\left(M_x,\left(0,\min\{\Theta_W(x)\}\right)\right)
\not=0,
\]
существует  $\varepsilon_0>0$ такое, что при всех $\varepsilon\in
(0,\varepsilon_0]$ система  (\ref{ps_2}) имеет $T$-периодическое
решение в  $W.$
\end{thm}

Предположения теоремы~\ref{thm3} означают, что множество
$\mathfrak{S}_W$ содержит только простые циклы системы
(\ref{np_2}). Следовательно, теорема~\ref{thm3} следует из
теоремы~\ref{MathNach2} и теоремы Лерэ-Шаудера о неподвижной точке
(см. \cite{lersch}, Гл.~1, \S~7 или \cite{krazab}, теорема~20.5).

\begin{cor}\label{cor3} Предположим, что все непостоянные
$T$-периодические решения системы (\ref{np_2}) являются простыми
циклами системы (\ref{np_2}). Предположим, что существует открытое
ограниченное множество $W\subset C([0,T],\mathbb{R}^n),$
содержащее все постоянные решения системы  (\ref{np_2}) и
удовлетворяющее условию
\begin{equation}\label{one}
\begin{array}{l}\mathfrak{S}_W\mbox{\it\ конечно}, \quad M_x(0)\cdot
M_x(\min\{\Theta_W(x)\})>0\\
\quad\mbox{\it для\ любого\ }x\in \mathfrak{S}_W\mbox{\it \
такого, что\ }\Theta_W(x)\not=\emptyset
\end{array}
\end{equation}
и условию
\[
(-1)^n d_{\mathbb{R}^n}(f,W\cap\mathbb{R}^n)\not=0.
\]
Тогда, при достаточно малых  $\varepsilon>0$ система  (\ref{ps_2})
имеет $T$-периодическое решение в $W.$
\end{cor}

Доказательство следствия~\ref{cor3} вытекает из утверждения
(\ref{one}), означающего (см. \cite{krazab}, \S~3.2), что $
d_{\mathbb{R}}\left(M_x,\left(0,\min\{\Theta_W(x)\}\right)\right)=
0$ для любого $x\in \mathfrak{S}_W$ такого, что
$\Theta_W(x)\not=\emptyset.$

\vskip0.2truecm Рассмотрим теперь некоторые приложения
теоремы~\ref{MathNach1} к задаче о существовании в системе
(\ref{ps_2})  $T$-периодических решений близких к простому циклу
системы (\ref{np_2}). Вначале установим следующий результат.

\begin{thm}\label{thm4} Пусть $\widetilde{x}$ -- простой  $T$-периодический
цикл системы (\ref{np_2}). Пусть $0\le\theta_1<\theta_2\le
\theta_1+\frac{T}{p},$ где  $p\in\mathbb{N}$ и $\frac{T}{p}$
является наименьшим периодом цикла  $\widetilde{x}.$ Предположим,
что
\begin{equation}\label{icond}
M_{\widetilde{x}}(\theta_1)\cdot M_{\widetilde{x}}(\theta_2)<0.
\end{equation}
Обозначим через  $\Theta$ множество всех нулей функции
$M_{\widetilde{x}}$ на $(\theta_1,\theta_2).$  Тогда, для
достаточно малых $\varepsilon>0,$ система (\ref{ps_2}) имеет
$T$-периодическое решение $x_\varepsilon$ такое, что
\begin{equation}\label{conve}
\rho\left(x_\varepsilon(t),\widetilde{x}(t+\Theta)\right)\to
0\quad \mbox{ при\
  }\varepsilon\to 0.
\end{equation}
\end{thm}

\noindent {\bf Доказательство.} Заметим, что условие (\ref{icond})
означает
\[
M_{\widetilde{x}}(\theta_1)\not=0\quad\mbox{\rm и}\quad
M_{\widetilde{x}}(\theta_2)\not=0
\]
и, таким образом, условия теоремы~\ref{MathNach1} удовлетворены.
Зафиксируем $\alpha>0,$ из теоремы~\ref{MathNach1} имеем
$\delta_0>0$ такое, что для любого
$\varepsilon\in(0,\delta_0^{1+\alpha})$ топологическая степень
$d(I-Q_\varepsilon,W_{V_{\delta(\varepsilon)}})$ определена при $
\delta(\varepsilon)=\varepsilon^{1/(1+\alpha)},$ и
\[
d(I-Q_\varepsilon,W_{V_{\delta(\varepsilon)}})=
(-1)^{\beta(\widetilde{x})}d_{\mathbb{R}}(M_{\widetilde{x}},(\theta_1,\theta_2)).
\]
Из (\ref{icond}) также имеем (см. \cite{krazab}, \S 3.2), что
$|d_{\mathbb{R}}(M_{\widetilde{x}},(\theta_1,\theta_2))|=1$ и,
таким образом, для любого $\varepsilon\in(0,\delta_0^{1+\alpha})$
система (\ref{ps_2}) имеет $T$-периодическое решение
$x_\varepsilon$ такое, что  $x_\varepsilon(0)\in
V_{\delta(\varepsilon)}.$ Более того, из свойства  1)
теоремы~\ref{MathNach1} заключаем
\begin{equation}\label{ss1}
\rho\left(x_\varepsilon(0),\widetilde{x}([\theta_1,\theta_2])\right)\le
\delta(\varepsilon)=\varepsilon^{1/(1+\alpha)}.
\end{equation}
Пусть $\nu_\varepsilon(t)=\Omega(0,t,x_\varepsilon(t)),$ тогда
согласно лемме~\ref{zamena}
\[
\dot \nu_\varepsilon(t)=\varepsilon \left(\Omega'_{(3)}(
t,0,\nu_\varepsilon(t))\right)^{-1}g(t,x(t,\nu_\varepsilon(t),\varepsilon)).
\]
Следовательно, существует $M_1>0$ такое, что
\begin{equation}\label{ss2}
  \|\nu_\varepsilon(0)-\nu_\varepsilon(t)\|\le M_1\varepsilon\quad \mbox{\rm для\ любых\
  }\varepsilon\in(0,\delta_0^{1+\alpha})\mbox{\rm\ и\ }t\in[0,T].
\end{equation}
С другой стороны, $\nu_\varepsilon(0)=x_\varepsilon(0)$ и, таким
образом, из (\ref{ss1}) и (\ref{ss2}) для любого
  $\varepsilon\in(0,\delta_0^{1+\alpha})$  и любого $t\in[0,T]$
  справедливо соотношение
$$
  \rho(\nu_\varepsilon(t),\widetilde{x}([\theta_1,\theta_2]))\le
  \|\nu_\varepsilon(t)-x_\varepsilon(0)\|+
$$
\begin{equation}\label{uu1}
  +\rho(x_\varepsilon(0),\widetilde{x}([\theta_1,\theta_2]))\le
  \varepsilon^{1/(1+\alpha)}\left(1+M_1
  \varepsilon^{\alpha/(1+\alpha)}\right).
\end{equation}
Так как для любого  $\theta\in[\theta_1,\theta_2]$ имеем
$\|x_\varepsilon(t)-\widetilde{x}(t+\theta)\|=\|x(t,\nu_\varepsilon(t))-
x(t,\widetilde{x}(\theta))\|,$ и, как это уже было замечено  в
доказательстве теоремы~\ref{MathNach1}, функция $x(\cdot,\cdot)$
непрерывно дифференцируема по обоим переменным, заключаем, что
существует  $M_2>0$ такое, что
\begin{equation}\label{uu2}
  \|x_\varepsilon(t)-\widetilde{x}(t+\theta)\|\le M_2\|\nu_\varepsilon(t)-\widetilde{x}(\theta)\|
\end{equation}
для любых $\varepsilon\in(0,\delta_0^{1+\alpha}),$ $t\in[0,T],$
  $\theta\in[\theta_1,\theta_2].$
Подставляя (\ref{uu1}) в (\ref{uu2}), получаем, что для любых
$\varepsilon\in(0,\delta_0^{1+\alpha})$ и $t\in[0,T]$ выполнено
\begin{equation}\label{fin0}
\rho(x_\varepsilon(t),\widetilde{x}(t+[\theta_1,\theta_2]))\le
  \varepsilon^{1/(1+\alpha)}M_2\left(1+M_1
  \varepsilon^{\alpha/(1+\alpha)}\right).
\end{equation}
Предположим теперь, что (\ref{conve}) не удовлетворено,
следовательно, существует  $\delta_*>0,$ а также
последовательности
$\{\varepsilon_k\}_{k\in\mathbb{N}}\subset(0,\delta_0^{1+\alpha}),$
$\varepsilon_k\to 0$ при $k\to\infty,$ и
$\{t_k\}_{k\in\mathbb{N}}\subset[0,T]$ такие, что
\begin{equation}\label{popo}
  x_{\varepsilon_k}(t_k)\not\in B_{\delta_*}(\widetilde{x}(t_k+\Theta))\quad\mbox{\rm
  для\ любого\ }k\in\mathbb{N}.
\end{equation}
Без ограничения общности можем считать, что последовательности
$\{x_k\}_{k\in\mathbb{N}}$ и $\{t_k\}_{k\in\mathbb{N}}$ сходятся.
Тогда пользуясь (\ref{fin0}), получаем существование
$\theta_*\in[\theta_1,\theta_2]$ такого, что
\begin{equation}\label{po}
  x_k(t)\to \widetilde{x}(t+\theta_*)\quad\mbox{\rm при\ }k\to\infty
\end{equation}
равномерно по отношению к $t\in[0,T].$ На основании леммы Малкина
(лемма~\ref{malkin_lem})  из (\ref{po}) заключаем, что
$M_{\widetilde{x}}(\theta_*)=0.$ С другой стороны, из (\ref{popo})
имеем  $\widetilde{x}(t_0+\theta_*)\not\in
B_{\delta_*/2}(\widetilde{x}(t_0+\Theta)),$ где
$t_0=\lim_{k\to\infty} t_k,$ и, значит,
$\widetilde{x}(\theta_*)\not\in
B_{\delta_*/2}(\widetilde{x}(\Theta)).$ Достигнутое противоречие
доказывает справедливость соотношения (\ref{conve}).

Теорема доказана.

\begin{rem}\label{remark2} Доказательство теоремы~\ref{thm4} предоставляет
информацию о скорости сходимости $T$-периодических решений системы
(\ref{ps_2}) к простому циклу системы (\ref{np_2}). Действительно,
из (\ref{fin0}) следует, что расстояние между траекторией
$T$-периодического решения $x_\varepsilon$ и простым циклом
$\widetilde{x}$ имеет порядок $\varepsilon^{1/(1+\alpha)},$ где
$\alpha>0$ -- произвольное заранее фиксированное число.
\end{rem}

В главе~\ref{sec_conv} будет установлено, что, на самом деле,
имеет место б\'ольшая скорость, чем указанная в замечании~1.4, но
это потребует значительно более долгих рассуждений.

\begin{cor}\label{cor4} Пусть $\widetilde{x}$ -- простой
$T$-периодический цикл ситемы (\ref{np_2}). Предположим, что
существует  $\theta_0\in[0,T]$ такое, что
\begin{equation}\label{monot}
M_{\widetilde{x}}(\theta_0)=0\mbox{\quad \it
и\quad}M_{\widetilde{x}}\mbox{\it\ строго\ монотонна\ в\
}\theta_0.
\end{equation}
Тогда для всех достаточно малых $\varepsilon>0$ система
(\ref{ps_2}) имеет $T$-периодическое решение $x_\varepsilon,$
удовлетворяющее
\begin{equation}\label{conv_2_}
  x_\varepsilon(t)\to \widetilde{x}(t+\theta_0)\quad\mbox{ при\ }\varepsilon\to
  0.
\end{equation}
\end{cor}

Следствие~\ref{cor4} непосредственно вытекает из
теоремы~\ref{thm4}, в которой  $\theta_1<\theta_0<\theta_2$ взяты
достаточно близкими к $\theta_0.$

\begin{cor}\label{cor5} Предположим, что $f$ -- один раз и $g$ -- три раза непрерывно дифференцируемые функции.
Пусть $\widetilde{x}$ -- простой  $T$-периодический цикл системы
(\ref{np_2}). Предположим, что для некоторого $\theta_0\in[0,T]$
выполнено
\[
  M_{\widetilde{x}}(\theta_0)=M_{\widetilde{x}}'(\theta_0)=M_{\widetilde{x}}''(\theta_0)=0,\quad
  M_{\widetilde{x}}'''(\theta_0)\not=0.
\]
Тогда для достаточно малых  $\varepsilon>0$ система (\ref{ps_2})
имеет  $T$-периодическое решение $x_{\varepsilon},$
удовлетворяющее свойству (\ref{conv_2_}).
\end{cor}

\subsubsection{Формула для фазы
 $T$-периодических решений синусоидально возмущенных систем}

 Предположим, что возмущенная
система (\ref{ps_2}) имеет вид
\begin{equation}\label{pss_2}
  \dot x=f(x)+\varepsilon\left(\begin{array}{c} 0 \\
  \sin\left(\dfrac{2\pi k}{T}t\right)g(x)\end{array}\right),
\end{equation}
где $g:\mathbb{R}^2\to\mathbb{R}$ -- непрерывная функция  и
$k\in\mathbb{N}.$

Тогда для функции Малкина (\ref{fun_mal}) имеем следующее
представление
$$
  M_{\widetilde{x}}(\theta)=\cos\left(\frac{2\pi
  k}{T}\theta\right)M_{\sin}-\sin\left(\frac{2\pi
  k}{T}\theta\right)M_{\cos},
$$
где
$$
  M_{\sin}={\rm
  sign}\left<\dot{\widetilde{x}}(0),\widetilde{z}(0)\right>\int\limits_0^T\widetilde{z}_2(\tau)
\sin\left(\frac{2\pi
  k}{T}\tau\right)g(\widetilde{x}(\tau))d\tau,
$$
$$
  M_{\cos}={\rm
  sign}\left<\dot{\widetilde{x}}(0),\widetilde{z}(0)\right>\int\limits_0^T\widetilde{z}_2(\tau)
\cos\left(\frac{2\pi
  k}{T}\tau\right)g(\widetilde{x}(\tau))d\tau.
$$

Пользуясь указанным представлением и следствием~\ref{cor4},
получаем следующее утверждение.

\begin{cor}\label{lem_phase}
  Пусть $\widetilde{x}$ -- простой
$T$-периодический цикл системы (\ref{np_2}). Предположим, что
$M_{\cos}\not=0.$ Тогда каждому числу
$$
  \theta_j=\frac{T\arctg\left(M_{\sin} /M_{\cos}\right)+T\pi j}{2\pi k},
  \quad j=1,...,2k,
$$
и всем достаточно малым $\varepsilon>0$ соответствует
$T$-периодическое решение $x_{j,\varepsilon}$ системы
(\ref{pss_2}) такое, что
$$
x_{j,\varepsilon}(t)\to \widetilde{x}(t+\theta_j)\quad\mbox{\rm
при\ }\varepsilon\to
  0,\quad j=1,...,2k.
$$
\end{cor}

Действительно, подстановкой проверяется, что числа $\theta_j,$
$j=1,...,2k,$ удовлетворяют уравнению
$M_{\widetilde{x}}(\theta)=0$ и свойству
$\left(M_{\widetilde{x}}\right)'(\theta)\not=0.$

Полученное следствие позволяет гарантировать существование
$T$-периодических решений в тех случаях, где аналитическое
вычисление решений порождающей системы (\ref{np_2})
затруднительно, а проверка справедливости неравенства
$M_{\cos}\not=0$ возможна.

\begin{exa} \label{preyexa} \rm Для обобщенной системы Хищник-Жертва (см. \cite{gor}, \S 5.3)
\begin{equation}\label{predat_ps}
 \begin{array}{lll}
  \dot x_1 & = & k_1 x_1-\dfrac{k_2 x_1 x_2}{k_0+x_1}-k_3 x_1^2,\\
  \dot x_2 & = & \dfrac{k_4 x_1 x_2}{k_0+x_1}-k_5 x_2
  +\varepsilon(\mu x_1^+ +\nu x_1^-)
  \sin\left(\dfrac{2\pi}{T}t\right),
 \end{array}
\end{equation}
\end{exa}
где $a^+:=\max\{a,0\},$ $a^-:=\max\{-a,0\},$ допускающей в
определенной области изменения параметров $k_0,...,k_5$ простой
цикл $\widetilde{x}$ некоторого периода $T>0,$
следствие~\ref{lem_phase} гарантирует существование по крайней
мере двух $T$-периодических решений вблизи $\widetilde{x}$ для
достаточно малых $\varepsilon>0$ и  таких $\mu,\nu\in\mathbb{R},$
при которых $M_{\cos}\not=0.$

\section{Сопоставление полученных результатов с имеющимися в
литературе}

В случае, когда
\begin{equation}\label{mawc1}
Q_0 x\not=x\quad \mbox{\rm для\ любого\ } x\in\partial W,
\end{equation}
и
\begin{equation}\label{mawc2}
  d_{\mathbb{R}^n}(f,W\cap\mathbb{R}^n)\not=0,
\end{equation}
задача о существовании $T$-периодических решений для (\ref{ps_2})
решена А.~Капетто, Ж.~Мавеном и Ф.~Занолином в \cite{camaza}. Они
установили (\cite{camaza}, следствие~1), что при условиях
(\ref{mawc1}) и (\ref{mawc2}) справедлива формула
\begin{equation}\label{mform}
  d(I-Q_0,W)=(-1)^n d_{\mathbb{R}^n}(f,W\cap\mathbb{R}^n),
\end{equation}
впервые полученная М.~А.~Красносельским и А.~И.~Перовым
 для общего случая
неавтономной порождающей системы, см. \cite{perov} (с.~108) или
\cite{kraper}. В случае, когда решения системы (\ref{ps_2})
удовлетворяют условиям единственности и нелокальной продолжимости,
формула (\ref{mform}) для частных случаев множеств $W$ установлена
И.~Берштейном и А.~Халанаем \cite{berhal}.

 Из (\ref{mform})  следует, что
\begin{equation}\label{mform1}
  d(I-Q_\varepsilon,W)=(-1)^n d_{\mathbb{R}^n}(f,W\cap\mathbb{R}^n)
\end{equation}
для любых достаточно малых $\varepsilon>0.$ Следовательно, при
условиях (\ref{mawc1}) и (\ref{mawc2}) система (\ref{ps_2}) имеет
$T$-периодическое решение в $W$ для любого $T$-периодического по
первой переменной возмущения $g$ и любого достаточно малого
$\varepsilon>0.$ Заметим, что предположение (\ref{mawc2})
означает, что множество  $W$ обязательно содержит постоянное
решение порождающей системы (\ref{np_2}).

В настоящей главе условие  (\ref{mawc1}) не требуется, то есть
резрешается, чтобы  $\partial W$ содержало неподвижные точки
оператора $Q_0,$ и полученная формула (\ref{kmnform})
теоремы~\ref{MathNach2} является обобщением формулы
(\ref{mform1}). Отметим, что теорема~\ref{MathNach2} может
гарантировать, что $d(I-Q_\varepsilon,W)\neq 0$ даже в случае,
когда $d_{\mathbb{R}^n}(f,W\cap\mathbb{R}^n)=0,$ то есть без
явного требования того, что множество  $W$ содержит постоянное
решение системы (\ref{np_2}).

Второй член в правой части формулы (\ref{kmnform}) схож с
аналогичным членом формулы Красносельского-Забрейко для подсчета
индекса вырожденной неподвижной точки оператора $Q_0$ на основе
сужения этого оператора на подпространство (в нашем случае на
одномерное), см. (\cite{krazab}, формула 24.13). Однако,
соответствующая теорема, полученная М.~А.~Красносельским и
П.~П.~Забрейко (\cite{krazab}, теорема~24.1), может быть применена
только в случае, когда оператор $Q_0$ имеет специальную форму,
гарантирующую, что  $Q_0$ имеет только изолированные неподвижные
точки. Такое свойство в рассматриваемом случае не выполнено, так
как любое $T$-периодическое решение автономной системы
(\ref{np_2}) является неизолированной неподвижной точкой оператора
$Q_0.$

Утверждение следствия~\ref{cor3} совпадает с утверждением
Капетто-Мавена-Занолина  (\cite{camaza}, теорема~2), но в
последней работе дополнительно требуется, чтобы множество
$\partial W$ не содержало $T$-периодических решений системы
(\ref{np_2}).

В работе \cite{mal} И.~Г.~Малкиным установлен следующий результат
(см. \cite{mal}, утверждение~с.~638 или \cite{malb}, теоремы
сс.~387 и 392). {\it Пусть
 $f:\mathbb{R}^n\to\mathbb{R}^n$ -- один раз и $g:\mathbb{R}\times\mathbb{R}^n\times[0,1]\to\mathbb{R}^n$
 -- два раза непрерывно
дифференцируемые функции. Пусть $\widetilde{x}$ -- простой
$T$-периодический цикл системы (\ref{np_2}). Предположим, что
существует $\theta_0\in[0,T]$ такое, что
$M_{\widetilde{x}}(\theta_0)=0$ и
\begin{equation}\label{mmm}
(M_{\widetilde{x}})'(\theta_0)\not=0.
\end{equation}
Тогда, для всех достаточно малых $\varepsilon>0$ система
(\ref{ps_2}) допускает $T$-периодическое решение $x_\varepsilon,$
удовлетворяющее свойству
\begin{equation}\label{conv_2}
  x_\varepsilon(t)\to \widetilde{x}(t+\theta_0)\quad\mbox{\rm при\ }\varepsilon\to
  0.
\end{equation}}
\noindent Таким образом, установленные в настоящей главе
теорема~1.6 и  следствие~\ref{cor4} являются обобщением теоремы
Малкина на случай, когда вместо (\ref{mmm}) имеется либо свойство
$M_{\widetilde{x}}(\theta_1)\cdot M_{\widetilde{x}}(\theta_2)<0,$
где $\theta_1<\theta_0<\theta_2,$ либо строгая монотонность
функции $M_{\widetilde{x}}$ в $\theta_0,$  а также на случай
недифференцируемых правых частей.

Случай, когда (\ref{mmm}) не выполнено, был исследован В.~Лудом в
\cite{loud}. Для того, чтобы сформулировать его результат, введем
некоторые обозначения. Вначале, повернем и перенесем координатные
оси так, что  $\widetilde{x}(0)=0$ и
$\dot{\widetilde{x}}(0)=\left(\widetilde{x}_1(0),0,...,0\right).$
Пусть $\Omega_\varepsilon$ -- оператор сдвига по траекториям
возмущенной системы (\ref{ps_2}). Положим
$F(\xi,\varepsilon)=\Omega_\varepsilon(T,0,\xi)-\xi,$ так как цикл
$\widetilde{x}$ простой, то  $n-1$ уравнений системы
$F(\xi,\varepsilon)=0$ могут быть решены вблизи $0$ по отношению к
$\xi^k,$ где $k\in\{ 1,2,...,n\},$ и в результате получим
скалярное уравнение $H(u,\varepsilon)=0.$ Пусть
$D_{\widetilde{x}}$ -- дескриминант уравнения
$$
\frac{1}{2}\frac{\partial^3 H}{\partial u^2\partial
\varepsilon}(0,0)m^2+ \frac{1}{2}\frac{\partial^3 H}{\partial
u\partial\varepsilon^2}(0,0)m+ \frac{1}{6}\frac{\partial^3
H}{\partial\varepsilon^3}(0,0)=0.
$$
В.~Луд (\cite{loud}, теорема~2) установил, что {\it если  $f$ --
три раза и $g$ -- два раза непрерывно дифференцируемые функции,
\begin{equation}\label
{desc} D_{\widetilde{x}}>0,
\end{equation}
и для некоторого  $\theta_0\in[0,T]$ удовлетворяющего
$M_{\widetilde{x}}(\theta_0)=0$ выполнено
$(M_{\widetilde{x}})'(\theta_0)=0$ и
\begin{equation}\label{llll}
(M_{\widetilde{x}})''(\theta_0)=0,
\end{equation}
то для всех достаточно малых $\varepsilon>0$ система (\ref{ps_2})
имеет $T$-периодическое решение $x_{\varepsilon},$ удовлетворяющее
условию сходимости (\ref{conv_2}).}

В случае, когда $M_{\widetilde{x}}(\cdot)$ -- тождественный нуль,
В.~Луд в \cite{loud} выводит из приведенного результата теорему о
существовании  $T$-периодических решений для (\ref{ps_2}) вблизи
$\widetilde{x}.$ Но даже в случае, когда
$(M_{\widetilde{x}})'''(\theta_0)\not=0,$ проверка условия
(\ref{desc}) -- далеко не очевидная задача (здесь предполагается,
что  $g$ три раза непрерывно дифференцируемая функция). Это
замечание обуславливает актуальность следствия~\ref{cor5}
настоящей главы.

Близкие следствию~\ref{lem_phase} результаты имеются в книге
Дж.~Гукенхеймера и Ф.~Холмса \cite{guck} (пример c.~250-251), но в
них, во-первых, рассматривается не функция Малкина, а функция
Мельникова (см. замечание~2.2 по поводу соответствующего
определения), во-вторых, предполагается непрерывная
дифференцируемость входящей в правую часть системы (\ref{ps_2})
функции $g.$ Вычисление функций Малкина для широкого класса
возмущенных систем проведено в книге И.~И.~Блехмана \cite{ble}, но
указанные в следствии~\ref{lem_phase} формулы там отсутствуют.

Отметим, что ситуация, когда периоды возмущения и порождающего
цикла несоизмеримы, изучена в монографии В.~И.~Арнольда
\cite{arno}.

Обсудим кратко публикации автора диссертации по результатам
настоящей главы. Преобразование системы, рассмотренное в
лемме~\ref{zamena}, указано в (\cite{makope}, формула~6). Основная
теорема главы (теорема~\ref{MathNach2}), связанная с вычислением
степени $d(I-Q_\varepsilon,W)$ и на которой основан геометрический
вариант решения задачи И.~Г.~Малкина (теорема~\ref{thm4}),
опубликована в \cite{mak2005}. Формулы для фазы
 $T$-периодических решений синусоидально возмущенных систем
 (следствие~\ref{lem_phase})  в случае дважды непрерывно
дифференцируемых систем получены в \cite{maknoc} (см. утверждения
с.~151--152). Рассмотрение примера~\ref{preyexa}, связанного с
моделью Хищник-Жертва, проведено автором в \cite{jap}.

\chapter{Возмущения систем, допускающих семейство $T$-периодических
решений, начальные условия которых заполняют границу некоторого
открытого множества $U\subset\mathbb{R}^n$}\label{gl1}

\setcounter{subsection}{0}

В этой главе  изучается
 задача о существовании $T$-периодических решений в
системе вида
\begin{equation}\label{ps_1}
  \dot x=f(t,x)+\varepsilon g(t,x,\varepsilon),
\end{equation}
лежащих в заданном открытом ограниченном множестве
$U\subset\mathbb{R}^n,$ граница $\partial U$ которого состоит из
начальных условий $T$-периодических решений порождающей системы
\begin{equation}\label{np_1}
  \dot x=f(t,x).
\end{equation}
Полученные результаты позволят также указать условия существования
таких $T$-периодических решений для (\ref{ps_1}), начальные
условия которых близки к границе множества $U.$ Рассмотрение этой
ситуации мотивировано задачей о рождении $T$-периодических решений
двумерной системы (\ref{ps_1}) из цикла $\widetilde{x}$ системы
(\ref{np_1}), в случае, когда последняя автономна. Цикл
$\widetilde{x}$ играет в этом случае роль границы множества
$U\subset\mathbb{R}^2$ (см. теорему~\ref{dimtwo_1} настоящей
главы). Всюду предполагается, что
$f:\mathbb{R}\times\mathbb{R}^n\to\mathbb{R}^n$ -- непрерывно
дифференцируемая и
$g:\mathbb{R}\times\mathbb{R}^n\times[0,1]\to\mathbb{R}^n$ --
непрерывная $T$-периодические по первой переменной функции.

Для решения поставленной задачи рассматривается вопрос о
вычислении топологической степени $d(I-Q_\varepsilon,W_U)$
эквивалентного интегрального оператора
\[
(Q_\varepsilon x)(t)=x(T)+\int\limits_0^t
f(\tau,x(\tau))d\tau+\varepsilon\int\limits_0^t
g(\tau,x(\tau),\varepsilon)d\tau, \quad t\in[0,T],
\]
на множестве
\[
W_U=\left\{\widehat {x}\in C([0,T],\mathbb{R}^n):\widehat{
x}(t)\in U\ \mbox{\rm для\ любого\ } t\in [0,T]\right\}.
\]

\section{Формула для вычисления топологической степени
интегрального оператора, эквивалентного задаче о $T$-периодических
решениях с начальными условиями в $U$}

Распространенным инструментом изучения систем вида (\ref{ps_1})
является следующая вспомогательная система (см. \cite{mal},
формула~3.8)
\begin{equation}\label{vsp_1}
  \dot y=f'_{(2)}(t,\Omega(t,0,\xi))y+g(t,\Omega(t,0,\xi),0),
\end{equation}
где $\Omega$ -- оператор сдвига по траекториям системы
(\ref{np_1}) и  $\xi\in\mathbb{R}^n$ -- параметр.

\begin{opr}
Функция $\Phi^s:\mathbb{R}^n\to\mathbb{R}^n,$ заданная для каждого
$s\in[0,T]$ как
$$
  \Phi^s(\xi)=\eta(T,s,\xi)-\eta(0,s,\xi),
$$
где $\eta(\cdot,s,\xi)$ -- решение системы (\ref{vsp_1}),
удовлетворяющее начальному условию $\eta(s,s,\xi)=0,$ называется
обобщенным оператором усреднения по отношению к $T$-периодической
системе (\ref{ps_1}).
\end{opr}

Пусть $\Phi^s$ -- обобщенный оператор усреднения, соответствующий
системе (\ref{ps_1}) и $\mathcal{P}_0$ -- оператор Пуанкаре
системы (\ref{np_1}). Нам понадобится следующее свойство.

\begin{lem}\label{form_eta}  Для любых $\xi\in\mathbb{R}^n$ и $s,\theta\in[0,T]$ справедлива
формула
$$ \eta(\theta,s,\xi)=\Omega'_{(3)}(\theta,0,\xi)\int\limits_s^\theta
\Omega'_{(3)}(0,\tau,\Omega(\tau,0,\xi))g(\tau,\Omega(\tau,0,\xi),0)
d\tau.
$$
Если дополнительно известно, что $\mathcal{P}_0(\xi)=\xi,$ то
$$
\Phi^s(\xi)=\int\limits_{s-T}^s
\Omega'_{(3)}(0,\tau,\Omega(\tau,0,\xi))g(\tau,\Omega(\tau,0,\xi),0)
d\tau
$$
при всех $s\in[0,T].$
\end{lem}

\noindent {\bf Доказательство.} \ Заметим (см. \cite{kraope},
теорема 2.1), что матрица $\Omega'_{(3)}(t,0,\xi)$ является
фундаментальной матрицей для линейной системы
$$
  \dot y =f'_{(2)}(t,\Omega(t,0,\xi))y,
$$
причем
$\left(\Omega'_{(3)}(t,0,\xi)\right)^{-1}=\Omega'_{(3)}(0,t,\Omega(t,0,\xi)).$
Действительно, дифференцируя по $\xi$ тождество
$$
  \Omega(0,t,\Omega(t,0,\xi))=\xi,\ \ \xi\in\mathbb{R}^n,
$$
получаем
\begin{equation}\label{f12}
  \Omega'_{(3)}(0,t,\Omega(t,0,\xi))\Omega(t,0,\xi)=I,\ \ \xi\in\mathbb{R}^n.
\end{equation}
Следовательно, по формуле вариации произвольной постоянной для
неоднородной системы (\ref{vsp_1}) получаем
\begin{eqnarray}
  \eta(t,s,\xi) & = & \int\limits_s^t \Omega'_{(3)}(t,0,\xi)\left(\Omega(\tau,0,\xi)\right)^{-1}g(\tau,\Omega(\tau,0,\xi),0)d\tau = \nonumber\\
      & = & \Omega'_{(3)}(t,0,\xi)\int\limits_s^t \Omega'_{(3)}(0,\tau,\Omega(\tau,0,\xi))
      g(\tau,\Omega(\tau,0,\xi),0)d\tau.\nonumber
\end{eqnarray}
Первая формула леммы доказана, переходим к доказательству второй
формулы. Имеем
$$
  \Phi^s(\xi)=\Omega'_{(3)}(T,0,\xi)\int\limits_s^T \left(\Omega(\tau,0,\xi)\right)^{-1}
      g(\tau,\Omega(\tau,0,\xi),0)d\tau-
$$
\begin{equation}\label{QQ}
  -\int\limits_s^0 \left(\Omega(\tau,0,\xi)\right)^{-1}
      g(\tau,\Omega(\tau,0,\xi),0)d\tau.
\end{equation}
 Производя замену переменных $\tau=u+T$ в интеграле
$$J=\Omega'_{(3)}(T,0,\xi)\int\limits_s^T \left(\Omega(\tau,0,\xi)\right)^{-1}
      g(\tau,\Omega(\tau,0,\xi),0)d\tau,
$$
и учитывая, что $\mathcal{P}_0(\xi)=\xi,$ получаем
$$
  J=\Omega'_{(3)}(T,0,\xi)\int\limits_{s-T}^0 \left(\Omega(T+u,0,\xi)\right)^{-1}
      g(T+u,\Omega(T+u,0,\xi),0)du=
$$
$$
=\Phi(T){\rm e}^{\Lambda T}\int\limits_{s-T}^0
\left(\Omega(T+u,0,\xi)\right)^{-1}
      g(u,\Omega(u,0,\xi),0)d\tau=
$$
$$
=\Phi(T)\int\limits_{s-T}^0 {\rm e}^{-\Lambda u}{\rm e}^{\Lambda(
T+u)} \left(\Omega(T+u,0,\xi)\right)^{-1}
      g(u,\Omega(u,0,\xi),0)d\tau,
$$
где $\Phi(t){\rm e}^{\Lambda t}$ -- представление Флоке
нормированной при $t=0$ фундаментальной матрицы
$\Omega'_{(3)}(t,0,\xi)$ (см. \cite{dem}, Гл.~III, \S~15). Так как
${\rm e}^{\Lambda
t}\left(\Omega(t,0,\xi)\right)^{-1}=\Phi^{-1}(t),$ то
$$
  {\rm e}^{\Lambda(T+u)}\left(\Omega(T+u,0,\xi)\right)^{-1}
={\rm e}^{\Lambda u}\left(\Omega(u,0,\xi)\right)^{-1}
$$
и
$$
J=\Phi(0)\int\limits_{s-T}^0 {\rm e}^{-\Lambda u}{\rm e}^{\Lambda
u} \left(\Omega(u,0,\xi)\right)^{-1}
      g(u,\Omega(u,0,\xi),0)d\tau=
$$
$$
=\int\limits_{s-T}^0 \left(\Omega(u,0,\xi)\right)^{-1}
      g(u,\Omega(u,0,\xi),0)d\tau.
$$
Подставляя полученное выражение в (\ref{QQ}), получаем
$$
  \Phi^s(\xi)=\int\limits_{s-T}^0 \left(\Omega(\tau,0,\xi)\right)^{-1}
      g(\tau,\Omega(\tau,0,\xi),0)d\tau+
$$
$$
  +\int\limits_0^s \left(\Omega(\tau,0,\xi)\right)^{-1}
      g(\tau,\Omega(\tau,0,\xi),0)d\tau=
$$
$$
  =\int\limits_{s-T}^s \left(\Omega(\tau,0,\xi)\right)^{-1}
      g(\tau,\Omega(\tau,0,\xi),0)d\tau.
$$

Лемма доказана.

Имеет место следующий результат.

\begin{thm}\label{fromnon} Пусть $\mathcal{P}_0(\xi)=\xi$ для любого
$\xi\in\partial U.$  Если

 $\Phi^s(\xi)\not=0$ для любого $\xi\in\partial U$ и любого
$s\in[0,T],$

\noindent то существует $\varepsilon_0>0$ такое, что при
$\varepsilon\in(0,\varepsilon_0]$ справедливы следующие
утверждения

1) для любого $x\in C([0,T],\mathbb{R}^n)$ такого, что
$Q_\varepsilon x=x$ имеем $x(t)\not\in\partial U$ для всех
$t\in[0,T],$ в частности, оператор $Q_\varepsilon$ не имеет
неподвижных точек на $\partial W_U;$

2) $
  d(I-Q_\varepsilon,W_U)=d_{\mathbb{R}^n}(-\Phi^T,U).
$
\end{thm}

\noindent {\bf Доказательство.} Положим
$$\Upsilon_\varepsilon(t,\xi)=
\Omega'_{(3)}(0,t,\Omega(t,0,\xi))g(t,\Omega(t,0,\xi),\varepsilon).
$$
Согласно лемме~\ref{zamena} каждой неподвижной точке из $W_U$
оператора $Q_\varepsilon$ соответствует неподвижная точка
(\ref{zp_}) оператора
$$
(G_\varepsilon
\nu)(t)=\Omega(T,0,\nu(T))+\int\limits_0^t\Upsilon_\varepsilon(\tau,\nu(\tau))d\tau,
$$
которая, как легко проверить, вновь принадлежит $W_U.$
Следовательно, если $d(I-G_\varepsilon,W_U)$ определен, то в силу
теоремы об эквивалентных векторных полях (см. \cite{krazab},
теорема~26.4) имеем
$$
  d(I-Q_\varepsilon,W_U)=d(I-G_\varepsilon,W_U).
$$

В пространстве $C([0,T],\mathbb{R}^n)$ рассмотрим вспомогательный
вполне непрерывный оператор
$$
  (A_{\varepsilon}\nu)(t)=\Omega(T,0,\nu(T))-
    \varepsilon\int \limits_0^T\Upsilon_\varepsilon(\tau,\nu(\tau))d\tau
$$
и покажем, что при достаточно малых $\varepsilon>0$ поля
$I-G_\varepsilon$ и $I-A_{\varepsilon}$ гомотопны на границе
множества $W_U$. Зададим следующую деформацию
$$
D_{\varepsilon}(\lambda,\nu)(t)=\nu(t)-\Omega(T,0,\nu(T))-\varepsilon\int
\limits_0^{\lambda
t+(1-\lambda)T}\Upsilon_\varepsilon(\tau,\nu(\tau))d\tau,
$$
соединяющую поля $G_\varepsilon$ и $G_{1,\varepsilon}$. Покажем,
что при достаточно малых $\varepsilon>0$ деформация
$D_\varepsilon$ невырождена на границе множества $W_U.$ Для этого
докажем более сильное утверждение, которое используем затем для
доказательства утверждения~1). Именно,  покажем, что существует
$\varepsilon_0>0$ такое, что при всех
$\varepsilon\in(0,\varepsilon_0]$ и $\lambda\in[0,1]$ каждое
решение уравнения $D_\varepsilon(\lambda,\nu)=\nu$ удовлетворяет
условию $\nu(t)\not\in\partial U$ при всех $t\in[0,T].$
Предположим, что это не так. Тогда для произвольной
последовательности
${\left\{\varepsilon_k\right\}}_{k\in\mathbb{N}}\subset(0,1]$
такой, что $\varepsilon\to 0$ при $k\to\infty$ найдутся
последовательности
${\left\{\lambda_k\right\}}_{k\in\mathbb{N}}\subset[0,1]$ и
${\left\{\nu_k\right\}}_{k\in\mathbb{N}}\subset
C([0,T],\mathbb{R}^n),$ при которых
\begin{equation}\label{f15}
  \nu_k(t)=\Omega(T,0,\nu_k(T))+\varepsilon_k\int \limits_0^{\lambda_k t+(1-\lambda_k)T}
  \Upsilon_{\varepsilon_k}(\tau,\nu_k(\tau))d\tau,\ \ t\in[0,T]
\end{equation}
и
\begin{equation}\label{nuzhno}
  \nu_k([0,T])\cap\partial U\not=\emptyset.
\end{equation}
 Так как последовательность чисел
${\left\{\lambda_k\right\}}_{k\in\mathbb{N}}$ ограничена, то из
нее можно выделить сходящуюся подпоследовательность. Поэтому, без
ограничения общности можем считать, что последовательность
${\left\{\lambda_k\right\}}_{k\in\mathbb{N}}$ сходится. Из
(\ref{nuzhno}) следует, что функции последовательности
${\left\{\nu_k\right\}}_{k\in\mathbb{N}}$ равномерно ограничены.
Поэтому, на основании непрерывности функции
$\Upsilon_{\varepsilon_k}$ найдется константа $c_0>0$ такая, что
$\Upsilon_{\varepsilon_k}(t,\nu_k(t))\le c_0$, $t\in[0,T],$
$k\in\mathbb{N},$ и для любых $t_1,\ t_2\in[0,T],$
$k\in\mathbb{N}$ имеем оценку
$$
  \|\nu_k(t_2)-\nu_k(t_1)\|=\varepsilon_k\left\|\int
    \limits_{\lambda_k t_1+(1-\lambda_k)T}^{\lambda_k t_2+
     (1-\lambda_k T)}\Upsilon_{\varepsilon_k}(\tau,z_k(\tau))d\tau\right\|\le
       \varepsilon_k\lambda_k(t_2-t_1)c_0,
$$
из которой следует, что функции последовательности
${\left\{\nu_k\right\}}_{k\in\mathbb{N}}$ равностепенно
непрерывны. Значит, применяя теорему Арцела, из этой
последовательности можно выделить сходящуюся
подпоследовательность. Поэтому мы без ограничения общности можем
считать, что последовательность
${\left\{\nu_k\right\}}_{k\in\mathbb{N}}$ сходится. Положим
$\lambda_0=\lim_{k\to\infty}\lambda_k$ и
$\nu_0=\lim_{k\to\infty}\nu_k.$ Тогда $\lambda_0\in[0,1]$ и
$\nu_0([0,T])\cap\partial U\not=\emptyset.$ Так как $\dot \nu_k\to
0$ при $k\to\infty,$ то функция $\nu_0$ постоянна. Соотношение
(\ref{nuzhno}) эквивалентно существованию числа $t_k\in[0,T]$
такого, что $\nu_k(t_k)\in\partial U.$ Тогда
\begin{equation}\label{f21}
 \Omega(T,0,\nu_k(t_k))=\nu_k(t_k), \quad k\in\mathbb{N}.
\end{equation}
Вычитая из равенства (\ref {f15}), записанного при $t=T$, это же
равенство, записанное при $t=t_n$, получаем
\begin{equation}\label{f22}
  \nu_k(T)-\nu_k(t_k)=\varepsilon_k\int \limits_{\lambda_k
  t_k+(1-\lambda_k)T}^T\Upsilon_{\varepsilon_k}(\tau,\nu_k(\tau))d\tau.
\end{equation}
На основании (\ref{f21}) выражение (\ref{f15}) при $t=T$ можно
переписать в виде
\begin{eqnarray}
  \nu_k(T)-\nu_k(t_k)=\Omega(T,0,\nu_k(T))-\Omega(T,0,\nu_k(t_k))+\nonumber\\
  +\varepsilon_k\int \limits_0^T\Upsilon_{\varepsilon_k}(\tau,z_k(\tau))d\tau\nonumber
\end{eqnarray}
или
\begin{eqnarray}\label{f23}
  (I-\Omega'_{(3)}(T,0,\nu_k(t_k)))(\nu_k(T)-\nu_k(t_k))= \nonumber\\
  =\varepsilon_k\int \limits_0^T\Upsilon_{\varepsilon_k}(\tau,\nu_k(\tau))d\tau+o(\nu_k(t_k),\nu_k(T)-\nu_k(t_k)),
\end{eqnarray}
где функция $o(\xi,h)$ удовлетворяет соотношению
\begin{equation}
  \frac{\|o(\xi,h)\|}{\|h\|} \to 0\quad \mbox{при}\ \|h\| \to 0,\ \xi\in {\rm R}^k.
\end{equation}
Подставляя (\ref {f22}) в (\ref {f23}), получаем равенство
\begin{eqnarray}\label{f25}
  (I-\Omega'_{(3)}(T,0,\nu_k(t_k)))\int \limits_{\lambda_k t_k+(1-\lambda_k)T}^T\Upsilon_{\varepsilon_k}(\tau,\nu_k(\tau))d\tau =
  \nonumber \\
  =\int
  \limits_0^T\Upsilon_{\varepsilon_k}(\tau,\nu_k(\tau))d\tau+\frac{o(\nu_k(t_k),\nu_k(T)-\nu_k(t_k))}{\varepsilon_k}.
\end{eqnarray}
Из (\ref {f22}) следует, что найдется константа $c>0$ такая, что
$$
\|\nu_k(T)-\nu_k(t_k)\| \le c\varepsilon_k.
$$
Откуда
\begin{equation}\label{f26}
  \left\|\frac{o(\nu_k(t_k),\nu_k(T)-\nu_k(t_k))}{\varepsilon_k} \right\| \le c \frac{\|o(\nu_k(t_k),\nu_k(T)-\nu_k(t_k))\|}
  {\|\nu_k(T)-\nu_k(t_k) \|}.
\end{equation}
Из (\ref {f26}), учитывая то, что значения функций $\nu_k$
ограничены равномерно по $k\in\mathbb{N},$ следует, что
\begin{equation}\label{f27}
  \left\|\frac{o(\nu_k(t_k),\nu_k(T)-\nu_k(t_k))}{\varepsilon_k} \right\| \to 0, \ {\rm при}\ k \to \infty.
\end{equation}
Совершив, учитывая (\ref{f27}), предельный переход при $k \to
\infty$ в (\ref {f25}), получим
$$
  (I-\Omega'_{(3)}(T,0,\xi_0))\int \limits_s^T\Upsilon_{0}(\tau,\xi_0)d\tau =
  \int \limits_0^T\Upsilon_{0}(\tau,\xi_0)d\tau
$$
или
\begin{eqnarray}\label{f28}
  \hskip-1cm\int\limits_s^T\Omega'_{(3)}(T,0,\xi_0)\Omega'_{(3)}(0,\tau,\Omega(\tau,0,\xi_0))g(\tau,\Omega(\tau,0,\xi_0),0)d\tau- \nonumber\\
  -\int\limits_s^0\Omega'_{(3)}(0,\tau,\Omega(\tau,0,\xi_0))g(\tau,\Omega(\tau,0,\xi_0),0)d\tau=0,
\end{eqnarray}
где $s=\lim_{k\to\infty}\left(\lambda_k
t_k+(1-\lambda_k)T\right)\in[0,T].$ Пользуясь
леммой~\ref{form_eta}, равенство (\ref{f28}) можно переписать в
виде
$$
  \eta(T,s,\xi_0)-\eta(0,s,\xi_0)=0,
$$
в чем противоречие с условием теоремы. Таким образом, существует
$\varepsilon_0 > 0$ такое, что при всех
$\varepsilon\in(0,\varepsilon_0]$ и $\lambda\in[0,1]$ каждое
решение уравнения $D_\varepsilon(\lambda,\nu)=\nu$ удовлетворяет
условию $\nu(t)\not\in\partial U$ при всех $t\in[0,T].$ При
$\lambda=1$ полученный результат совпадает с утверждением 1)
теоремы. Перейдем к доказательству утверждения 2). Как уже
говорилось, доказанное свойство означает, в частности, что
\begin{equation}\label{f31}
  D_\varepsilon(\lambda,\nu)\not= 0,\  \lambda \in [0,1] , \ \nu \in \partial W_U,
  \ \
  \varepsilon\in(0,\varepsilon_0],
\end{equation}
то есть поля $I-G_{\varepsilon}$ и $I-A_{\varepsilon}$ гомотопны
на границе множества $W_U$ при $\varepsilon\in(0,\varepsilon_0].$
Обозначим через $C_0([0,T],\mathbb{R}^n)$ подпространство
пространства $C([0,T],\mathbb{R}^n),$ состоящее из всех постоянных
функций, определенных на отрезке $[0,T]$ и принимающих значения в
$\mathbb{R}^n.$ Имеем $A_{\varepsilon}(\partial W_U)\subset
C_0([0,T],\mathbb{R}^n).$ Далее, по построению множество $W_U$
содержит функции, тождественно равные произвольному фиксированному
элементу из $U.$ Наконец, из (\ref{f31}) при $\lambda=0$ получаем
$$
  A_{\varepsilon}(\nu)\not= \nu,\  \nu \in \partial W_U, \ \varepsilon\in(0,\varepsilon_0],
$$
откуда, учитывая соотношение $\partial (W_U\cap
C_0([0,T],\mathbb{R}^n))\subset\partial W_U,$ следует, что при
$\varepsilon\in(0,\varepsilon_0]$ поле $I-A_{\varepsilon}$ не
имеет нулей на границе множества $W_U\cap
C_0([0,T],\mathbb{R}^n).$ Поэтому, при
$\varepsilon\in(0,\varepsilon_0]$ законно сужение поля
$I-A_{\varepsilon}$ на подпространство $C_0([0,T],\mathbb{R}^n),$
что означает (см. \cite{krazab}, теорема~27.1)
$$
  {d}_{C([0,T],\mathbb{R}^n)}(I-A_{\varepsilon},W_U)=
$$
\begin{equation} \label{vr1}
=
    {d}_{C_0([0,T],\mathbb{R}^n)}(I-A_{\varepsilon},
       W_U\cap C_0([0,T],\mathbb{R}^n)),
\end{equation}
где в левой и правой частях равенства записаны топологические
степени в пространствах $C([0,T],\mathbb{R}^n)$ и
$C_0([0,T],\mathbb{R}^n)$ соответственно.

Заметим, что постоянная функция $\nu\in W_U\cap
C_0([0,T],\mathbb{R}^n)$ тогда и только тогда является решением
уравнения $A_{\varepsilon}\nu=\nu,$ когда элемент $\xi=\nu(0)$
является решением уравнения $A^0_{\varepsilon}\xi=\xi,$ где
$$
  A_{\varepsilon}^0\xi=\Omega(T,0,\xi)+\varepsilon\int
  \limits_0^{T}\Upsilon_\varepsilon(\tau,\xi)d\tau.
$$
Применяя теорему об эквивалентных уравнениях к уравнениям
$A_{\varepsilon}\nu=\nu$ и $A_{\varepsilon}^0\xi=\xi$ (см.
\cite{krazab}, теорема~26.4), при
$\varepsilon\in(0,\varepsilon_0]$ получаем
\begin{equation}\label{vr2}
    {d}_{C_0([0,T],\mathbb{R}^n)}(I-A_{\varepsilon},W_U\cap C_0([0,T],{\rm R}^n))
     = {d}_{\mathbb{R}^n}(I-A_\varepsilon^0,U).
\end{equation}
Для вычисления топологической степени ${\rm
d}_{\mathbb{R}^n}(I-A_\varepsilon^0,U)$ положим
$$
  A_{1,\varepsilon}\xi=-\varepsilon\int \limits_0^{T}\Upsilon_\varepsilon(\tau,\xi)d\tau.
$$
Из условия теоремы следует, что $
(I-A_{\varepsilon}^0)(\xi)=A_{1,\varepsilon}\xi,\ $
$\xi\in\partial U,$ поэтому при $\varepsilon\in(0,\varepsilon_0]$
имеем
\begin{equation}\label{vr3}
  {d}_{\mathbb{R}^n}(I-A_{\varepsilon}^0,U)=
    {d}_{\mathbb{R}^n}(A_{1,\varepsilon},U).
\end{equation}
Покажем, что векторные поля $A_{1,\varepsilon}$ и $A_{1,1}$
гомотопны на границе множества $U$ при
$\varepsilon\in(0,\varepsilon_0].$ Зададим линейную деформацию
$$
  D_{1,\varepsilon}(\lambda,\xi)=(\lambda\varepsilon+1-\lambda)\int \limits_0^{T}\Upsilon_\varepsilon(\tau,\xi)d\tau,\ \ \xi\in U.
$$
Покажем, что эта деформация невырождена на границе множества $U.$
Предположим, что это не так, тогда для некоторых
$\lambda\in[0,1],\ $ $\xi\in\partial U$ и
$\varepsilon\in(0,\varepsilon_0)$ будем иметь
$$
  (\lambda_0\varepsilon+1-\lambda)\int \limits_0^{T}\Upsilon_\varepsilon(\tau,\xi)d\tau=0,
$$
откуда
$$
  \int \limits_0^{T}\Upsilon_\varepsilon(\tau,\xi)d\tau=0,
$$
в чем противоречие с невырожденностью поля $A_{1,\varepsilon}$ при
$\varepsilon\in(0,\varepsilon_0]$ на $\partial U.$ Таким образом,
пользуясь леммой~\ref{form_eta},
\begin{equation}\label{vr4}
  {d}_{\mathbb{R}^n}(A_{1,\varepsilon},U)=
     {d}_{\mathbb{R}^n}(A_{1,1},U)={d}_{\mathbb{R}^n}(\eta(0,T,\cdot),U).
\end{equation}
Подставляя (\ref{vr4}) в (\ref{vr3}), получаем
\begin{equation}\label{res_}
  {d}_{\mathbb{R}^n}(I-A_\varepsilon^0,U)=
    {d}_{\mathbb{R}^n}(\eta(0,T,\cdot),U),\ \ \varepsilon\in(0,\varepsilon_0).
\end{equation}
Подставляя (\ref{vr2}) в (\ref{vr1}), пользуясь гомотопностью
полей $I-A_\varepsilon$ и $I-G_\varepsilon$ и соотношением
(\ref{res_}), окончательно имеем
$$
  {d}_{C([0,T],\mathbb{R}^n)}(I-G_\varepsilon,W_U)=
      {d}_{\mathbb{R}^n}(-\eta(T,0,\cdot),U), \ \ \varepsilon\in(0,\varepsilon_0].
$$

Теорема доказана.

Предположим теперь, что

$(A_{\mathcal{P}})$ решение $x$ системы (\ref{ps_1}) с начальным
условием $x(t_0)=\xi$ существует, единственно и продолжимо на
отрезок $[0,T]$ при любых $t_0\in[0,T],$ $\xi\in\mathbb{R}^n$ и
$\varepsilon>0.$

При условии $(A_{\mathcal{P}})$ для системы (\ref{ps_1}) при любых
$\varepsilon>0$  определен оператор $\mathcal{P}_\varepsilon$
Пуанкаре, соответствующий задаче о $T$-периодических решениях для
(\ref{ps_1}) (см. определение~1.2).

\begin{thm}\label{andr} Пусть $\mathcal{P}_0(\xi)=\xi$ для любого
$\xi\in\partial U.$ Пусть выполнено условие $(A_{\mathcal{P}}),$ и

 $\Phi^s(\xi)\not=0$ для любого $\xi\in\partial U$ и любого
$s\in[0,T].$

\noindent Тогда существует $\varepsilon_0>0$ такое, что при
$\varepsilon\in(0,\varepsilon_0]$ оператор
$\mathcal{P}_\varepsilon$ невырожден на $\partial U,$ и

$
  d_{\mathbb{R}^n}(I-\mathcal{P}_\varepsilon,U)=d_{\mathbb{R}^n}(-\Phi^T,U),\quad\varepsilon\in(0,\varepsilon_0].
$
\end{thm}

\noindent{\bf Доказательство.} Пусть $Q_\varepsilon$ --
интегральный оператор, соответствующий задаче о $T$-периодических
решениях для (\ref{ps_1}). Положим
$$
  W_\varepsilon=\left\{\hat x:
  C([0,T],\mathbb{R}^n):\Omega_\varepsilon(0,t,\hat x(t))\in
  U,\ \mbox{\rm для\ любого\ }t\in[0,T]\right\},
$$
где через $\Omega_\varepsilon$ обозначен оператор сдвига по
траекториям системы (\ref{ps_1}). Покажем, что существует
$\varepsilon_0>0$ такое, что для любого
$\varepsilon\in(0,\varepsilon_0]$ и любого
$\alpha\in[0,\varepsilon_0]$ выполнено:
\begin{equation}\label{ggg}
  \mbox{\rm если}\ Q_\varepsilon x=x\ \mbox{\rm и}\ x\in \overline{W}_\alpha,\ \mbox{\rm то}\
  x\in W_0.
\end{equation}

Предположим противное, следовательно, существуют
последовательности
$\{\varepsilon_n\}_{n\in\mathbb{N}}\subset(0,\varepsilon_0],$
$\varepsilon_n\to 0$ при $n\to\infty,$
$\{\alpha_n\}_{n\in\mathbb{N}}\subset(0,\varepsilon_0],$
$\{x_n\}_{n\in\mathbb{N}}\subset C([0,T],\mathbb{R}^n),$ $x_n\in
\overline{W}_{\alpha_n}$ такие, что $Q_{\varepsilon_n}x_n=x_n$ и
$x_n\not\in W_0.$ Так как $x_n\in \overline{W}_{\varepsilon_n},$
то $x_n(0)\in U.$ С другой стороны, из соотношения $x_n\not\in
W_0$ заключаем, что при любом $n\in\mathbb{N}$ существует
$t_n\in(0,T]$ такое, что $x_n(t_n)\in\partial U,$ в чем
противоречие с утверждением~1) теоремы~\ref{fromnon}.

 Из (\ref{ggg}) заключаем, что степень $d(I-Q_\lambda,W_\lambda)$ определена при
 любом  $\varepsilon\in
(0,\varepsilon_0],$ и
$$
  d(I-Q_\varepsilon,W_\varepsilon)=d(I-Q_\varepsilon,W_0),\quad\varepsilon\in(0,\varepsilon_0].
$$
Из принципа родственности (см. \cite{krazab}, теорема~28.5)
следует, что
$$
  d(I-Q_\varepsilon,W_\varepsilon)=d_{\mathbb{R}^n}(I-\mathcal{P}_\varepsilon,U).
$$
С другой стороны, в силу утверждения~2) теоремы~\ref{fromnon}
имеем
$$
  d(I-Q_\varepsilon,W_0)=d_{\mathbb{R}^n}(-\Phi^T,U).
$$

Теорема доказана.

\section{Теоремы о продолжении $T$-периодических решений из
$\overline{U}$ по параметру}

В этом пункте будут указаны приложения теоремы~\ref{fromnon} к
задаче о продолжении лежащих в $\overline{U}$ $T$-периодических
решений системы (\ref{ps_1}) при увеличении параметра
$\varepsilon$ от нуля до некоторого достаточно малого
положительного числа.

\subsubsection{Случай систем любого числа уравнений}

\begin{thm}\label{fromdan}
Если
 $\Phi^s(\xi)\not=0$ для любого $\xi\in\partial U$ и любого
$s\in[0,T],$
 и
$$
  d_{\mathbb{R}^n}(-\Phi^T,U)\not=0,
$$
то при всех достаточно малых $\varepsilon>0$ система (\ref{ps_1})
имеет по крайней мере одно $T$-периодическое решение с начальным
условием в $U.$
\end{thm}

Предположим теперь, что выполнено условие $(A_{\mathcal{P}})$ (см.
теорему~\ref{andr}), и пусть $\mathcal{P}_\varepsilon$ -- оператор
Пуанкаре,
 соответствующий задаче о
$T$-периодических решениях для возмущенной системы  (\ref{ps_1}).
Предположим далее, что

$(A_*)$ отображение $\mathcal{P}_0$ имеет конечное число
$a_1,...,a_k$ неподвижных точек в $U,$ при этом для любого
$i\in\overline{1,k}$ матрица
$\left(I-(\mathcal{P}_0)'\right)(a_k)$ невырождена.

\noindent При условии $(A_*)$ каждому $a_i,$ $i\in\overline{1,k}$
можно поставить в соответствие индекс Пуанкаре ${\rm
ind}(a_i,I-\mathcal{P}_0)$ (см. \cite{lef}, Гл.~IX, \S~4), который
в этом случае определяется как $(-1)^\beta,$ где $\beta$ -- сумма
кратностей вещественных и отрицательных собственных значений
матрицы $\left(I-(\mathcal{P}_0)'\right)(a_k)$ (см. \cite{krazab},
теорема~6.1).

Справедлив следующий результат.

\begin{thm}\label{ifac}
Предположим, что выполнены условия теоремы~\ref{andr}, а также
предположение~$(A_*).$  Предположим, что
\begin{equation}\label{as3_}
  {\rm ind}(a_1,I-\mathcal{P}_0)+\cdot\cdot+{\rm ind}(a_k,I-\mathcal{P}_0)\not=
  d_{\mathbb{R}^n}(-\Phi^T,U).
\end{equation}
Тогда при всех достаточно малых $\varepsilon>0$ система
(\ref{ps_1}) имеет $T$-периодическое решение $x_\varepsilon(t)$
такое, что
\begin{equation}\label{llk}
  \rho(x_\varepsilon(0),\partial U)\to 0\quad\mbox{при}\quad \varepsilon\to
  0.
\end{equation}
\end{thm}

\noindent {\bf Доказательство.} Из $(A_*)$ следует, что точки
$a_1,...,a_k$ являются изолированными. Поэтому (см. \cite{kraope},
теорема~6.1),  существует $\delta_0>0$ такое, что
$$
d_{\mathbb{R}^n}(I-\mathcal{P}_0,B_{\delta_0}(a_i))={\rm
ind}(a_i,I-\mathcal{P}_0),\quad i\in\overline{1,k}.
$$
Далее, при условии $(A_*)$ в силу теоремы Пуанкаре (см.
\cite{malb}, утверждение~с.~378) $\delta_0>0$ может быть выбрано
еще настолько малым, что для некоторого $\varepsilon_0>0$ и любых
$\varepsilon\in(0,\varepsilon_0],$ $i\in\overline{1,k}$ существует
единственное $T$-периодическое решение системы (\ref{ps_1}) в
$\delta_0$-окрестности точки $a_i.$ Из теоремы~\ref{andr} следует,
что $\varepsilon_0>0$ может быть выбрано настолько малым, что
$$
d_{\mathbb{R}^n}(I-\mathcal{P}_\varepsilon,U)=d_{\mathbb{R}^n}(-\Phi^T,U),\quad
\varepsilon\in(0,\varepsilon_0].
$$
Наконец, уменьшим $\varepsilon_0>0$ еще и так, что
$$
d_{\mathbb{R}^n}(I-\mathcal{P}_0,B_{\delta_0}(a_i))=
d_{\mathbb{R}^n}(I-\mathcal{P}_\varepsilon,B_{\delta_0}(a_i)),\quad
i\in\overline{1,k},\ \varepsilon\in(0,\varepsilon_0].
$$
Суммируя полученные соотношения и учитывая (\ref{as3_}), приходим
к выводу, что
$$
  d_{\mathbb{R}^n}\left(I-\mathcal{P}_\varepsilon,U\backslash\bigcup\limits_{i=1}^k
  B_{\delta_0}(a_i)\right)\not=0,\quad\varepsilon\in(0,\varepsilon_0].
$$
Следовательно, при любом $\varepsilon\in(0,\varepsilon_0]$
существует решение $x_\varepsilon$ системы (\ref{ps_1}) такое, что
$$
  x_\varepsilon(0)\in U\backslash\bigcup\limits_{i=1}^k
  B_{\delta_0}(a_i),\quad \varepsilon\in(0,\varepsilon_0].
$$
Заметим, что
\begin{equation}\label{plm}
  \rho(x_\varepsilon(0),\partial U)\to 0\quad\mbox{при}\ \varepsilon\to 0.
\end{equation} Действительно, в противном случае имели бы
противоречие с единственностью $T$-периодических решений с
рожденными из точек $a_1,...,a_k$ начальными условиями.

Теорема доказана.

\

\subsubsection{Случай двумерных систем, допускающих цикл в отсутствии возмущения}

Следующее приложение теоремы~2.1 связано со случаем, когда система
(\ref{np_1}) двумерна и автономна, то есть имеет вид
\begin{equation}\label{npa_1}
 \dot x=f(x),\qquad x\in\mathbb{R}^2,
\end{equation}
соответствующая возмущенная система (\ref{ps_1}) запишется в этом
случае как
\begin{equation}\label{psa_1}
 \dot x=f(x)+\varepsilon g(t,x,\varepsilon),\qquad
 x\in\mathbb{R}^2.
\end{equation}
 Предположим, что система (\ref{npa_1}) имеет
$T$-периодический цикл $\widetilde{x},$ и обозначим через $U$
внутренность цикла. Очевидно, $\Omega(T,0,\xi)=\xi$ для любого
$\xi\in\partial U.$

\begin{thm}\label{dimtwo_1}
Предположим, что

$\Phi^s(\xi)\not=0$ для любого $\xi\in\partial U$ и любого
$s\in[0,T],$

$(A_1)$ $\widetilde{x}$ -- единственный $T$-периодический цикл
системы (\ref{npa_1}) в некоторой окрестности $\partial U.$

\noindent Тогда, если

$d_{\mathbb{R}^2}(-\Phi^T,U)\not=1,$

\noindent то существует $\varepsilon_0>0$ такое, что при всех
$\varepsilon\in(0,\varepsilon_0]$ система (\ref{psa_1}) имеет по
крайней мере два $T$-периодических решения $x_{\varepsilon,1}$ и
$x_{\varepsilon_2}$ таких, что
$$
  \rho(x_{\varepsilon,1}(t),\partial U)+\rho(x_{\varepsilon,2}(t),\partial
  U)\to 0\quad\mbox{при}\ \varepsilon\to 0.
$$
Более того, $x_{\varepsilon,1}(t)\in U,$
$x_{\varepsilon,2}(t)\not\in U$ при всех $t\in[0,T],$ и каждое
прочее $T$-периодическое решение $x$ системы (\ref{psa_1})
удовлетворяет условию $x(t)\not\in\partial U$ при всех
$t\in[0,T].$
\end{thm}

Условие $(A_1)$ не отрицает существования вблизи $\partial U$
циклов системы (\ref{npa_1}) с отличными от $T>0$ периодами.

\noindent {\bf Доказательство теоремы 2.5.} Пусть
$\varepsilon_0>0$ -- то, о котором говорится в
теореме~\ref{fromnon}, тогда
$$
  d(Q_\varepsilon,W_U)=d_{\mathbb{R}^2}(-\Phi^T,U)\quad \mbox{для любых
  }\varepsilon\in(0,\varepsilon_0],
$$
или, учитывая условие $d_{\mathbb{R}^2}(-\Phi^T,U)\not=1,$
\begin{equation}\label{dee1}
  d(Q_\varepsilon,W_U)\not=1\quad \mbox{для любых
  }\varepsilon\in(0,\varepsilon_0].
\end{equation}
 Положим $U_\delta^-=U\backslash B_\delta(\partial
U),$ $U_\delta^+=U\cup B_\delta(\partial U).$ На основании условия
$(A_1)$ можно зафиксировать такое $\delta_0>0,$ что система
(\ref{npa_1}) не имеет $T$-периодических решений с начальными
условиями из $\partial U_\delta^-\cup\partial U_\delta^+$ при всех
$\delta\in(0,\delta_0].$ Без ограничения общности можем считать,
что $\delta_0>0$ выбрано достаточно малым так, что
$U_{\delta_0}^-\not=\emptyset.$ По теореме Капетто-Мавена-Занолина
(\cite{camaza}, cледствие~1) имеем
$$
   d(Q_0,W_{U_\delta^-})=d_{\mathbb{R}^2}(f,U_\delta^-)\mbox{ и }
   d(Q_0,W_{U_\delta^+})=d_{\mathbb{R}^2}(f,U_\delta^+),\quad\delta\in(0,\delta_0].
$$
Без ограничения общности можно считать, что малость $\delta_0>0$
достаточна для того, чтобы
$$
  d_{\mathbb{R}^2}(f,U_\delta^-)=d_{\mathbb{R}^2}(f,U_\delta^+)=d_{\mathbb{R}^2}(f,U),\quad\delta\in(0,\delta_0].
$$
По теореме Пункаре (см. С.~Лефшец \cite{lef}, теорема~11.1 или
М.~А.~Красносельский и др.  \cite{krapla}, теорема~2.3) имеем
$d_{\mathbb{R}^2}(f,U)=1,$ поэтому
$$
   d(Q_0,W_{U_\delta^-})=
   d(Q_0,W_{U_\delta^+})=1\quad\mbox{при
   всех }\delta\in(0,\delta_0].
$$
Таким образом, каждому $\delta\in(0,\delta_0]$ соответствует
$\varepsilon_\delta>0$ такое, что
\begin{equation}\label{dee2}
   d(Q_\varepsilon,W_{U_\delta^-})=
   d(Q_\varepsilon,W_{U_\delta^+})=1,\qquad\varepsilon\in(0,\varepsilon_\delta],
    \ \delta\in(0,\delta_0].
\end{equation}
Без ограничения общности можно считать, что
$\varepsilon_\delta<\varepsilon_0$ при всех
$\delta\in(0,\delta_0].$ Тогда из (\ref{dee1}) и (\ref{dee2})
получаем, что при всех $\delta\in(0,\delta_0]$ и
$\varepsilon\in(0,\varepsilon_\delta]$ система (\ref{psa_1}) имеет
по крайней мере два $T$-периодических решения
$x_{\varepsilon,1}\in W_U\backslash W_{U_\delta^-}$ и
$x_{\varepsilon,2}(t)\in W_{U_\delta^+}\backslash W_U.$ Из этого,
в частности, имеем $x_{\varepsilon,1}(0)\in U,$
$x_{\varepsilon,2}(t)(0)\not\in U$ и, используя утверждение~1)
теоремы~\ref{fromnon}, заключаем, что $x_{\varepsilon,1}(t)\in U,$
$x_{\varepsilon,2}(t)(t)\not\in U$ при всех $t\in[0,T].$

 Теорема доказана.

\subsubsection{Случай двумерных систем, допускающих семейство циклов в отсутствии возмущения}

Обозначим через $\widetilde {x}$ периодическое решение системы
(\ref{npa_1}) наименьшего периода ${T}.$ В этом подпункте
предполагается, что

$(C)$ алгебраическая кратность мультипликатора $+1$ системы
\begin{equation}\label{lsg}
  \dot y=f'(\widetilde {x}(t))y
\end{equation}
равна $2.$

\noindent Последнее имеет место, в частности, в случае, когда
$\widetilde{x}$ вложен в семейство циклов системы (\ref{npa_1}).
 Обозначим через
$U\subset\mathbb{R}^2$ внутренность цикла $\widetilde {x}.$ Нас
интересует вопрос о том, существуют ли вблизи границы множества
$U$ периодические решения системы (\ref{psa_1}) с периодом $T.$

\begin{lem}\label{prop_u}
 Предположим, что $T$-периодическая система
\begin{equation}\label{t1}
  \dot u =A(t)u
\end{equation}
имеет мультипликатор $+1$ алгебраической кратности $2,$ и
$\widetilde {u}$ -- $T$-периодическое решение этой системы такое,
что
$${\widetilde {u}}_1(0)\not=0,\ {\widetilde {u}}_2(0)=0\quad ({\widetilde {u}}_1(0)=0,\ {\widetilde {u}}_2(0)\not=0).$$
 Тогда для решения $\widehat {u}$ системы
(\ref{t1}), удовлетворяющего условию
$$\widehat {u}_1(0)=0,\ \widehat{
u}_2(0)\not=0\quad (\widehat {u}_1(0)\not=0,\ \widehat
{u}_2(0)=0),$$ справедлива формула
$$
\widehat {u}(t+T)=\widehat {u}(t)+\frac{\widehat{
u}_1(T)}{{\widetilde {u}}_1(0)}{\widetilde {u}}(t)\ \
  \left(\widehat {u}(t+T)=\widehat {u}(t)+\frac{\widehat {u}_2(T)}{{\widetilde {u}}_2(0)}{\widetilde{u}}(t)\right),\quad t\in\mathbb{R}.
$$
\end{lem}

\noindent {\bf Доказательство.} Установим справедливость
утверждения в случае, когда ${\widetilde{u}}_1(0)\not=0,\
{\widetilde{u}}_2(0)=0.$
 Обозначим через $X$ нормированную ($X(0)=I$) фундаментальную
матрицу системы (\ref{t1}). Так как
$X(T)\left(\begin{array}{c}1\\0\end{array}\right)=
\left(\begin{array}{c}1\\0\end{array}\right),$ то $X(T)=
\left(\begin{array}{cc}1&a\\0&b\end{array}\right).$ По условию
леммы алгебраическая кратность собственного значения $+1$ матрицы
$X(T)$ равна двум, значит  $X(T)=
\left(\begin{array}{cc}1&a\\0&1\end{array}\right),$ где
$a\in\mathbb{R}$ -- некоторое число. Имеем
\begin{eqnarray}
  X(t+T)\widehat {u}(0)&=&X(t)X(T)\widehat{
  u}(0)=X(t)\left(\begin{array}{cc}1&a\\0&1\end{array}\right)\widehat{
  u}(0)=\nonumber\\
  &=&X(t)\widehat {u}(0)+X(t)\left(\begin{array}{c}a\widehat{
  u}_2(0)\\0\end{array}\right)=\nonumber\\
  &=&X(t)\widehat {u}(0)+\frac{a\widehat {u}_2(0)}{{\widetilde {u}}_1(0)}
  {\widetilde{
  u}}(t).\nonumber
\end{eqnarray}
В тоже время
$$
  X(T)\widehat{
  u}(0)=\left(\begin{array}{cc}1&a\\0&1\end{array}\right)\widehat{
  u}(0)=\widehat {u}(0)+\left(\begin{array}{c}a\widehat{
  u}_2(0)\\0\end{array}\right),
$$
откуда $a\widehat {u}_2(0)=\widehat {u}_1(T).$

Справедливость утверждения леммы в случае, когда ${\widetilde{
u}}_1(0)=0,\ {\widetilde{u}}_2(0)\not=0$ устанавливается
аналогично.

Лемма доказана.

Не ограничивая общности решения поставленной задачи, можем
считать, что
\begin{equation}\label{init}
  \dot {{\widetilde {x}}}_1(0)\not=0\mbox{ и }\dot {{\widetilde {x}}}_2(0)=0.
\end{equation}

Пусть $\widehat {y}$ -- решение системы (\ref{lsg}),
удовлетворяющее условию
\begin{equation}\label{inity}
\widehat {y}_1(0)=0,\ \widehat {y}_2(0)\not=0.
\end{equation} Обозначим через
$\widetilde {z}$ и $\widehat {z}$ решения сопряженной системы
\begin{equation}\label{ssg}
  \dot z=(f'(\widetilde{x}(t)))^*z,
\end{equation}
удовлетворяющие начальным условиям
\begin{equation}\label{initz}
  \widehat {z}(0)=\left(\begin{array}{c}1/\dot{\widetilde{
  x}}_1(0)\\0\end{array}\right)\mbox{\ \ \ и \ \ \ }
\widetilde {z}(0)=\left(\begin{array}{c}0\\1/\widehat{
  y}_2(0)\end{array}\right).
\end{equation}
В силу леммы Перрона (см. \cite{perron} или \cite{dem}, Гл.~III,
\S~12) имеем
\begin{equation}\label{perr}
  \left(\dot{\widetilde {x}}(t)\ \widehat {y}(t)\right)^*\left(\widehat {z}(t)\
  \widetilde {z}(t)\right)=I, \quad t\in\mathbb{R}
\end{equation}

\begin{lem}\label{prop_z} Пусть выполнено условие $(C).$
Тогда решение $\widetilde {z}$ является ${T}$-периодическим.
\end{lem}

\noindent {\bf Доказательство.} Если $\widehat {y}_1(T)=0,$ то в
силу леммы~\ref{prop_u} каждое решение системы (\ref{lsg}), а
значит и системы (\ref{ssg}), является $T$-периодическим.
Рассмотрим случай, когда
\begin{equation}\label{ca}
\widehat{y}_1(T)\not=0.
\end{equation} В силу теоремы о периодических
решениях сопряженной системы (см. \cite{dem}, Гл.~III, \S~23,
теорема~2) система (\ref{lsg}) имеет по крайней мере одно
$T$-периодическое решение, обозначим это решение через
$\widetilde{z}.$ На основании леммы~\ref{prop_u} имеем
$$
 \left<\widehat {y}(T),\widetilde{\widetilde {z}}(T)\right>=
 \left<\widehat {y}(0),\widetilde{\widetilde
 {z}}(T)\right>+\frac{\widehat{
 y}_1(T)}{\dot{\widetilde {x}}_1(0)}\left<\dot{\widetilde
 {x}}(0),\widetilde{\widetilde{
 z}}(T)\right>=$$
 $$=\left<\widehat {y}(0),\widetilde{\widetilde {z}}(0)\right>+
\frac{\widehat{
 y}_1(T)}{\dot{\widetilde {x}}_1(0)} \dot{\widetilde {x}}_1(0)
 \widetilde{\widetilde{
 z}}_1(0).
$$
Но, в силу леммы Перрона (см. \cite{perron} или \cite{dem},
Гл.~III, \S~12) $ \left<\widehat{ y}(T),\widetilde{\widetilde{
z}}(T)\right>=\left<\widehat{ y}(0),\widetilde{\widetilde
{z}}(0)\right>,$ поэтому $\widetilde{\widetilde{
 z}}_1(0)=0$ и, следовательно, решения $\widetilde{
 z}$ и $\widetilde{\widetilde{
 z}}$ линейно зависимы.

Лемма доказана.

Нижеследующая лемма дает разложение поля $\Phi^s(\xi)$ по $\dot
{\widetilde {x}}$ и ${\widehat{y}}$ для случая, когда цикл
$\widetilde {x}$ удовлетворяет условию $(C).$

\begin{lem}\label{lem_Ftheta} Пусть выполнено условие $(C).$
Тогда для любых $s,\theta\in\mathbb{R}$ имеет место формула
\begin{equation}\label{Ftheta}
\Phi^s(\widetilde {x}(\theta))=\left(\widehat{f}(\theta)-
\frac{\widehat {z}_2(T)}{\widetilde{
z}_2(0)}\widetilde{f}(\theta,s+\theta)\right)\dot{\widetilde{
x}}(\theta)+ \widetilde{f}(\theta,0)\widehat {y}(\theta),
\end{equation} где
$$
   \widetilde{f}(\theta,t)=\int\limits_t^{T}
\left<\widetilde{z} (\tau),g(\tau-\theta,\widetilde{
  x}(\tau),0)\right>d\tau,
  $$
  $$
   \widehat{f}(\theta)=\int\limits_0^{T}
\left<\widehat{z} (\tau),g(\tau-\theta,\widetilde{
  x}(\tau),0)\right>d\tau.
$$
\end{lem}

\noindent {\bf Доказательство.} Положим $\widetilde
{Z}(t)=(\widehat {z}(t),\widetilde {z}(t))$ и обозначим через
$Z(t)$ фундаметальную матрицу системы (\ref{ssg}) такую, что
$Z(0)=I,$ имеем $Z(t)=\widetilde {Z}(t)\widetilde {Z}^{-1}(0).$ По
лемме Перрона (см. \cite{perron} или \cite{dem}, Гл.~III, \S~12)
$Y^{-1}(t)=Z^*(t),$ и, учитывая лемму~\ref{form_eta},
\begin{eqnarray}
\Phi^s(\widetilde {x}(\theta))&=&\int\limits_{s-T}^s
  (\Omega)'_{(3)}(0,\tau,\Omega(\tau,0,\widetilde
  {x}(\theta)))g(\tau,\widetilde{
  x}(\tau+\theta),0)d\tau=\nonumber\\
  &=&
Y(\theta)\int\limits_{s-T+\theta}^{s+\theta}
  Y^{-1}(\tau)g(\tau-\theta,\widetilde{
  x}(\tau),0)d\tau=\nonumber\\
& =& (\widetilde{
Z}^*(\theta))^{-1}\int\limits_{s-T+\theta}^{s+\theta}
  \widetilde {Z}^*(\tau)g(\tau-\theta,\widetilde{
  x}(\tau),0)d\tau= \nonumber\\
&=& (\widetilde {x}(t)\ \widetilde {y}(t)) \left(\begin{array}{c}
\int\limits_{s-T+\theta}^{s+\theta}
  \left<\widehat {z}(\tau),g(\tau-\theta,\widetilde {x}(\tau),0)\right>d\tau\\
  \widetilde{ f}(\theta)
\end{array}\right).\nonumber
\end{eqnarray}
Но, используя замену переменных $t=\tau+T$ в интеграле и
лемму~\ref{prop_u}, можем провести следующее преобразование
$$
\int\limits_{s-T+\theta}^{s+\theta}
  \left<\widehat {z}(\tau),g(\tau-\theta,\widetilde
  {x}(\tau),0)\right>d\tau=$$
$$=
  \int\limits_0^{s+\theta}\left<\widehat {z}(\tau),g(\tau-\theta,\widetilde{
  x}(\tau),0)\right>d\tau+\int\limits_{s-T+\theta}^0\left<\widehat{
  z}(\tau),g(\tau-\theta,\widetilde{
  x}(\tau),0)\right>d\tau=
$$
$$
= \int\limits_0^{s+\theta}\left<\widehat{
z}(\tau),g(\tau-\theta,\widetilde{
  x}(\tau),0)\right>d\tau+\int\limits_{s+\theta}^{T}\left<\widehat{
  z}(t-T),g(t-\theta,\widetilde{
  x}(t),0)\right>dt=
$$
$$
= \int\limits_0^{s+\theta}\left<\widehat{
z}(\tau),g(\tau-\theta,\widetilde{
  x}(\tau),0)\right>d\tau+$$
$$+\int\limits_{s+\theta}^{T}\left<\left(\widehat {z}(t)-\frac{\widehat
{z}_2(T)}{\widetilde {z}_2(0)}\widetilde
{z}(t)\right),g(t-\theta,\widetilde
  {x}(t),0)\right>dt=
$$
$$
= \int\limits_0^{T}\left<\widehat
{z}(\tau),g(\tau-\theta,\widetilde{
  x}(\tau),0)\right>d\tau-$$
  $$-
  \frac{\widehat {z}_2(T)}{\widetilde {z}_2(0)}\int\limits_{s+\theta}^{T}\left<\widetilde
  {z}(t),g(t-\theta,\widetilde{
  x}(t),0)\right>dt=\widehat{f}(\theta)-
  \frac{\widehat {z}_2(T)}{\widetilde {z}_2(0)}\widetilde{f}(\theta,s+\theta).
$$

Лемма доказана.

Лемма \ref{lem_Ftheta} позволяет получить следующее следствие из
теоремы~\ref{dimtwo_1}, в котором предполагается, что цикл
$\widetilde{x}$ удовлетворяет условию (\ref{init}).

\begin{thm}\label{dimtwoF_1}
Пусть выполнены условия $(A_1)$ и $(C).$ Предположим, что для
каждого $\theta_0\in[0,{T}]$ такого, что
$\widetilde{f}(\theta_0,0)=0,$ имеем $$
  \left|\widehat{f}(\theta_0)\right|>\left|\frac{\widehat{z}_2(T)}{\widetilde{z}_2(0)}\widetilde{f}(\theta_0,s+\theta_0)\right|
  \quad\mbox{при всех }s\in[0,T].
$$
Тогда при каждом достаточно малом $\varepsilon>0$ всякое
$T$-периодическое решение $x$ системы (\ref{psa_1}) необходимо
таково, что $x(t)\not\in\partial U$ при любом
$t\in\left[0,T\right].$ Более того, если дополнительно известно,
что

$ (A_2)\quad  d_{\mathbb{R}^2}(-\Phi^{T},U)\not=1, $

\noindent то при всех достаточно малых $\varepsilon>0$ система
(\ref{psa_1}) имеет по крайней мере два $T$-периодических решения
$\widetilde{x}_\varepsilon$ и
$\widetilde{\widetilde{x}}_\varepsilon$ таких, что
$\widetilde{x}_\varepsilon(t)\in U,$
$\widetilde{\widetilde{x}}_\varepsilon(t)\not\in U$ для любого
$t\in[0,T],$ и
$$\rho\left(\widetilde{x}_\varepsilon(t),\partial U\right)+
\rho\left(\widetilde{\widetilde{x}}_\varepsilon(t),\partial U
\right)\to 0\quad\mbox{при }\varepsilon\to 0
$$
равномерно по $t\in[0,T].$
\end{thm}


В следующих двух пунктах даются приложения разработанных теорем к
конкретным классам систем (\ref{psa_1}).

\section{Модификация теоремы Борсука-Улама и новые свойства периодических решений уравнения Дуффинга }

\begin{exa}\label{duf_exa} \rm В качестве примера рассмотрим
задачу о существовании периодических решений у уравнения Дуффинга
\begin{equation}\label{so}
  \ddot u+u+u^3=\varepsilon\cos((1+\delta)t).
\end{equation}
\end{exa}

Как известно (см., например, \cite{guck}, пример~с.~250), при
$\varepsilon=0$ уравнение (\ref{so}) допускает семейство
периодических решений, период которых изменяется монотонно от
$2\pi$ до $0,$ когда начальное условие решения изменяется от
$(0,0)$ до $(+\infty,0).$ На основании теоремы~\ref{dimtwoF_1} мы
покажем, что в каждой "полуокрестности"\ порождающего
периодического решения имеется по крайней мере одно периодическое
решение уравнения (\ref{so}) того же периода, что и порождающее.

\begin{lem}\label{duf_lem}
Обозначим через $u_\delta$ единственное с точностью до сдвига
периодическое решение невозмущенного (при $\varepsilon=0$)
уравнения Дуффинга (\ref{so}) с наименьшим периодом
$\frac{2\pi}{1+\delta}.$ Существует $\delta_0>0,$ при котором
каждому $\delta\in[0,\delta_0]$ соответствует $\varepsilon_0>0$
такое, что:

1) при любом $\varepsilon\in(0,\varepsilon_0]$ уравнение Дуффинга
(\ref{so}) имеет по крайней мере два
$\frac{2\pi}{1+\delta}$-периодических решения
$\widetilde{u}_{\delta,\varepsilon}$ и
$\widetilde{\widetilde{u}}_{\delta,\varepsilon}$ таких, что
значения функции
$\left(\widetilde{u}_{\delta,\varepsilon},\dot{\widetilde{u}}_{\delta,\varepsilon}\right)$
лежат строго внутри области $U\subset\mathbb{R}^2,$ ограниченной
кривой
$\widetilde{x}(t)=\left(u_\delta(t),\dot{{u}}_\delta(t)\right),$ а
значения функции
$\left(\widetilde{\widetilde{u}}_{\delta,\varepsilon},\dot{\widetilde{\widetilde{u}}}_{\delta,\varepsilon}
\right)$ лежат строго снаружи $U;$

2) всякое $\frac{2\pi}{1+\delta}$-периодическое решение $u$
системы (\ref{so}) при $\varepsilon\in(0,\varepsilon_0]$ таково,
что кривая $t\to\left(u(t),\dot{{u}}(t)\right)$ не имеет точек
пересечения с кривой
$t\to\left(u_\delta(t),\dot{{u}}_\delta(t)\right);$

3) решения $\widetilde{u}_{\delta,\varepsilon}$  и
$\widetilde{\widetilde{u}}_{\delta,\varepsilon}$ удовлетворяют
условию
$$
  \widetilde{u}_{\delta,\varepsilon}(t)\to u_\delta\left(t-\widetilde{\theta}\right)
  \quad\mbox{и}\quad
  \widetilde{\widetilde{u}}_{\delta,\varepsilon}(t)\to u_\delta\left(t-\widetilde{\widetilde{\theta}}\right)
  \quad\mbox{при}\ \varepsilon\to 0
$$
для некоторых
$\widetilde{\theta},\widetilde{\widetilde{\theta}}\in\left[0,\frac{2\pi}{1+\delta}\right].$
\end{lem}
 Для
доказательства леммы~\ref{duf_lem} нам понадобятся некоторые
дополнительные утверждения. Первое из них -- модификация теоремы,
доказанной К.~Борсуком в \cite{borsuk}, предположение о
справедливости которой ранее высказал С.~Улам (см. также
\cite{krapla}, теорема~2.2).

\begin{lem}\label{lem_barsuk}
 Пусть
$\widetilde{x}:[0,T]\to\mathbb{R}^2,$
$\widetilde{x}(0)=\widetilde{x}(T)$ -- жорданова кривая,
ограничивающая множество $U\subset\mathbb{R}^2.$ Пусть
$F:\mathbb{R}^2\to\mathbb{R}^2$ -- непрерывное векторное поле
такое, что $F(\xi)\not=0$ для каждого $\xi\in\partial U.$
Предположим, что существует направляющая функция
$z:[0,T]\to\mathbb{R}^2,$ $z(0)=z(T)$ такая, что:

1) $\left<z(\theta),\dot {\widetilde{x}}(\theta)\right>\not=0$ для
каждого $\theta\in[0,T],$

2) скалярная функция
$\left<F(\widetilde{x}(\cdot)),z(\cdot)\right>$ имеет ровно два
нуля $\theta_1,\theta_2$ на интервале $[0,T)$ и строго монотонна в
этих точках,

3) ${\rm sign}\left<F(\widetilde{x}(\theta_1)),
\left(\begin{array}{c} z_2(\theta_1) \\ -z_1(\theta_1)
\end{array}\right)
\right>=$ \\ $\mbox{\ }$ \hskip3cm$=-{\rm
sign}\left<F(\widetilde{x}(\theta_2)),\left(\begin{array}{c}
z_2(\theta_2) \\ -z_1(\theta_2)
\end{array}\right)\right>.$

Тогда либо $d(F,U)=0,$ либо $d(F,U)=2.$

\end{lem}

\noindent {\bf Доказательство.} Предположим, что параметризация
кривой $\widetilde{x}$ положительна, то есть множество $U$
расположено по левую сторону относительно наблюдателя, движущегося
по $\partial U$ вместе с $\widetilde{x}(t),$ когда $t$ возрастает
от $0$ до $T$, в противном случае, мы рассмотрим противоположную
параметризацию $\widetilde
{\widetilde{x}}(\theta)=\widetilde{x}(-\theta).$ Пусть
$\Gamma_{\dot {\widetilde{x}}}:[0,T]\to\mathbb{R}$ -- некоторая
однозначная ветвь угловой функции, связанной с вектором $\dot
{\widetilde{x}}(t),$ $t\in[0,T],$ (см., например, \cite{krapla},
\S 1.2), причем такая, что $\Gamma_{\dot {\widetilde{x}}}(\theta)$
возрастает, когда вектор $\dot {\widetilde{x}}(\theta)$ вращается
против часовой стрелки. На основе $\Gamma_{\dot {\widetilde{x}}}$
сейчас будет определена угловая функция $\Gamma_{F\circ
{\widetilde{x}}}$ для вектор-функции $\theta\to
F(\widetilde{x}(\theta)).$

Без ограничения общности можем считать, что
$\theta_1,\theta_2\in(0,T),$ в противном случае мы могли бы
сдвинуть время в функциях $\widetilde{x}$ и $z.$ Также можем
предполагать, что $\left<\dot
{\widetilde{x}}(\theta),z(\theta)\right>>0$ для каждого
$\theta\in[0,T],$ иначе мы рассмотрели бы $\widetilde{
z}(\theta)=z(-\theta)$ вместо $z(\theta).$ Обозначим через
$\widehat{h_1,h_2}\in[0,2\pi)$ угол между векторами  $h_1$ и
$h_2,$ посчитанный в направлении против часовой стрелки, то есть
$\widehat{h_1,h_2}+\widehat{h_2,h_1}=2\pi.$ Пусть ${\rm
ind}(\theta_i,f)=+1$ или ${\rm ind}(\theta_i,f)=-1$ в зависимости
от того возрастает или убывает  $f$ в $\theta_i,$ $i=1,2.$ Введем
функции $\angle:\mathbb{R}^2\times\mathbb{R}^2\to[-\pi,\pi]$ и
$\widetilde{H}_{\theta_i}:[0,T]\to\{-1,0,1\}$ следующим образом
$$
  \angle(h_1,h_2)=\arccos\frac{\left<h_1,h_2\right>}{\|h_1\|\cdot\|h_2\|}{\rm
  sign}\left(\pi-\widehat{h_1,h_2}\right),
$$
$$
  \widetilde{H}_{\theta_i}(\theta)=\left\{\begin{array}{ll} 0 & \mbox{\rm при\ \
  }\theta\in[0,\theta_i), \\  {\rm ind}(\theta_i,f){\rm sign}\left<z(\theta_i)^\bot,F(\widetilde{x}(\theta_i))\right>
  & \mbox{\rm при\ \
  }\theta\in[\theta_i,T].\end{array}\right.
$$
Для каждого $\theta\in[0,T]$ определим
$$
   \Gamma_{F\circ \widetilde{x}}(\theta)=\Gamma_{\dot {\widetilde{x}}}(\theta)+\angle(\dot
   {\widetilde{x}}(\theta),z(\theta))+\angle({\rm
   sign}\left<z(\theta),F(\widetilde{x}(\theta))\right>z(\theta),F(\widetilde{x}(\theta)))+
$$
\begin{equation}\label{put}
   +\pi\widetilde{H}_{\theta_1}(\theta)+
   \pi\widetilde{H}_{\theta_2}(\theta).
\end{equation}
Можно проверить, что функция $\Gamma_{F\circ \widetilde{x}}$
непрерывна на $[0,T]$ (для этого необходимо проверить это условие
только в точках $\theta_1$ и $\theta_2$), и что $\Gamma_{F\circ
\widetilde{x}}(\theta)$ возрастает, когда вектор
$F(\widetilde{x}(\theta))$ вращается против часовой стрелки,
следовательно, $\Gamma_{F\circ \widetilde{x}}(\cdot)$ является
угловой функция вектора $F(\widetilde{x}(\theta))$ для
$\theta\in[0,T].$ По определению числа вращения двумерных
векторных полей на границах односвязных множеств (см.
\cite{krapla}, \S~1.3, формула~1.11) имеем
\begin{equation}\label{wehave}
d_{\mathbb{R}^2}(F,U)=\frac{1}{2\pi}[\Gamma_{F\circ
\widetilde{x}}(T)-\Gamma_{F\circ q}(0)].
\end{equation}
По теореме Пуанкаре о полной вариации угловой функции, связанной с
касательным вектором  $\dot {\widetilde{x}}(\theta),$ на кривой
$\widetilde{x}$ (см. \cite{krapla}, теорема~2.4) имеем
$$
  \frac{1}{2\pi}[\Gamma_{\dot {\widetilde{x}}}(T)-\Gamma_{\dot {\widetilde{x}}}(0)]=1
$$
и, так как второй и третий члены в (\ref{put}) $T$-периодичны,
(\ref{wehave}) может быть переписано как
$$
  d_{\mathbb{R}^2}(F,
  U)= \frac{1}{2}\left[{\rm ind}(\theta_1,f){\rm
sign}\left<z(\theta_1)^\bot,F(\widetilde{x}(\theta_1))\right>+\right.
$$
\begin{equation}\left.
+
  {\rm ind}(\theta_2,f){\rm
  sign}\left<z(\theta_2)^\bot,F(\widetilde{x}(\theta_2))\right>\right].\label{lf}
\end{equation}
Так как функция $f$ $T$-периодична, то
\begin{equation}\label{llf}
{\rm ind}(\theta_1,f)=-{\rm ind}(\theta_2,f).
\end{equation}
На основании соотношений 3) и (\ref{llf}) утверждение леммы
вытекает из (\ref{lf}).

Лемма доказана.

\begin{lem}\label{lem_nu_mu}
 Рассмотрим систему
\begin{equation}\label{lin}
  \begin{array}{lll}
    \dot x_1&=&x_2\\
    \dot x_2&=&-x_1+\varepsilon(\mu x^+_1+\nu x_1^-+\cos(t)),
  \end{array}
\end{equation}
где $a^+:=\max\{a,0\},$ $a^-:=\max\{-a,0\}.$
 Обозначим через $U\subset\mathbb{R}^2$ внутренность цикла
$\widetilde{x}=(\sin t,\cos t)$ системы (\ref{lin}) с
$\varepsilon=0.$ Тогда, если $|\mu-\nu|\not=2,$ то соответствующий
$2\pi$-периодической системе (\ref{lin}) обобщенный оператор
усреднения $\Phi^s$ невырожден на $\partial U$ при $s\in[0,2\pi].$
Если же $|\mu-\nu|<2,$ то $d_{\mathbb{R}^2}(-\Phi^T,U)\in\{0,2\}.$
\end{lem}

\noindent {\bf Доказательство.} Имеем
$$
  \widehat{f}=-\int\limits_0^{2\pi}\sin\tau\left(\mu \widetilde{x}^+_1(\tau)+\nu
  \widetilde{x}_1^-(\tau)+\cos(\tau-\theta)\right)d\tau=
$$
$$
  =-\mu\int\limits_0^\pi\sin\tau\sin\tau
  d\tau+\nu\int\limits_\pi^{2\pi}\sin\tau\sin\tau
  d\tau-\int\limits_0^{2\pi}\sin\tau\cos\tau
  d\tau\cos\theta-
$$
$$
-\int\limits_0^{2\pi}\sin\tau\sin\tau d\tau\sin\theta
  = -\mu\frac{\pi}{2}+\nu\frac{\pi}{2}-\pi\sin\theta,
$$
$$
  \widetilde{f}=\int\limits_0^{2\pi}\cos\tau\left(\mu \widetilde{x}^+_1(\tau)+\nu
  \widetilde{x}_1^-(\tau)+\cos(\tau-\theta)\right)d\tau=
$$
$$
  =\mu\int\limits_0^\pi\cos\tau\sin\tau
  d\tau-\nu\int\limits_\pi^{2\pi}\cos\tau\sin\tau
  d\tau+\int\limits_0^{2\pi}\cos\tau\cos\tau
  d\tau\cos\theta+
$$
$$
+\int\limits_0^{2\pi}\cos\tau\sin\tau d\tau\sin\theta
  =\pi\cos\theta,
$$
поэтому
$$
  \Phi^s(\widetilde{x}(\theta))=\frac{\pi}{2}\left(-\mu+\nu-2\sin\theta\right)\left(
  \begin{array}{c} \cos\theta\\
  -\sin\theta\end{array}\right)+\pi\cos\theta
  \left(
  \begin{array}{c} \sin\theta\\
  \cos\theta\end{array}\right).
$$
Из последней формулы следует, что необходимым и достаточным
условием невырожденности поля $\Phi^s$ является $|\mu-\nu|\not=2.$
Для доказательства того, что $d_{\mathbb{R}^2}(-\Phi^T,U)\not=0$
используем лемму~\ref{lem_barsuk} с $F=-\Phi^T$ и направляющей
функцией $z(\theta)=(\cos\theta,\sin\theta).$ Проверим выполнение
условий этой леммы, имеем
$$
  \left<-\Phi^T(\widetilde{x}(\theta)),z(\theta)\right>=-\frac{\pi}{2}\left(-\mu+\nu-2\sin\theta\right).
$$
Так как по условию леммы $|\mu-\nu|<2,$ то
$\theta_0=\arcsin\dfrac{-\mu+\nu}{2}$ будет единственным корнем
уравнения
$\left<-\Phi^T(\widetilde{x}(\theta)),z(\theta)\right>=0$ на
интервале $[-\pi/2,\pi/2].$ Поэтому на интервале $[0,2\pi)$
уравнение
$\left<-\Phi^T(\widetilde{x}(\theta)),z(\theta)\right>=0$ имеет
ровно два корня
$$
\begin{array}{l}
  \theta_1=\left\{\begin{array}{l}\theta_1=\theta_0,\quad\mbox{если
  }\theta_0\ge 0,\\
  \theta_1=\theta_0+\pi,\quad\mbox{в противном
  случае,}\end{array}\right.\\
  \theta_2=\theta_1+\pi.
\end{array}
$$
Далее, имеем
$$
  \left(\left<-\Phi^T(\widetilde{x}(\cdot)),z(\cdot)\right>\right)'(\theta)=\cos\theta,
$$
поэтому, если
$\left(\left<-\Phi^T(\widetilde{x}(\cdot)),z(\cdot)\right>\right)'(\theta_i)=0,$
то $|-\mu+\nu|=2,$ противореча условиям леммы. Таким образом,
условие~2) леммы~\ref{lem_barsuk} удовлетворено. Наконец, так как
$${\rm
sign}\left<-\Phi(\widetilde{x}(\theta_i)),\left(
  \begin{array}{c} z_2(\theta_i)\\
  -z_1(\theta_i)\end{array}\right)\right>=-\pi\cos\theta_i,
$$
то условие~3) леммы~\ref{lem_barsuk} также выполнено.

Лемма доказана.

В следующем пункте главы лемма~\ref{lem_barsuk}
 используется для анализа более общего нелинейного
случая.

\noindent {\bf Доказательство леммы \ref{duf_lem}.} Если $u$ --
решение уравнения (\ref{so}), то  $v=(u,\dot u)$ удовлетворяет
системе
\begin{equation}\label{duf}
  \begin{array}{lll}
    \dot v_1&=&v_2\\
    \dot v_2&=&-v_1- v_1^3+\varepsilon \cos((1+\delta)t),
  \end{array}
\end{equation}
обратно, если $v$ -- решение системы (\ref{duf}), то $v_1$ --
решение уравнения (\ref{so}). Без ограничения общности можно
считать, что $\dot u_\delta(0)=0,$ тогда $u_\delta(0)\not=0.$
Заменой переменных
$$
  x(t)=\frac{v(t)}{u_\delta(0)}
$$
перейдем от системы (\ref{duf}) к системе
\begin{equation}\label{duf1}
  \begin{array}{lll}
    \dot x_1&=&x_2\\
    \dot x_2&=&-x_1-(u_\delta(0))^2x_1^3+\varepsilon
    \dfrac{1}{u_\delta(0)}\cos((1+\delta)t).
  \end{array}
\end{equation}
Таким образом, для доказательства леммы~\ref{duf_lem} достаточно
установить, что существует $\delta_0>0$ такое, что при
$\delta\in(0,\delta_0]$ условия теоремы~\ref{dimtwo_1}, связанные
с оператором $\Phi^s$ системы (\ref{duf1}),  выполнены с
$\widetilde{x}(t)=\dfrac{v_\delta(t)}{u_\delta(0)}$ и
$T=2\pi/(1+\delta).$ Для этого, в свою очередь, достаточно
установить аналогичное утверждение для системы
\begin{equation}\label{duf11}
  \begin{array}{lll}
    \dot x_1&=&x_2\\
    \dot x_2&=&-x_1-(u_\delta(0))^2x_1^3+\varepsilon
    \cos((1+\delta)t).
  \end{array}
\end{equation}
Так как период периодических решений порождающей системы
(\ref{duf}) изменяется монотонно от $2\pi$ до $0,$ когда начальное
условие решения изменяется от $(0,0)$ до $(+\infty,0)$ (см.
\cite{guck}, пример~с.~250), то $u_\delta(0)\to 0,$ когда
$\delta\to 0.$ Но для $\delta=0$ справедливость желаемого для
системы (\ref{duf11}) утверждения следует из леммы
\ref{lem_nu_mu}, следовательно, это утверждение остается
справедливым и при малых $\delta>0.$

Лемма доказана.

\section{Симметричные и вырожденные двумерные случаи }\label{sym}

В этом пункте будут рассмотрены два случая, в которых условия
теорем~2.5 и 2.6 упрощаются.

\subsubsection{Случай, когда рассматриваемая система удовлетворяет условиям симметрии}

Рассмотрим систему
\begin{equation}\label{ps2}
 \dot x=
f(x)+\varepsilon \sin(wt)g(x),
\end{equation}
где при каждом $\xi\in\mathbb{R}^2$ имеют место соотношения
\begin{eqnarray}
  f_1(\xi)&=&f_1(-\xi_1,\xi_2),\label{m1}\\
  f_2(\xi)&=&-f_2(-\xi_1,\xi_2),\label{m2}\\
  f_1(\xi)&=&-f_1(\xi_1,-\xi_2),\label{K1}\\
  f_2(\xi)&=&f_2(\xi_1,-\xi_2),\label{K2}\\
  (f_1)'_{(1)}(\xi)&=&-(f_2)'_{(2)}(\xi),\label{m3}
\end{eqnarray}
\begin{equation}\label{symg}
  g(\xi)=\left(\begin{array}{c}-g_1(\xi_1,-\xi_2) \\
  g_2(\xi_1,-\xi_2)\end{array}\right),
\end{equation}
\begin{equation}\label{K3}
  g(\xi)=\left(\begin{array}{c}-g_1(-\xi_1,\xi_2) \\
  g_2(-\xi_1,\xi_2)\end{array}\right).
\end{equation}
 Будем считать, что  $\widetilde {x}$ -- периодический цикл порождающей системы
 (\ref{npa_1})
наименьшего периода $2\pi/w,$ удовлетворяющий условиям
(\ref{init}), $(C)$ (см. подпункт~2.2.3), и
\begin{equation}\label{U}
  \widetilde{x}_1(0)=0,\quad \widetilde{x}_2(0)\not=0.
\end{equation}
Обозначим через $U\subset\mathbb{R}^2$ внутренность цикла
$\widetilde {x}.$ Пусть $\widehat {y}$ -- решение линеаризованной
системы (\ref{lsg}), удовлетворяющее начальному условию
\begin{equation}\label{U1}
\widehat{y}(0)=\left(0,1/\dot{\widetilde{x}}_1(0)\right).
\end{equation}
Предположим далее, что
\begin{equation}\label{STARR}
  \partial U\cap\left(\left(0,+\infty\right)\times 0\right)\not=0,
\end{equation}
\begin{equation}\label{positf}
f_1(\xi)>0,\ f_2(\xi)<0,\qquad\xi\in\partial
U\cap\left((0,+\infty)\times(0,+\infty)\right),
\end{equation}
\begin{equation}\label{positg}
g_1(\xi)>0\mbox{ и }g_2(\xi)>0,\qquad\xi\in\partial
U\cap\left((0,+\infty)\times(0,+\infty)\right).
\end{equation}

Введем $\widetilde{\xi},\ \widehat{\xi}\in\mathbb{R}^2$ как
\begin{eqnarray}
  \widetilde{\xi} &=& 4\int\limits_0^{\pi/(2w)}\left<
\left(\begin{array}{c} -\dot{\widetilde{x}}_2(\tau) \\
\dot{\widetilde{x}}_1(\tau)\end{array}\right),
  g(\widetilde{x}(\tau))\right>\left(\begin{array}{c}\cos(w\tau) \\ \sin(w\tau)\end{array}\right)d\tau,\nonumber\\
  \widehat{\xi} &=&\int\limits_0^{2\pi/w}\left<
\left(\begin{array}{c} \widehat{y}_2(\tau) \\
-\widehat{y}_1(\tau)\end{array}\right),
  g(\widetilde{x}(\tau))\right>\left(\begin{array}{c}\cos(w\tau) \\
  \sin(w\tau)\end{array}\right)d\tau.\nonumber
\end{eqnarray}

Теорема 2.6 позволяет доказать следующее достаточное условие
существования $T$-периодических решений в системе (\ref{ps2}).

\begin{thm}\label{thm_sym} Пусть
$\widetilde {x}$ -- периодический цикл порождающей системы
 (\ref{npa_1})
наименьшего периода $2\pi/w,$ удовлетворяющий условиям $(A_1),$
$(C),$ (\ref{init})   и (\ref{U}). Пусть выполнены условия
(\ref{m1})-(\ref{K3}), (\ref{STARR})-(\ref{positg}). Тогда, если
\begin{equation}\label{lincond}
\left|\widehat{\xi_1}\right|+\dfrac{\left|\widehat{y}_1(2\pi/w)\right|}{\dot{\widetilde{x}}_1(0)}
\widetilde{\xi_1}<\min\left\{\left|\widehat{\xi}_2\right|-
\dfrac{\left|\widehat{y}_1(2\pi/w)\right|}{4\dot{\widetilde{x}}_1(0)}\widetilde{\xi}_2,\widetilde{\xi}_1
\right\},
\end{equation}
то при всех достаточно малых $\varepsilon>0$ система (\ref{ps2})
имеет по крайней мере два $2\pi/w$-периодических решения
$\widetilde{x}_\varepsilon,$
$\widetilde{\widetilde{x}}_\varepsilon$
 таких, что
$\widetilde{x}_\varepsilon(t)\in U,$
$\widetilde{\widetilde{x}}_\varepsilon(t)\not\in U$ для любого
$t\in[0,2\pi/w],$ и
$$\rho\left(\widetilde{x}_\varepsilon(t),\partial U\right)+
\rho\left(\widetilde{\widetilde{x}}_\varepsilon(t),\partial U
\right)\to 0\quad\mbox{при }\varepsilon\to 0
$$
равномерно по $t\in[0,2\pi/w].$ Прочие $2\pi/w$-периодические
решения $x$ системы (\ref{ps2}) удовлетворяют условию
$x(t)\not\in\partial U$ для любых $t\in\left[0,2\pi/w\right]$ и
достаточно малых $\varepsilon>0.$
\end{thm}

Нам понадобится следующая лемма.

\begin{lem}\label{lem_sopr}
Рассмотрим линейную систему
\begin{equation}\label{lin1}
  \left(\begin{array}{l}
    \dot y_1\\
    \dot y_2
  \end{array}\right)=\left(\begin{array}{cc}
    a(t) & d(t)\\
    b(t) & -a(t)
  \end{array}\right)\left(\begin{array}{l}
    y_1\\
    y_2
  \end{array}\right).
\end{equation}
Предположим, что $a(-t)=-a(t),$ $b(-t)=b(t),$ $d(-t)=d(t)$ для
любых $t\in[c_1,c_2].$ Тогда, если $y$ -- некоторое решение
системы (\ref{lin1}), то функция $z(t)=(y_2(-t),y_1(-t))$
удовлетворяет на отрезке $[c_1,c_2]$ сопряженной к (\ref{lin1})
системе.
\end{lem}

Справедливость леммы \ref{lem_sopr} проверяется непосредственной
подстановкой решения $z(t)=(y_2(-t),y_1(-t))$ в сопряженную к
(\ref{lin1}) систему.

\noindent {\bf Доказательство теоремы \ref{thm_sym}.} Пользуясь
условиями (\ref{m1}) и (\ref{m2}), легко проверить, что функция
$p(t)=(-\widetilde{x}_1(-t),\widetilde{x}_2(-t))$ является
решением системы (\ref{npa_1}). Но из (\ref{U}) имеем
$p(0)=\widetilde{x}(0),$ следовательно,
\begin{equation}\label{sym1}
 (-\widetilde{x}_1(-t),\widetilde{x}_2(-t))=\widetilde{x}(t) \mbox{ для любого
 }t\in\mathbb{R}.
\end{equation}
Линеаризуя   систему (\ref{ps2}) при $\varepsilon=0$ на цикле
$\widetilde{x},$ имеем систему
\begin{equation}\label{slin}
  \left(\begin{array}{l}
    \dot y_1\\
    \dot y_2
  \end{array}\right)=\left(\begin{array}{cc}
    (f_1)'_{(1)}(\widetilde{x}(t)) & (f_1)'_{(2)}(\widetilde{x}(t))\\
    (f_2)'_{(1)}(\widetilde{x}(t)) & (f_2)'_{(2)}(\widetilde{x}(t))
  \end{array}\right)\left(\begin{array}{l}
    y_1\\
    y_2
  \end{array}\right).
\end{equation}
Из условия (\ref{m3}) следует, что
\begin{equation}\label{U3}
  (f_1)'_{(1)}(\widetilde{x}(t))=-(f_2)'_{(2)}(\widetilde{x}(t)),\quad
  t\in\mathbb{R}.
\end{equation}
Из (\ref{m1}) имеем
$-(f_1)'_{(1)}(-\xi_1,\xi_2)=(f_1)'_{(1)}(\xi_1,\xi_2),$
$\in\in\mathbb{R}^2,$ и, учитывая (\ref{sym1}), получаем
\begin{equation}\label{m5}
  (f_1)'_{(1)}(\widetilde{x}(t))=-(f_1)'_{(1)}(\widetilde{x}(-t)),
  \quad t\in\mathbb{R}.
\end{equation}
Из (\ref{m2}) имеем
$(f_2)'_{(1)}(\xi_1,\xi_2)=(f_2)'_{(1)}(-\xi_1,\xi_2),$
$\xi\in\mathbb{R}^2,$ и, учитывая (\ref{sym1}), получаем
\begin{equation}\label{m4}
  (f_2)'_{(1)}(\widetilde{x}(-t))=(f_2)'_{(1)}(\widetilde{x}(t)),\quad
  t\in\mathbb{R}.
\end{equation}
Наконец, из (\ref{m1}) имеем
$(f_1)'_{(2)}(-\xi_1,\xi_2)=(f_1)'_{(2)}(\xi_1,\xi_2),$
$\xi\in\mathbb{R}^2,$ и, учитывая (\ref{sym1}), получаем
\begin{equation}\label{m6}
  (f_1)'_{(2)}(\widetilde{x}(t))=(f_1)'_{(2)}(\widetilde{x}(-t)),\quad
  t\in\mathbb{R}.
\end{equation}
Таким образом, выполнены условия леммы~\ref{lem_sopr}, на
основании которой, учитывая также (\ref{init}) и (\ref{U1}),
заключаем, что функции
\begin{equation}\label{funkcii}
 \widehat {z}(t)=\left(\begin{array}{c} \widehat {y}_2(-t)
 \\ \widehat {y}_1(-t)\end{array}\right)\quad\mbox{и}\quad
\widetilde {z}(t)=\left(\begin{array}{c} \dot {\widetilde
{x}}_2(-t)
 \\\dot {\widetilde {x}}_1(-t)\end{array}\right)
\end{equation}
удовлетворяют сопряженной к (\ref{slin}) системе и, вместе с тем,
условию (\ref{initz}). Из (\ref{sym1}) для любого $t\in\mathbb{R}$
имеем
$$\dot {\widetilde{x}}_1(-t)=\dot {\widetilde{x}}_1(t),\quad -\dot {\widetilde{x}}_2(-t)=\dot{\widetilde{x}}_2(t).$$
Покажем, что вместе с $\widehat{y}$ решением системы (\ref{slin})
является функция $p(t)=(-\widehat{y}_1(-t),\widehat{y}_2(-t)).$
Действительно, из (\ref{m5}) и (\ref{m6}) имеем
$$
 \dot
 p_1(t)=-(f_1)'_{(1)}(\widetilde{x}(t))y_1(-t)+(f_1)'_{(2)}(\widetilde{x}(t))y_2(-t),
$$
и из (\ref{m4}), (\ref{U3}) и (\ref{m5}) заключаем
$$
 \dot
 p_2(t)=-(f_2)'_{(1)}(\widetilde{x}(t))y_1(-t)+(f_2)'_{(2)}(\widetilde{x}(t))y_2(-t).
$$
Но из (\ref{U1}) следует, что $p(0)=\widehat{y}(0),$ поэтому
\begin{equation}\label{symy}
  (-\widehat{y}_1(-t),\widehat{y}_2(-t))=\widehat{y}(t), \quad
  t\in\mathbb{R}.
\end{equation}
Таким образом, функции $\widehat{z}$ и $\widetilde{z}$ можно
переписать в виде
\begin{equation}\label{STAR2}
 \widehat {z}(t)=\left(\begin{array}{c} \widehat {y}_2(t)
 \\ -\widehat {y}_1(t)\end{array}\right)\quad\mbox{и}\quad
\widetilde {z}(t)=\left(\begin{array}{c} -\dot {\widetilde
{x}}_2(t)
 \\\dot {\widetilde {x}}_1(t)\end{array}\right).
\end{equation}
Введем в рассмотрение вспомогательное векторное поле
$F:\mathbb{R}^2\to\mathbb{R}^2,$ определив его на $\partial U$ как
$$
  F\left(\widetilde{x}(\theta)\right)={\rm
  sign}\left(\widehat{\xi_2}\right)\cos(w\theta)
\left(\begin{array}{c}
  {\widehat{y}}_2(\theta)\\
  -{\widehat{y}}_1(\theta)\end{array}\right)-\sin(w\theta)\left(\begin{array}{c}
  -\dot{\widetilde{x}}_2(\theta)\\
  \dot{\widetilde{x}}_1(\theta)\end{array}\right).
$$
Очевидно, поле $F$ невырожденно на $\partial U,$ покажем, что для
него выполнены условия леммы~\ref{lem_barsuk} с направляющей
функцией $z(t):=\dot{\widetilde{x}}(t).$ Во-первых, заметим, что
\begin{equation}\label{detx}
  \left<\left(\begin{array}{c}
  \dot{\widetilde{x}}_1(\theta)\\
  \dot{\widetilde{x}}_2(\theta)\end{array}\right),
  \left(\begin{array}{c}
  {\widehat{y}}_2(\theta)\\
  -{\widehat{y}}_1(\theta)\end{array}\right)\right>=
  {\rm det}\left\| \left(\begin{array}{cc}
  \dot{\widetilde{x}}_1(\theta) & {\widehat{y}}_1(\theta)\\
  \dot{\widetilde{x}}_2(\theta) & {\widehat{y}}_2(\theta)\end{array}\right)\right\|.
\end{equation}
Поэтому функция $\left<F(\widetilde{x}(\cdot)),z(\cdot)\right>$ из
условия~2 леммы~\ref{lem_barsuk} имеет вид
$$
  \left<F(\widetilde{x}(\theta)),\dot{\widetilde{x}}(\theta)\right>={\rm
  sign}\left(\widehat{\xi_2}\right)\cos(w\theta){\rm det}\left\| \left(\begin{array}{cc}
  \dot{\widetilde{x}}_1(\theta) & {\widehat{y}}_1(\theta)\\
  \dot{\widetilde{x}}_2(\theta) & {\widehat{y}}_2(\theta)\end{array}\right)\right\|,
$$
пользуясь которым легко установить, что эта функция  допускает
ровно два нуля $\theta_1=\pi/(2w)$ и $\theta_2=3\pi/w$ на
интервале $[0,2\pi/w)$ и строго монотонна в указанных точках. В
тоже время
$$
  \left<F(\widetilde{x}(\theta_i)),\left(\begin{array}{c}
  \dot{\widetilde{x}}_2(\theta_i)\\
  -\dot{\widetilde{x}}_1(\theta_i)\end{array}\right)\right>=\sin(w\theta_i)\left\|\dot{\widetilde{x}}(\theta_i)\right\|
$$
и, значит,
$$
  {\rm sign}\left<F(\widetilde{x}(\pi/(2w))),\left(\begin{array}{c}
  \dot{\widetilde{x}}_2(\pi/(2w))\\
  -\dot{\widetilde{x}}_1(\pi/(2w))\end{array}\right)\right>=
$$
$$
=-{\rm sign}
  \left<F(\widetilde{x}(3\pi/w)),\left(\begin{array}{c}
  \dot{\widetilde{x}}_2(3\pi/w)\\
  -\dot{\widetilde{x}}_1(3\pi/w)\end{array}\right)\right>.
$$
Таким образом, все условия леммы~\ref{lem_barsuk} удовлетворены и,
следовательно, $d_{\mathbb{R}^2}(F,U)\in\{0,2\}.$

Условие $(C)$ и лемма~2.4 позволяют утверждать, что для
обобщенного оператора усреднения $\Phi^s,$ соответствующего
системе (\ref{ps2}), справедлива формула (\ref{Ftheta}). Покажем,
что условия теоремы гарантируют
\begin{equation}\label{gar}
  \left<\Phi^s(\widetilde{x}(\theta)),F(\widetilde{x}(\theta))\right>\not=0,\quad
  s,\theta\in[0,2\pi/w].
\end{equation}
Из (\ref{detx}) следует,что
$$
\left<\dot{\widetilde{x}}(\theta),
  \left(\begin{array}{c}
  {\widehat{y}}_2(\theta)\\
  -{\widehat{y}}_1(\theta)\end{array}\right)\right>=1,\quad\theta\in[0,2\pi/w],
$$
откуда имеем также
$$
\left<\widehat{y}(\theta),
  \left(\begin{array}{c}
  -\dot{\widetilde{x}}_2(\theta)\\
  \dot{\widetilde{x}}_1(\theta)\end{array}\right)\right>=1,\quad\theta\in[0,2\pi/w].
$$
Поэтому, используя формулу (\ref{Ftheta}), можем записать
$$
  \left<\Phi^s(\widetilde{x}(\theta)),F(\widetilde{x}(\theta))\right>=
$$
$$=
  \left(\int\limits_0^{2\pi/w}\left<
\left(\begin{array}{c} \widehat{y}_2(\tau) \\
-\widehat{y}_1(\tau)\end{array}\right),
  g(\widetilde{x}(\tau))\right>\sin(w\tau-w\theta)d\tau+\right.
$$
$$+\left.
  \frac{\widehat{y}_1(2\pi/w)}{\dot{\widetilde{x}}_1(0)}
\int\limits_{s+\theta}^{2\pi/w}\left<
\left(\begin{array}{c} -\dot{\widetilde{x}}_2(\tau) \\
\dot{\widetilde{x}}_1(\tau)\end{array}\right),
  g(\widetilde{x}(\tau))\right>\sin(w\tau-w\theta)d\tau\right)\cdot
$$
$$\cdot{\rm
  sign}\left(\widehat{\xi_2}\right)\cos(w\tau)
  -
$$
$$
-\left(\int\limits_{0}^{2\pi/w}\left<
\left(\begin{array}{c} -\dot{\widetilde{x}}_2(\tau) \\
\dot{\widetilde{x}}_1(\tau)\end{array}\right),
  g(\widetilde{x}(\tau))\right>\sin(w\tau-w\theta)d\tau\right)\cdot
$$
$$
\cdot\sin(w\tau).
$$
Покажем, что
\begin{equation}\label{pok}
\int\limits_{0}^{2\pi/w}\left<
\left(\begin{array}{c} -\dot{\widetilde{x}}_2(\tau) \\
\dot{\widetilde{x}}_1(\tau)\end{array}\right),
  g(\widetilde{x}(\tau))\right>\left(\begin{array}{c}\cos(w\tau) \\
\sin(w\tau)\end{array}\right)d\tau=\left(\begin{array}{c}\widetilde{\xi}_1 \\
0\end{array}\right).
\end{equation}
Из (\ref{K3}) и (\ref{sym1}) для любого $t\in\mathbb{R}$ имеем
\begin{equation}\label{usl}
\begin{array}{c}\dot {\widetilde{x}}_1(\pi/w-t)=\dot
{\widetilde{x}}_1(\pi/w+t),\\ g_1\left(\widetilde{x}(\pi/w-t)\right)=-g_1\left(\widetilde{x}(\pi/w+t)\right),\\
\dot {\widetilde{x}}_2(\pi/w-t)=-\dot{\widetilde{x}}_2(\pi/w+t),
\\
 g_2\left(\widetilde{x}(\pi/w-t)\right)=g_2\left(\widetilde{x}(\pi/w+t)\right).\end{array}
\end{equation}
В тоже время, {$\cos(\pi-wt)=\cos(\pi+wt)$ и
$\sin(\pi-wt)=-\sin(\pi+wt),$ поэтому при всех $t\in\mathbb{R}$
\begin{equation}\label{N}
\begin{array}{l}
  \dot{\widetilde{x}}_2(\pi/w-t)\sin(\pi-wt)g_1(\widetilde{x}(\pi/w-t))
   =\\
  \qquad=-\dot{\widetilde{x}}_2(\pi/w+t)\sin(\pi+wt)g_1(\widetilde{x}(\pi/w+t)),\\
  \dot{\widetilde{x}}_1(\pi/w-t)\sin(\pi-wt)g_2(\widetilde{x}(\pi/w-t))
  = \\
  \qquad=-\dot{\widetilde{x}}_1(\pi/w+t)\sin(\pi+wt)g_2(\widetilde{x}(\pi/w+t)),
\end{array}
\end{equation}
и
\begin{equation}\label{f0}
\int\limits_0^{2\pi/w}\dot{\widetilde{
  x}}_2(\tau)\sin(w\tau)g_1(\widetilde{x}(\tau))d\tau=
\int\limits_0^{2\pi/w}\dot{\widetilde{
  x}}_1(\tau)\sin(w\tau)g_2(\widetilde{x}(\tau))d\tau=0.
\end{equation}
Аналогично при всех $t\in\mathbb{R}$
\begin{equation}\label{O1}
\begin{array}{l}
  \dot{\widetilde{x}}_2(\pi/w-t)\cos(\pi-wt)g_1(\widetilde{x}(\pi/w-t))
   =\\
  \qquad=\dot{\widetilde{x}}_2(\pi/w+t)\cos(\pi+wt)g_1(\widetilde{x}(\pi/w+t)),\\
  \dot{\widetilde{x}}_1(\pi/w-t)\cos(\pi-wt)g_2(\widetilde{x}(\pi/w-t))
  = \\
  \qquad=\dot{\widetilde{x}}_1(\pi/w+t)\cos(\pi+wt)g_2(\widetilde{x}(\pi/w+t)),
\end{array}
\end{equation}
 и
$$
 \int\limits_{0}^{2\pi/w}\left<
\left(\begin{array}{c} -\dot{\widetilde{x}}_2(\tau) \\
\dot{\widetilde{x}}_1(\tau)\end{array}\right),
  g(\widetilde{x}(\tau))\right>\cos(w\tau)d\tau=
$$
\begin{equation}\label{pre1}
=  2\int\limits_{0}^{\pi/w}\left<
\left(\begin{array}{c} -\dot{\widetilde{x}}_2(\tau) \\
\dot{\widetilde{x}}_1(\tau)\end{array}\right),
  g(\widetilde{x}(\tau))\right>\cos(w\tau)d\tau.
\end{equation}
Далее, по условию (\ref{STARR}) теоремы существует $t_*>0$ такое,
что $\widetilde{x}_2(t_*)=0,$ без ограничения общности можно
считать, что
\begin{equation}\label{posit}
\widetilde{x}_2(t)\not=0\mbox{ при всех }t\in[0,t_*).
\end{equation}
Из (\ref{sym1}) следует, что $t_*\in[0,\pi/w].$ Пользуясь
(\ref{K1}) и (\ref{K2}), легко установить, что системе
(\ref{npa_1}) удовлетворяет функция
$p(t)=(\widetilde{x}_1(t_*-t),-\widetilde{x}_2(t_*-t)).$ В тоже
время, в силу автономности системе (\ref{npa_1}) удовлетворяет
также и функция
$q(t)=(\widetilde{x}_1(t_*+t),\widetilde{x}_2(t_*+t)).$ Но
$p(0)=q(0),$ значит
\begin{equation}\label{pere}
(\widetilde{x}_1(t_*-t),-\widetilde{x}_2(t_*-t))=(\widetilde{x}_1(t_*+t),\widetilde{x}_2(t_*+t)),\quad
t\in\mathbb{R}
\end{equation}
и, в частности,
$(0,-\widetilde{x}_2(0))=(\widetilde{x}_1(2
t_*),\widetilde{x}_2(2t_*)).$ Из $0=\widetilde{x}_1(2 t_*),$ на
основании (\ref{sym1}), следует, что число $4 t_*$ является
периодом $\widetilde{x}.$ С другой стороны, из
$-\widetilde{x}_2(0)=\widetilde{x}_2(2t_*)$ имеем, что
$t_*\not=\pi/w$ и, следовательно, $4t_*\in[0,4\pi/w).$ Но на
интервале $[0,4\pi/w)$ есть только одно число $2\pi/w,$ являющееся
периодом $\widetilde{x},$ поэтому $4t_*=2\pi/w$ и $t_*=\pi/(2w).$
Таким образом, для любого $t\in\left[0,\mathbb{R}\right]$
равенство (\ref{pere}) можно переписать следующим образом
$$
\left(\widetilde{x}_1\left(\frac{\pi}{2w}-t\right),-\widetilde{x}_2\left(\frac{\pi}{2w}-t\right)\right)=
\left(\widetilde{x}_1\left(\frac{\pi}{2w}+t\right),\widetilde{x}_2\left(\frac{\pi}{2w}+t\right)\right),
$$
откуда, учитывая (\ref{symg}), для любого $t\in\mathbb{R}$ имеем
\begin{equation}\label{usl2}
\begin{array}{c}
\dot {\widetilde{x}}_1(\pi/(2w)-t)=-\dot
{\widetilde{x}}_1(\pi/(2w)+t), \\
g_1\left(\widetilde{x}(\pi/(2w)-t)\right)=-g_1\left(\widetilde{x}(\pi/(2w)+t)\right),\\
\dot
{\widetilde{x}}_2(\pi/(2w)-t)=\dot{\widetilde{x}}_2(\pi/(2w)+t),
\\
g_2\left(\widetilde{x}(\pi/(2w)-t)\right)=g_2\left(\widetilde{x}(\pi/(2w)+t)\right).
\end{array}
\end{equation}
Пользуясь полученными соотношениями и учитывая, что
$\sin(\pi/2-wt)=\sin(\pi/2+wt),$ при всех $t\in\mathbb{R}$ имеем
\begin{equation}\label{N1}
\begin{array}{l}
  \dot{\widetilde{x}}_2(\pi/w-t)\sin(\pi/2-wt)g_1(\widetilde{x}(\pi/w-t))
   =\\
  \qquad=-\dot{\widetilde{x}}_2(\pi/w+t)\sin(\pi/2+wt)g_1(\widetilde{x}(\pi/w+t)),\\
  \dot{\widetilde{x}}_1(\pi/w-t)\sin(\pi/2-wt)g_2(\widetilde{x}(\pi/w-t))
  = \\
  \qquad=-\dot{\widetilde{x}}_1(\pi/w+t)\sin(\pi/2+wt)g_2(\widetilde{x}(\pi/w+t)).
\end{array}
\end{equation}
Аналогично, учитывая, что $\cos(\pi/2-wt)=-\cos(\pi/2+wt),$ при
всех $t\in\mathbb{R}$ имеем
\begin{equation}\label{O2}
\begin{array}{l}
  \dot{\widetilde{x}}_2(\pi/w-t)\cos(\pi/2-wt)g_1(\widetilde{x}(\pi/w-t))
   =\\
  \qquad=\dot{\widetilde{x}}_2(\pi/w+t)\cos(\pi/2+wt)g_1(\widetilde{x}(\pi/w+t)),\\
  \dot{\widetilde{x}}_1(\pi/w-t)\cos(\pi/2-wt)g_2(\widetilde{x}(\pi/w-t))
  = \\
  \qquad=\dot{\widetilde{x}}_1(\pi/w+t)\cos(\pi/2+wt)g_2(\widetilde{x}(\pi/w+t)).
\end{array}
\end{equation}
Из (\ref{O2}) следует, что
$$
 \int\limits_{0}^{\pi/w}\left<
\left(\begin{array}{c} -\dot{\widetilde{x}}_2(\tau) \\
\dot{\widetilde{x}}_1(\tau)\end{array}\right),
  g(\widetilde{x}(\tau))\right>\cos(w\tau)d\tau=
$$
\begin{equation}\label{pre2}
=  2\int\limits_{0}^{\pi/(2w)}\left<
\left(\begin{array}{c} -\dot{\widetilde{x}}_2(\tau) \\
\dot{\widetilde{x}}_1(\tau)\end{array}\right),
  g(\widetilde{x}(\tau))\right>\cos(w\tau)d\tau.
\end{equation}
Равенства (\ref{f0}), (\ref{pre1}) и (\ref{pre2}) означают, что
соотношение (\ref{pok}) выполнено и, вводя
$\mathcal{F}:\mathbb{R}\to\mathbb{R}^2$ как
\begin{equation}\label{N2}
\mathcal{F}(t)= \int\limits_t^{2\pi/w}\left<
\left(\begin{array}{c} -\dot{\widetilde{x}}_2(\tau) \\
\dot{\widetilde{x}}_1(\tau)\end{array}\right),
g(\widetilde{x}(\tau))\right>\left(\begin{array}{c}\cos(w\tau) \\
\sin(w\tau)\end{array}\right)d\tau,
\end{equation}
можем переписать
$\left<\Phi^s(\widetilde{x}(\theta)),F(\widetilde{x}(\theta))\right>$
в виде
$$
\left<\Phi^s(\widetilde{x}(\theta)),F(\widetilde{x}(\theta))\right>
=$$
$$
=\left(\left<\widehat{\xi},\left(\begin{array}{c}-\sin(w\theta) \\
\cos(w\theta)\end{array}\right)\right>+
\frac{\widehat{y}_1(2\pi/w)}{\dot{\widetilde{x}}_1(0)}\left<\mathcal{F}(s+\theta),
\left(\begin{array}{c}-\sin(w\theta) \\
\cos(w\theta)\end{array}\right)\right>\right)\cdot $$ $$\cdot{\rm
  sign}\left(\widehat{\xi_2}\right)\cos(w\theta)-
$$
$$
  -\left<\left(\begin{array}{c} \widetilde{\xi}_1\\0
\end{array}\right)
,\left(\begin{array}{c}-\sin(w\theta) \\
\cos(w\theta)\end{array}\right)\right>\sin(w\theta)=
$$
$$=\left(\left(-\widehat{\xi}_1-\frac{\widehat{y}_1(2\pi/w)}{\dot{\widetilde{x}}_1(0)}\mathcal{F}_1(s+\theta)\right)
  \sin(w\theta){\rm
  sign}\left(\widehat{\xi_2}\right)\cos(w\theta)+\right.
$$
$$
\left.+\left(\widehat{\xi}_2{\rm
  sign}\left(\widehat{\xi_2}\right)+\frac{\widehat{y}_1(2\pi/w)}{\dot{\widetilde{x}}_1(0)}\mathcal{F}_2(s+\theta){\rm
  sign}\left(\widehat{\xi_2}\right)\right)
  \cos^2(w\theta)+\right.$$
  $$\left.+\widetilde{\xi}_1\sin^2(w\theta)\right).
$$
Таким образом, для доказательства желаемого факта (\ref{gar})
достаточно установить, что
$$
  \min_{s,\theta\in[0,2\pi/w]}\left|
\left(\left|\widehat{\xi}_2\right|+\frac{\widehat{y}_1(2\pi/w)}{\dot{\widetilde{x}}_1(0)}\mathcal{F}_2(s+\theta){\rm
  sign}\left(\widehat{\xi_2}\right)\right)
  \cos^2(w\theta)+\right.
$$
$$\left.+\widetilde{\xi}_1\sin^2(w\theta)
  \right|>
$$
\begin{equation}\label{ner}
>
  \max_{s,\theta\in[0,2\pi/w]}\left|
\left(\widehat{\xi}_1+\frac{\widehat{y}_1(2\pi/w)}{\dot{\widetilde{x}}_1(0)}\mathcal{F}_1(s+\theta)\right)
  \sin(w\theta)\cos(w\theta)\right|.
\end{equation}
 Из условий (\ref{positf}) и (\ref{posit}) следует,
что
\begin{equation}\label{positdx1}
\dot { \widetilde{x}}_1(t)>0\mbox{ при }t\in[0,\pi/(2w)).
\end{equation}
Из последнего, в частности, следует, что функция $\widetilde{x}_1$
строго возрастает на интервале $[0,\pi/(2w)),$ поэтому
\begin{equation}\label{positx1}
  \widetilde{x}_1(t)>0\mbox{ при }t\in(0,\pi/(2w)].
\end{equation}
Поскольку имеем (\ref{positg}), то оценки (\ref{posit}) и
(\ref{positx1}) влекут
\begin{equation}\label{positgg}
  g_1(\widetilde{x}(t))>0\mbox{ и }g_2(\widetilde{x}(t))>0\mbox{ при
  }t\in(0,\pi/(2w)).
\end{equation}
Наконец, из (\ref{positf}) и (\ref{positx1}) имеем
\begin{equation}\label{positdx2}
\dot { \widetilde{x}}_2(t)<0\mbox{ при }t\in(0,\pi/(2w)].
\end{equation}
Оценки (\ref{positgg}) и (\ref{positdx2}) позволяют заключить, что
\begin{equation}\label{positxi}
\widetilde{\xi}_1>0.
\end{equation}
В силу соотношений (\ref{N}), и учитывая $2\pi/w$-периодичность
подинтегральных функций в (\ref{N2}), получаем
$$
  \max_{t\in[0,4\pi/w]}\left|\mathcal{F}_2(t)\right|=\max_{t\in[2\pi/w,3\pi/w]}\left|\mathcal{F}_2(t)\right|.
$$
На основании соотношений (\ref{N1})
$$
  \max_{t\in[2\pi/w,3\pi/w]}\left|\mathcal{F}_2(t)\right|=\max_{t\in[4\pi/(2w),5\pi/(2w)]}\left|\mathcal{F}_2(t)\right|.
$$
Но из (\ref{positdx1}), (\ref{positgg}) и  (\ref{positdx2})
следует, что функции
$\tau\to-\dot{\widetilde{x}}_2(\tau)\sin(w\tau)g_1(\widetilde{x}(\tau))$
и $\tau\to\dot{\widetilde{x}}_1(\tau)
\sin(w\tau)g_2(\widetilde{x}(\tau))$ неотрицательны на
$[0,\pi/(2w)].$ Поэтому
$$
  \max_{t\in[0,4\pi/w]}\left|\mathcal{F}_2(t)\right|=\left|\mathcal{F}_2(5\pi/(2w))\right|=$$
  $$=
  \int\limits_{4\pi/(2w)}^{5\pi/(2w)}\left<
\left(\begin{array}{c} -\dot{\widetilde{x}}_2(\tau) \\
\dot{\widetilde{x}}_1(\tau)\end{array}\right),
g(\widetilde{x}(\tau))\right>\sin(w\tau)d\tau
$$
и, суммируя (\ref{N}) и (\ref{N1}), получаем
\begin{equation}\label{L}
  \max\limits_{t\in[0,4\pi/w]}\left|\mathcal{F}_2(t)\right|=\widetilde{\xi}_2/4.
\end{equation}
 Из условия (\ref{lincond}) теоремы и (\ref{L})
 следует, что
$$
  \left|\widehat{\xi}_2\right|>\max_{s,\theta\in[0,2\pi/w]}\left|\frac{\widehat{y}_1(2\pi/w)}
  {\dot{\widetilde{x}}_1(0)}\mathcal{F}_2(s+\theta)\right|.
$$
Поэтому имеем
$$
  \min_{s,\theta\in[0,2\pi/w]}\left|
\left(\left|\widehat{\xi}_2\right|-\frac{\widehat{y}_1(2\pi/w)}{\dot{\widetilde{x}}_1(0)}\mathcal{F}_2(s+\theta){\rm
  sign}\left(\widehat{\xi_2}\right)\right)
  \cos^2(w\theta)+\right.
$$
$$\left.+\widetilde{\xi}_1\sin^2(w\theta)
  \right|\ge
$$
$$
  \ge\min\left\{\min_{s,\theta\in[0,2\pi/w]}\left(
\left|\widehat{\xi}_2\right|+\frac{\widehat{y}_1(2\pi/w)}{\dot{\widetilde{x}}_1(0)}\mathcal{F}_2(s+\theta){\rm
  sign}\left(\widehat{\xi_2}\right)\right),\widetilde{\xi}_1
  \right\}\cdot
$$
$$\cdot\left(\cos^2(w\theta)+\sin^2(w\theta)\right)\ge
$$
$$
  \ge\min\left\{\min_{s,\theta\in[0,2\pi/w]}\left(
\left|\widehat{\xi}_2\right|-\frac{\left|\widehat{y}_1(2\pi/w)\right|}{\dot{\widetilde{x}}_1(0)}\left|\mathcal{F}_2(s+\theta)\right|
\right),\widetilde{\xi}_1
  \right\}=
$$
$$ =\min\left\{
\left|\widehat{\xi}_2\right|-\frac{\left|\widehat{y}_1(2\pi/w)\right|}{\dot{\widetilde{x}}_1(0)}\max_{s,\theta\in[0,2\pi/w]}\left|\mathcal{F}_2(s+\theta)\right|
,\widetilde{\xi}_1
  \right\}
$$
 Подставляя (\ref{L}) в полученное неравенство и используя
условие (\ref{lincond}), имеем
$$
  \min_{s,\theta\in[0,2\pi/w]}\left|
\left(\left|\widehat{\xi}_2\right|+\frac{\widehat{y}_1(2\pi/w)}{\dot{\widetilde{x}}_1(0)}\mathcal{F}_2(s+\theta){\rm
  sign}\left(\widehat{\xi_2}\right)\right)
  \cos^2(w\theta)+\right.
  $$
  $$\left.+\widetilde{\xi}_1\sin^2(w\theta)
  \right|>\left|\widehat{\xi_1}\right|+\dfrac{\left|\widehat{y}_1(2\pi/w)\right|}{\dot{\widetilde{x}}_1(0)}
\widetilde{\xi_1}.
$$
Таким образом, для завершения доказательства неравенства
(\ref{ner}) остается показать, что
$$
  \max_{s,\theta\in[0,2\pi/w]}\left|
\left(\widehat{\xi}_1+\frac{\widehat{y}_1(2\pi/w)}{\dot{\widetilde{x}}_1(0)}\mathcal{F}_1(s+\theta)\right)
  \sin(w\theta)\cos(w\theta)\right|\le $$
  $$\le\left|\widehat{\xi_1}\right|+\dfrac{\left|\widehat{y}_1(2\pi/w)\right|}{\dot{\widetilde{x}}_1(0)}
\widetilde{\xi_1},
$$
для чего, в свою очередь, достаточно установить, что
\begin{equation}\label{dost}
  \max_{s,\theta\in[0,2\pi/w]}\left|\mathcal{F}_1(s+\theta)\right|=2\widetilde{\xi}_1.
\end{equation}
Действительно, оценки (\ref{positdx1}), (\ref{positgg}) и
(\ref{positdx2})
 влекут
\begin{equation}\label{vvv} \left<
\left(\begin{array}{c} -\dot{\widetilde{x}}_2(\tau) \\
\dot{\widetilde{x}}_1(\tau)\end{array}\right),
g(\widetilde{x}(\tau))\right>\cos(w\tau)d\tau\ge 0,\quad
\tau\in[0,\pi/(2w)].
\end{equation}
Далее, на основании неравенств (\ref{O2}) множество значений
$\tau,$ при которых неравенство (\ref{vvv}) верно, можем расширить
c $[0,\pi/(2w)]$ до $[0,\pi/w]$ и затем, на основании неравенств
(\ref{O1}) с $[0,\pi/w]$ до $[0,2\pi/w].$ Наконец, в силу
$2\pi/w$-периодичности по $\tau$ левой части неравенства
(\ref{vvv}) заключаем, что оно верно при всех $\tau\in\mathbb{R}.$
Полученное свойство означает, что
$$
\max_{s,\theta\in[0,2\pi/w]}\left|\mathcal{F}_1(s+\theta)\right|=2\mathcal{F}_1(0),
$$
то есть  имеем (\ref{dost}) и, следовательно, справедливость
неравенства (\ref{ner}) установлена. Как отмечалось ранее,
справедливость этого неравенства означает, что выполнено
(\ref{gar}). Поэтому при любом $s\in [0,2\pi/w]$ оператор $\Phi^s$
невырожден на $\partial U,$ и
$$
  d_{\mathbb{R}^2}(-\Phi^s,U)=d_{\mathbb{R}^2}(\Phi^s,U)=d_{\mathbb{R}^2}(F,U)\in\{0,2\}.
$$

Таким образом выполнены все условия теоремы~2.5, применяя которую
получаем утверждение доказываемой теоремы.

Теорема доказана.

\begin{rem}
В условиях теоремы \ref{thm_sym} число $\widetilde{\xi}_1$
положительно (см. формулу~\ref{positxi} доказательства).
\end{rem}

\begin{rem}
В условиях теоремы \ref{thm_sym} функция $\widetilde{f}(\cdot,0)$
имеет вид
$$
  \widetilde{f}(\theta,0)=\int\limits_0^T\left\|\left(\dot{\widetilde{x}}(\tau),\sin(w(\tau-\theta))
  g(\widetilde{x}(\tau))\right)
  \right\|d\tau,\quad \theta\in[0,T]
$$
(см. формулу~\ref{STAR2} доказательства), то есть совпадает с
субгармонической порядка $1/1$ функцией Мельникова $M^{1/1}$
системы~(\ref{ps_1}) (см. \cite{mel}, формула для $A_0(v),$ c.~42
или Дж.~Гукенхеймер и Ф.~Холмс \cite{guck}, формула~4.6.2). В
частности, $M^{1/1}$ имеет ровно два простых нуля на интервале
$[0,2\pi/w).$
\end{rem}
Для того, чтобы заметить, что в условиях теоремы~2.7 функция
$M^{1/1}$ имеет ровно два простых нуля на $[0,2\pi/w),$ достаточно
обратиться к формуле (\ref{pok}) доказательства этой теоремы, из
которой следует, что
$M^{1/1}(\theta)=-\widetilde{\xi}_1\sin(w\theta).$

\begin{exa}\label{exa_green}  \rm Для иллюстрации работы теоремы \ref{thm_sym}
рассмотрим следующую модифицированную систему Гринспана-Холмса
(см. \cite{green})
\begin{equation}\label{green}
  \begin{array}{lll}
    \dot x_1&=&x_2(1-\delta(x_1^2+x_2^2))\\
    \dot
    x_2&=&-x_1(1-\delta(x_1^2+x_2^2))+\varepsilon\sin((1-\delta)t),
  \end{array}
\end{equation}
\end{exa}
на примере которой выяснение смысла условий (\ref{lincond})
представляется наиболее наглядным.

Легко проверить, что при $\delta\in(0,1)$ система (\ref{green})
удовлетворяет условиям (\ref{m1})-(\ref{K3}),
(\ref{STARR})-(\ref{positg}), и порождающая система
\begin{equation}\label{greennp}
  \begin{array}{lll}
    \dot x_1&=&x_2(1-\delta(x_1^2+x_2^2)),\\
    \dot
    x_2&=&-x_1(1-\delta(x_1^2+x_2^2))
  \end{array}
\end{equation}
допускает семейство циклов
$$
  x(t)=\left(\begin{array}{l} \alpha\sin((1-\delta\alpha^2)t) \\
  \alpha\cos((1-\delta\alpha^2)t)\end{array}\right),\quad\alpha>0.
$$
Рассмотрим задачу о возмущении цикла
\begin{equation}\label{orbit}
  \widetilde{x}(t)=\left(\begin{array}{l} \sin((1-\delta)t) \\
  \cos((1-\delta)t)\end{array}\right)
\end{equation}
периода $2\pi/(1-\delta),$ совпадающего с периодом возмущения.
Легко проверить, что цикл $\widetilde{x}$ удовлетворяет условиям
$(A_1),$ $(C),$ (\ref{init}) и (\ref{U}).
 Линеаризованная на
$\widetilde{x}$ система (\ref{greennp}) имеет вид
$$
\hskip-13cm\left(\begin{array}{c} \dot y_1 \\ \dot y_2
  \end{array}\right)=$$
  $$=\left(\begin{array}{cc} -2\delta\sin((1-\delta)t)\cos((1-\delta)t) & 1-\delta-2\delta\cos^2((1-\delta)t) \\
 -1+\delta+2\delta\sin^2((1-\delta)t) & 2\delta\sin((1-\delta)t)\cos((1-\delta)t)
 \end{array}\right)\circ
$$ $$\hskip13cm\circ
 \left(\begin{array}{c} y_1 \\ y_2
  \end{array}\right),
$$
и кроме $2\pi/(1-\delta)$-периодического решения
$$
  \dot{\widetilde{x}}(t)=(1-\delta)\left(\begin{array}{l} \cos((1-\delta)t) \\
 -\sin((1-\delta)t) \end{array}\right),
$$
удовлетворяющего начальному условию
$\dot{\widetilde{x}}(0)=(1-\delta,0),$ допускает следующее решение
$$
  \widehat{y}(t)=\frac{1}{1-\delta}\left(\begin{array}{c} -2\delta t \cos((1-\delta)t) +\sin((1-\delta)t)\\
 2\delta t \sin((1-\delta)t) +\cos((1-\delta)t)
 \end{array}\right),
$$
удовлетворяющее начальному условию
$\widehat{y}(0)=(0,1/(1-\delta)).$ После некоторых преобразований
для $\widetilde{\xi}$ и $\widehat{\xi}$ получаем следующие
выражения
$$
  \widetilde{\xi}=\left(\begin{array}{c} 4/3\\
4/3
 \end{array}\right),\qquad   \widetilde{\xi}=\left(\begin{array}{c} \dfrac{2\delta\pi^2}{(1-\delta)^3}\\
-\dfrac{\pi}{(1-\delta)^3}
 \end{array}\right),
$$
на основании которых условие (\ref{lincond}) теоремы~\ref{thm_sym}
для системы (\ref{green}) записывается в виде
$$
  \dfrac{2\delta\pi^2}{(1-\delta)^3}+\dfrac{16\delta\pi^2}{3(1-\delta)^2}<\min\left\{\dfrac{\pi}{(1-\delta)^3}
  -\dfrac{4\delta\pi}{3(1-\delta)^3},\dfrac{4}{3}\right\}.
$$
Нетрудно заключить, что последнему неравенству удовлетворяют все
$\delta>0,$ при которых
\begin{equation}\label{PU}
  2(1-\delta)^3-(3\pi^2+8\pi)\delta>0.
\end{equation}
Точное решение неравенства (\ref{PU}) может быть получено при
помощи известных формул Кардано для корней многочленов третьей
степени (см. \cite{mish}, Гл.~III, \S~3.2). В настоящей
диссертации мы ограничимся заключением о том, что неравенству
(\ref{PU}) удовлетворяет отрезок $[0,1/40],$ проверяемым
непосредственной постановкой чисел $\delta=0$ и $\delta=1/40$ в
данное неравенство. В любом случае, из неравенства (\ref{PU})
следует,  что условие (\ref{lincond}) теоремы~2.7 связано с
ограничением на нелинейность порождающей системы (\ref{greennp}).

Таким образом, теорема~\ref{thm_sym} позволяет получить для
системы (\ref{green}) следующее утверждение, более сильное, чем
утверждение леммы~\ref{duf_lem} о $2\pi/(1-\delta)$-периодических
решениях уравнения Дуффинга.

\begin{pro}\label{pro_green}
Пусть $
  \delta\in(0,1/40].$ Существует $\varepsilon_0>0$ такое, что:

1) при каждом $\varepsilon\in(0,\varepsilon_0]$ возмущенное
уравнение (\ref{green}) имеет по крайней мере два различных
$\frac{2\pi}{1+\delta}$-периодических решения
$\widetilde{x}_{\varepsilon},$
$\widetilde{\widetilde{x}}_{\varepsilon}$ таких, что значения
функции $\widetilde{x}_{\varepsilon}$ лежат строго внутри
единичного круга, а значения функции
$\widetilde{\widetilde{x}}_\varepsilon$ лежат строго снаружи
единичного круга;

2) всякое $\frac{2\pi}{1+\delta}$-периодическое решение $x$
системы (\ref{green}) с $\varepsilon\in(0,\varepsilon_0]$ таково,
что $\|x(t)\|\not=1$ при любом
$t\in\left[0,\frac{2\pi}{1+\delta}\right];$

3) решения $\widetilde{x}_{\varepsilon}$ и
$\widetilde{\widetilde{x}}_{\varepsilon}$
 удовлетворяют
условиям \\$
  \widetilde{x}_{\varepsilon}(t)\to
  \left(\sin\left((1-\delta)\left(t-\widetilde{\theta}\right)\right),
  \cos\left((1-\delta)\left(t-\widetilde{\theta}\right)\right)\right)
  $ и \\$
  \widetilde{\widetilde{x}}_{\varepsilon}(t)\to \left(\sin\left((1-\delta)
  \left(t-\widetilde{\widetilde{\theta}}\right)\right),
  \cos\left((1-\delta)\left(t-\widetilde{\widetilde{\theta}}\right)\right)\right)
  $ при $\varepsilon\to 0$
для некоторых
$\widetilde{\theta},\widetilde{\widetilde{\theta}}\in\left[0,\frac{2\pi}{1+\delta}\right].$
\end{pro}

\subsubsection{Случай, когда порождающий цикл является
вырожденным}

 В этом подпункте предполагается, что $\{x(\cdot,\alpha)\}_{\alpha>0}$ --
 семейство циклов порождающей системы (\ref{npa_1})
таких, что
\begin{equation}\label{ST1}
x(0,\alpha)=(0,\alpha),
\end{equation}
и векторы $x(0,\alpha)$ и $x'_{(1)}(0,\alpha)(0)$
линейно-независимы для всех $\alpha>0,$ то есть
\begin{equation}\label{ST2}
x(0,\alpha)\nparallel x'_{(1)}(0,\alpha)(0),\quad \alpha>0.
\end{equation}
Через $T(\alpha)$ обозначается наименьший период цикла
$x(\cdot,\alpha).$
\begin{opr}
Цикл $x(\cdot,\alpha_0)$ семейства
$\{x(\cdot,\alpha)\}_{\alpha>0},$ удовлетво-ряющий условиям
(\ref{ST1})-(\ref{ST2}), будем называть вырожденным, если
$$
T'(\alpha_0)=0.
$$
\end{opr}

Следующий пример показывает, что если цикл $x(\cdot,\alpha_0)$
является вырожденным, то применение теоремы~\ref{dimtwo_1}
значительно упрощается. Затем это наблюдение обосновывается, и
соответствующие упрощения теорем~\ref{dimtwo_1} и \ref{dimtwoF_1}
формулируются для общего вырожденного случая.

\begin{exa}\label{exa_mak}\rm
Рассмотрим систему
\begin{equation}\label{mak}
  \begin{array}{lll}
    \dot x_1&=&x_2\left(\left(\sqrt{x_1^2+x_2^2}-1\right)^2+1\right)\\
    \dot
    x_2&=&-x_1\left(\left(\sqrt{x_1^2+x_2^2}-1\right)^2+1\right)+\\
    & &+\varepsilon(\mu x^+_1+\nu
  x_1^-+\cos((1+\delta)t)).
  \end{array}
\end{equation}
\end{exa}
Порождающая система
\begin{equation}\label{maknp}
  \begin{array}{lll}
    \dot x_1&=&x_2\left(\left(\sqrt{x_1^2+x_2^2}-1\right)^2+1\right)\\
    \dot
    x_2&=&-x_1\left(\left(\sqrt{x_1^2+x_2^2}-1\right)^2+1\right)
  \end{array}
\end{equation}
допускает семейство циклов
$$
  x(t,\alpha)=\alpha\left(\begin{array}{c}
  \sin\left(\left((\alpha-1)^2+1\right)t\right) \\
  \cos\left(\left((\alpha-1)^2+1\right)t\right)
  \end{array}\right)
$$
периода
$$
  T(\alpha)=\frac{2\pi}{(\alpha-1)^2+1}.
$$
Таким образом, имеем
$$
  T'(1)=0
$$
и, значит, цикл $\widetilde{x}(t)=x(t,1)=(\sin t,\cos t)$ системы
(\ref{maknp}) является вырожденным. Оказывается, линеаризованная
на цикле $\widetilde{x}$ система (\ref{maknp})
 имеет вид
$$
 \left(\begin{array}{l} \dot y_1 \\
 \dot y_2 \end{array}\right)=\left(\begin{array}{cc} 0 & 1 \\
 -1 & 0 \end{array}\right)\left(\begin{array}{l} y_1 \\
  y_2 \end{array}\right)
$$
и, значит, решение $\widehat{y},$ участвующее в формуле для
обобщенного оператора усреднения $\Phi^s$ системы (\ref{mak}),
дается формулой
$$
   \widehat{y}(t)=\left(\begin{array}{l} \sin(t)\\
 \cos(t) \end{array}\right).
$$
Отсюда следует, что оператор $\Phi^s$ для системы (\ref{mak})
совпадает с оператором $\Phi^s$ для косинусоидально возмущенной
линейной системы (\ref{lin}), рассмотренной в
лемме~\ref{lem_nu_mu}. Получаем следующее утверждение

\begin{pro}\label{pro_mak} Пусть $|\mu-\nu|<2.$ Тогда существует
$\varepsilon_0>0$ такое, что:

1) при каждом $\varepsilon\in(0,\varepsilon_0]$ система
(\ref{mak}) имеет по крайней мере два $2\pi$-периодических решения
$\widetilde{x}_{\varepsilon},$
$\widetilde{\widetilde{x}}_{\varepsilon}$ таких, что значения
функции $\widetilde{x}_{\varepsilon}$ лежат строго внутри
единичного круга с центром в нуле, а значения функции
$\widetilde{\widetilde{x}}$ лежат строго снаружи этого круга;

2) всякое $2\pi$-периодическое решение $x$ системы (\ref{mak}) с
$\varepsilon\in(0,\varepsilon_0]$ таково, что $\|x(t)\|\not=1$ при
любом $t\in\left[0,2\pi\right];$

3) решения $\widetilde{x}_{\varepsilon}$ и
$\widetilde{\widetilde{x}}_{\varepsilon}$
 удовлетворяют
условиям \\$
  \widetilde{x}_{\varepsilon}(t)\to
  \left(\sin\left(t-\widetilde{\theta}\right),
  \cos\left(t-\widetilde{\theta}\right)\right)
  $ и \\ $
  \widetilde{\widetilde{x}}_{\varepsilon}(t)\to \left(\sin\left(t-\widetilde{\widetilde{\theta}}\right),
  \cos\left(t-\widetilde{\widetilde{\theta}}\right)\right)
  $ при $\varepsilon\to 0$
для некоторых
$\widetilde{\theta},\widetilde{\widetilde{\theta}}\in\left[0,2\pi\right].$
\end{pro}

В общем же вырожденном случае справедливо следующее утверждение.

\begin{lem}\label{lem_per_sol} Если цикл $x(\cdot,\alpha_0)$ является вырожденным, то каждое решение
линеаризованной системы (\ref{lsg}) с
$\widetilde{x}:=x(\cdot,\alpha_0)$ является
$T(\alpha_0)$-периодическим.
\end{lem}

\noindent {\bf Доказательство.} Так как правая часть порождающей
системы (\ref{npa_1}) непрерывно дифференцируема, то (см.,
например, Л.~С.~Понтрягин \cite{pont}, Гл.~4, \S~24) функция
$(t,\alpha)\to x(t,\alpha)$ непрерывно дифферен-цируема по
совокупности переменных. В пределах данного доказательства через
$x'_t$ и $x'_\alpha$ обозначаются производные функции $x$ по
первой и второй переменным соответственно. Дифференцируя тождество
$$x'_t=f(x(t,\alpha))$$
по $\alpha,$ получаем
$$x'_t {}'_\alpha=f'(x(t,\alpha))x'_\alpha,$$
следовательно, $\widehat{y}=x'_\alpha(\cdot,\alpha_0)$ является
решением линеаризованной системы (\ref{lsg}) c
$\widetilde{x}=x(\cdot,\alpha_0).$ Имеем
\begin{eqnarray}
  \widehat{y}(T(\alpha_0))-\widehat{y}(0)&=&\lim_{\Delta\to 0}\frac{x(T(\alpha_0),\alpha_0+\Delta)-
  x(T(\alpha_0),\alpha_0)}{\Delta}-\nonumber\\
  & & -\lim_{\Delta\to
  0}\frac{x(0,\alpha_0+\Delta)-x(0,\alpha_0)}{\Delta}=\nonumber\\
  & =&\lim_{\Delta\to
  0}\frac{x(T(\alpha_0),\alpha_0+\Delta)-x(0,\alpha_0+\Delta)}{\Delta}+\nonumber\\
  & &
+\frac{-  x(T(\alpha_0),\alpha_0)+x(0,\alpha_0)}{\Delta}
  =\nonumber\\
  &=&\lim_{\Delta\to
  0}\frac{x(T(\alpha_0),\alpha_0+\Delta)-x(0,\alpha_0+\Delta)}{\Delta}=\nonumber\\
  &=&\left(x(T(\cdot),\alpha_0)\right)'(\alpha_0)=x'_t(T(\alpha_0),\alpha_0)T'(\alpha_0)=0,\nonumber
\end{eqnarray}
то есть $\widehat{y}$ -- $T(\alpha_0)$-периодическое решение
системы (\ref{lsg}). Но
$$
  \widehat{y}(0)=\lim_{\Delta\to
  0}\frac{x(0,\alpha_0+\Delta)-x(0,\alpha_0)}{\Delta}=\lim_{\Delta\to
  0}\left(\begin{array}{c} 0\\1\end{array}\right)=\left(\begin{array}{c}
  0\\1\end{array}\right),
$$
следовательно, $\widehat{y}$ и $x'_t(\cdot,\alpha_0)$ -- два
линейно-независимых $T(\alpha_0)$-периодических решений
(двумерной) системы (\ref{lsg}).

Лемма доказана.

Из леммы~\ref{lem_per_sol} следует, что если цикл
$\widetilde{x}=x(\cdot,\alpha_0)$ является вырожденным, то формула
(\ref{Ftheta}) для оператора $\Phi^s$ принимает значительно более
простой вид:
\begin{equation}\label{Fthetas}
\Phi^s(\widetilde {x}(\theta))=\widehat{f}(\theta)\dot{\widetilde{
x}}(\theta)+ \widetilde{f}(\theta,0)\widehat {y}(\theta) \quad
\mbox{для любых }s,\theta\in\mathbb{R},
\end{equation} где
$$
   \widetilde{f}(\theta,0)=\int\limits_0^{T}
\left<\widetilde{z} (\tau),g(\tau-\theta,\widetilde{
  x}(\tau),0)\right>d\tau,$$
  $$
   \widehat{f}(\theta)=\int\limits_0^{T}
\left<\widehat{z} (\tau),g(\tau-\theta,\widetilde{
  x}(\tau),0)\right>d\tau.
$$

Соответственно, получаем следующие следствия из
теорем~\ref{dimtwoF_1} и \ref{thm_sym} для случая вырожденных
циклов.

\begin{cor}\label{cor_dimtwo} Пусть выполнены условия $(A_1)$ и $(C).$
 Пусть $T$-периодический цикл $\widetilde{x}=x(\cdot,\alpha_0)$ является вырожденным.
Предположим, что для каждого
$\theta_0\in[0,\widetilde {T}]$ такого, что
$\widetilde{f}(\theta_0,0)=0,$ имеем
$$
  \widehat{f}(\theta_0)\not=0.
$$
Тогда существует $\varepsilon>0$ такое, что при каждом
$\varepsilon\in(0,\varepsilon_0]$ всякое $T$-периодическое решение
$x$ системы (\ref{psa_1}) таково, что $x(t)\not\in\partial U$ при
любом $t\in\left[0,T\right],$ где $U\subset\mathbb{R}^2$ --
внутренность цикла $\widetilde{x}.$ Более того, если известно, что
$$
  d_{\mathbb{R}^2}(-\Phi^{T},U)\not=1,
$$
то при всех достаточно малых $\varepsilon>0$ система (\ref{psa_1})
имеет по крайней мере два $T$-периодических решения
$\widetilde{x}_\varepsilon$ и
$\widetilde{\widetilde{x}}_\varepsilon$ таких, что
$\widetilde{x}_\varepsilon(t)\in U,$
$\widetilde{\widetilde{x}}_\varepsilon(t)\not\in U$ для любого
$t\in[0,T],$ и
$$\rho\left(\widetilde{x}_\varepsilon(t),\partial U\right)+
\rho\left(\widetilde{\widetilde{x}}_\varepsilon(t),\partial U
\right)\to 0\quad\mbox{при }\varepsilon\to 0
$$
равномерно по $t\in[0,T].$
\end{cor}

Напомним, что число $\widetilde{\xi}_1\in\mathbb{R}$ и вектор
$\widehat{\xi}\in\mathbb{R}^2$ в теореме~\ref{thm_sym}  введены
как
\begin{eqnarray}
  \widetilde{\xi}_1 &=& 4\int\limits_0^{\pi/(2w)}\left<
\left(\begin{array}{c} -\dot{\widetilde{x}}_2(\tau) \\
\dot{\widetilde{x}}_1(\tau)\end{array}\right),
  g(\widetilde{x}(\tau))\right>\cos(w\tau)d\tau,\nonumber\\
  \widehat{\xi} &=&\int\limits_0^{2\pi/w}\left<
\left(\begin{array}{c} \widehat{y}_2(\tau) \\
-\widehat{y}_1(\tau)\end{array}\right),
  g(\widetilde{x}(\tau))\right>\left(\begin{array}{c}\cos(w\tau) \\
  \sin(w\tau)\end{array}\right)d\tau.\nonumber
\end{eqnarray}

\begin{cor}\label{cor_dimtwoF}

Пусть $\widetilde{x}=x(\cdot,\alpha_0)$ -- вырожденный
периодический цикл порождающей системы
 (\ref{npa_1})
наименьшего периода $2\pi/w,$ удовлетворяющий условиям $(A_1),$
$(C),$ (\ref{init})   и (\ref{U}). Пусть выполнены условия
симметрии (\ref{m1})-(\ref{K3}), (\ref{STARR})-(\ref{positg}).
 Тогда, если
\begin{equation}\label{lincond1}
  \left|\widehat{\xi}_1\right|<
   \min\left\{\left|\widehat{\xi}_2\right|,
\left|\widetilde{\xi}_1 \right|\right\},
\end{equation}
то при всех достаточно малых $\varepsilon>0$ синусоидально
возмущенная система (\ref{ps2}) имеет по крайней мере два
$2\pi/w$-периодических решения $\widetilde{x}_\varepsilon,$
$\widetilde{\widetilde{x}}_\varepsilon$
 таких, что
$\widetilde{x}_\varepsilon(t)\in U,$
$\widetilde{\widetilde{x}}_\varepsilon(t)\not\in U$ для любого
$t\in[0,2\pi/w],$ и
$$\rho\left(\widetilde{x}_\varepsilon(t),\partial U\right)+
\rho\left(\widetilde{\widetilde{x}}_\varepsilon(t),\partial U
\right)\to 0\quad\mbox{при }\varepsilon\to 0
$$
равномерно по $t\in[0,2\pi/w].$ Прочие $2\pi/w$-периодические
решения $x$ системы (\ref{ps2}) при
$\varepsilon\in(0,\varepsilon_0]$ удовлетворяют условию
$x(t)\not\in\partial U$ для любого $t\in\left[0,2\pi/w\right].$
\end{cor}

\section{Сопоставление полученных результатов с имеющимися в
литературе}

В силу леммы~\ref{form_eta} для обобщенного оператора сдвига
$\Phi^s$ системы (\ref{ps_1})  справедливо соотношение
$$
  \eta(T,s,\xi)-\eta(0,s,\xi)=\int\limits_{s-T}^s
  \Omega'_\xi(0,\tau,\Omega(\tau,0,\xi))g(\tau,\Omega(\tau,0,\xi),0)d\tau=
$$
\begin{equation}\label{soo}
  =\int\limits_{s-T}^s
  Y^{-1}(\tau,\xi)g(\tau,\Omega(\tau,0,\xi),0)d\tau,
\end{equation}
где $Y(t,\xi)$ -- нормированная $(Y(0,\xi)=I)$ фундаментальная
матрица системы $
  \dot y=f'_{(2)}(t,\Omega(t,0,\xi))y.
$ Из (\ref{soo}) следует, что при $f=0$ оператор $\Phi^s$ имеет
вид
$$
  \Phi^s(\xi)=\int\limits_0^T g(\tau,\xi,0)d\tau,
$$
то есть не зависит от $s$ и совпадает с классическим оператором
усреднения Крылова-Боголюбова-Митропольского, см.~\cite{bog},
формула~7.112.
 Соответствующая теорема принципа
усреднения Крылова-Боголюбова-Митропольского о $T$-периодических
решениях системы (\ref{ps_1}) была впервые сформулирована в
терминах теории топологической степени Ж.~Мавеном в его
диссертационной работе \cite{mawphd1} (позже опубликована в
\cite{mawphd2}) и утверждает, что если $
 \int\limits_0^T g(\tau,\xi,0)d\tau\not=0
$ для всех $\xi\in\partial U,$ то при достаточно малых
$\varepsilon>0$  оператор
$$(Q_\varepsilon x)(t)=x(T)+\varepsilon\int\limits_0^T g(\tau,\xi,\varepsilon)d\tau$$
не имеет неподвижных точек на $\partial W_U,$ и
$
  d(I-Q_\varepsilon,W_U)=d_{\mathbb{R}^n}\left(-\int\limits_0^T g(\tau,\cdot,0)d\tau,U\right).
$ Таким образом, доказанная в настоящей главе
теорема~\ref{fromnon} является обобщением указанной теоремы
Мавена.

В случае, когда порождающая система (\ref{np_1}) имеет
невырожденный цикл $\widetilde{x}$ периода ${T},$
теорема~\ref{dimtwoF_1} дополняет классический метод
В.~К.~Мельникова (\cite{mel}, лемма~7), который утверждает, что
{\it если так называемая функция Мельникова
$$
  M^{1/1}(\theta)=\int\limits_0^{T}\left\|\left(\dot{\widetilde{x}}(\tau),g(\tau-\theta,\widetilde{x}(\tau),0)\right)\right\|d\tau
$$
имеет нуль $\theta_0\in[0,T]$ такой, что
$\left(M^{1/1}\right)'(\theta_0)\not=0,$ то при всех достаточно
малых $\varepsilon>0$ возмущенная система (\ref{ps_1}) допускает
$T$-периодическое решение $x_\varepsilon,$ удовлетворяющее условию
$$
  x_\varepsilon(t)\to\widetilde{x}(t+\theta_0)\quad \mbox{при}\
  \varepsilon\to 0.
$$}
Таким образом, теорема~\ref{dimtwoF_1} дает условия, при которых
$T$-периодические решения Мельникова не пересекают порождающий
цикл. Более того, траектории различных $T$-периодических решений
Мельникова могут совпадать, в то время как теорема~\ref{dimtwoF_1}
гарантирует существование для системы (\ref{ps_1}) по крайней мере
двух $T$-периодических решений с различными траекториями.

Сказанное позволило установить (лемма~\ref{duf_lem}), что
периодические решения
 возмущенного уравнения Дуффинга не пересекают порождающих циклов малой амплитуды.
 Такое свойство отсутствует в  классических результатах о
 периодических решениях уравнения Дуффинга, полученных
А.~Д.~Морозовым \cite{mord} и Б.~Гринспаном-Ф.~Холмсом
\cite{green}. Из доказательства леммы~\ref{duf_lem} следует, что
ее утверждения справедливы также для уравнения
$$
  \ddot u+u+u^3=\varepsilon(\mu x^+_1+\nu
  x_1^-+\cos((1+\delta)t)),
$$
где $|\mu-\nu|<2,$ к которому в силу недифференцируемости правой
части метод Мельникова не применим. Также установлены аналогичные
свойства $T$-периодических решений для системы Гринспана-Холмса
(пример~\ref{exa_green}), что стало возможным благодаря
предложенной для симметричного случая теореме~\ref{thm_sym}.
Отметим (см. замечание~2.2), что в случае симметричных систем,
рассматриваемых в теореме~\ref{thm_sym}, соответствующая функция
Мельникова $M^{1/1}(\theta)$ имеет ровно два простых нуля на
интервале $[0,T).$ Это означает, что в общих условиях
теоремы~\ref{thm_sym} и указанной выше теоремы Мельникова,
последняя теорема всегда гарантирует существование для системы
(\ref{ps_1}) точно такого же количества периодических решений
(двух) вблизи порождающего цикла $\widetilde{x},$ что и
теорема~\ref{thm_sym}, но теорема~\ref{thm_sym} дополнительно
утверждает, что траектории полученных решений не пересекаются.

Наконец отметим, что в рассмотренном в пункте~2.4 случае
вырожденного порождающего цикла метод Мельникова не работает, а
имеющиеся его модификации требуют выполнения целого ряда
дополнительных условий, см. К.~Йагасаки (\cite{yaga},
теорема~3.5). Условия же предложенной теоремы~\ref{dimtwoF_1}, как
продемонстрировано в примере~\ref{exa_mak}, в вырожденном случае
наоборот упрощаются, см. также соответствующие
следствия~\ref{cor_dimtwo} и \ref{cor_dimtwoF} из
теорем~\ref{dimtwoF_1} и \ref{thm_sym}. Для изучения существования
в возмущенной системе (\ref{psa_1}) периодических решений близких
к вырожденным циклам порождающей системы может, вообще говоря,
использоваться общая теорема Рума-Чиконе (\cite{chic},
теорема~4.1), но она работает только в случае, когда возмущение
зависит от фазовой переменной (см. \cite{chic}, формула~2.7), что
не требуется в указанных утверждениях пункта~2.4. Другие
качественные результаты о поведении периодических решений
возмущенных систем вблизи вырожденного порождающего цикла получены
А.~Д.~Морозовым и Л.~П.~Шильниковым в \cite{mor}.

Обсудим кратко публикации автора по результатам настоящей главы.
 Обобщенный оператор усреднения $\Phi^s$ предложен
М.~И.~Каменским, О.~Ю.~Макаренковым и П.~Нистри в \cite{kmndan}.
Лемма~\ref{form_eta} о виде оператора $\Phi^s$ доказана автором в
\cite{kmnnon}, им же в \cite{kmnnon} проведено доказательство
обобщенной формулы Мавена (утверждение~1 теоремы~\ref{fromnon}),
причем в несколько расширенной формулировке, чем данное в
настоящей главе. Доказательство утверждения~1
теоремы~\ref{fromnon} для частного класса систем (\ref{ps2}), но
дающее явное значение $\varepsilon_0,$ сделано
 в \cite{vest}. Теорема~\ref{fromdan}, связанная с приложением
обобщенной формулы Мавена к существованию периодических решений,
предложенная автором, опубликована в \cite{kmndan}.
Теорема~\ref{dimtwo_1} и формула~(\ref{Ftheta})
леммы~\ref{lem_Ftheta}, составляющие геометрический подход в
решении задачи В.~К.~Мельникова, опубликованы в \cite{barcelona}.
Доказательство утверждения~1) леммы~\ref{duf_lem} о существовании
периодических решений в уравнении Дуффинга для случая несколько
более общего уравнения сделано в \cite{ndes}.

\chapter{Скорость сходимости полученных $T$-периодических решений
при уменьшении амплитуды возмущения}\label{sec_conv}

\setcounter{subsection}{0}

Как отмечалось во введении, результаты о
 существовании
$T$-периодических решений  в $T$-периодических системах
обыкновенных дифференциальных уравнений
\begin{equation}\label{ps_3}
   \dot x=f(t,x)+\varepsilon g(t,x,\varepsilon),
\end{equation}
где $\varepsilon>0$ --  малый параметр, основанные на
геометрических методах, используют всего лишь непрерывность
функции $g$ (от $f$ может требоваться непрерывная
дифференцируемость), см., например, \cite{kam}, \cite{alex3},
\cite{alex2}, \cite{kraper}, \cite{krastr}, \cite{mit},
\cite{muh3}-\cite{muh1}, \cite{sam}, \cite{camaza}-\cite{dancer},
\cite{zan}-\cite{laz}, \cite{mawphd1}-\cite{ortega},
 а также результаты глав~1 и 2.

 В настоящей главе
изучается задача, поставленная Дж.~Хейлом и П.~Тбоас в
\cite{hale}, о скорости сходимости и поведении $T$-периодических
решений системы (\ref{ps_3}) при $\varepsilon\to 0$ для непрерывно
дифференцируемой функции
$f:\mathbb{R}\times\mathbb{R}^n\to\mathbb{R}^n$ и непрерывной
функции
$g:\mathbb{R}\times\mathbb{R}^n\times[0,1]\to\mathbb{R}^n.$

Особенность исследования данной задачи состоит в том, что во
многих случаях неизвестно, что порождающее $T$-периодическое
решение $\widetilde{x}$ системы
\begin{equation}\label{np_3}
   \dot x=f(t,x)
\end{equation}
удовлетворяет условиям невырожденности, то есть, что
линеаризованная система
\begin{equation}\label{ls_3}
  \dot y=f'_{(2)}(t,\widetilde{x}(t))
\end{equation}
не имеет мультипликаторов $+1.$ Более того, как следует из
\cite{hale}, случай, когда $+1$ является мультипликатором системы
(\ref{ls_3}) представляет особый интерес.

Пусть $\{\varepsilon_k\}_{k\in\mathbb{N}}$ -- сходящаяся к нулю
последовательность значений параметра системы (\ref{ps_3}) и
$\{x_k\}_{k\in\mathbb{N}}$ -- соответствующая последовательность
$T$-периодических решений этой системы такая, что
\begin{equation}\label{as1_3}
  x_k(t)\to\widetilde{x}(t)\quad\mbox{при}\ \ k\to \infty,
\end{equation}
где $\widetilde{x}$ -- $T$-периодическое решение порождающей
системы (\ref{np_3}).

\section{Одна альтернатива для общего случая}

 Следующая альтернатива утверждает, что
 либо
начальные условия $x_k(0)$  сходятся  к начальному условию
$\widetilde{x}(0)$ порождающего решения $\widetilde{x}$ вдоль
плоскости\linebreak
$\left\{l\in\mathbb{R}^n:\left(\Omega'_{(3)}(T,0,\widetilde{x}(0))-I\right)l=0\right\},$
либо сходимость имеет скорость $\varepsilon_k>0.$ При этом, в
последнем случае описание поведения решений $x_k$ при $k\to\infty$
может быть уточнено на основании обобщенного оператора усреднения
$\Phi^s,$ соответствующего задаче о $T$-периодических решениях для
системы (\ref{ps_3}).

\begin{thm}\label{alter} Пусть выполнено условие (\ref{as1_3}) и
\begin{equation}\label{as2_3}
  \frac{\widetilde{x}(0)-x_k(0)}{\left\|\widetilde{x}(0)-x_k(0)\right\|}\to l\quad\mbox{при}\ \ k\to \infty,
\end{equation}
где $l\in\mathbb{R}^n.$ Тогда либо
\begin{equation}\label{alt1_3}
\left(\Omega'_{(3)}(T,0,\widetilde{x}(0))-I\right)l=0,
\end{equation}
либо существует константа $c>0$ такая, что
\begin{equation}\label{alt2_3}
  \left\|\widetilde{x}(0)-x_k(0)\right\|\le
  c\varepsilon_k,\quad k\in\mathbb{N}.
\end{equation}
В  последнем случае для любой сходящейся последовательности $
\left\{\dfrac{\widetilde{x}(0)-x_k(0)}{\varepsilon_k}\right\}_{k\in\mathbb{N}}
$ последовательность $
\left\{\dfrac{\widetilde{x}(0)-\Omega(0,t,x_k(t))}{\varepsilon_{k}}
\right\}_{k\in\mathbb{N}}$ также сходится, и справедливо
соотношение
$$
 (\Omega'_{(3)}(T,0,\widetilde{x}(0))-I)\lim\limits_{k\to\infty}\frac{\widetilde{x}(0)-\Omega(0,t,x_{k}(t))}{\varepsilon_{k}}=
$$
\begin{equation} \label{claa_3}
  \hskip1cm=\Phi^t(\widetilde{x}(0)),\quad t\in
 [0,T].
\end{equation}
\end{thm}

\noindent {\bf Доказательство.}  Положим
\begin{equation} \label{4_3}
 \nu_k(t)=\Omega(0,t,x_k(t)).
\end{equation}
Тогда, в силу леммы~\ref{zamena}
\begin{equation}\label{10_3}
  \nu_k(t)=\Omega(T,0,\nu_k(T))+\varepsilon_k\int\limits_0^t\Upsilon_{\varepsilon_k}(\tau,\nu_k(\tau))d\tau,
\end{equation}
где
$$
 \Upsilon_\varepsilon(t,\xi)={\Omega'}_{\xi}(0,t,\Omega(t,0,\xi))g(t,\Omega(t,0,\xi),\varepsilon).
$$
 Заметим, что
$\widetilde{x}(0)=\lim_{n\to\infty}\nu_k(0)=\lim_{n\to\infty}\nu_k(T)$
и из (\ref{10_3}) при $k=\infty$ получаем
$\widetilde{x}(0)=\Omega(T,0,\widetilde{x}(0)).$ Перепишем
(\ref{10_3}) в следующем виде
\begin{eqnarray}\label{10v_3}
  -(\widetilde{x}(0)-\nu_k(t))&=&\left(\Omega(T,0,\nu_k(t))-\Omega(T,0,\widetilde{x}(0))\right)+\nonumber\\
  & &+\Omega(T,0,\nu_k(T))-\Omega(T,0,\nu_k(t))
  +\nonumber\\
  & &+\varepsilon_k\int\limits_0^t\Upsilon_{\varepsilon_k}(\tau,\nu_k(\tau))d\tau\nonumber
\end{eqnarray}
или
\begin{eqnarray}\label{10vv_3}
  -(\widetilde{x}(0)-\nu_k(t))&=&-\Omega'_{(3)}(T,0,\widetilde{x}(0))\left(\widetilde{x}(0)-\nu_k(t)\right)+o\left(\widetilde{x}(0)-\nu_k(t)\right)\nonumber\\
  & &+\Omega'_{(3)}(T,0,\nu_k(T))\left(\nu_k(T)-\nu_k(t)\right)+\nonumber\\
  & &+o\left(\nu_k(T)-\nu_k(t)\right)
  +\varepsilon_k\int\limits_0^t\Upsilon_{\varepsilon_k}(\tau,\nu_k(\tau))d\tau.
\end{eqnarray}
Существует две возможности: либо
\begin{equation}\label{caI_3}
\dfrac{\left\|\widetilde{x}(0)-x_k(0)\right\|}{\varepsilon_k}\to\infty\quad\mbox{при}\
k\to\infty,
\end{equation}
либо существует $c>0$ такое, что
\begin{equation}\label{caII_3}
\dfrac{\left\|\widetilde{x}(0)-x_k(0)\right\|}{\varepsilon_k}<c\quad\mbox{при}\
k\in\mathbb{N}.
\end{equation}
В случае (\ref{caI_3}) имеем
$$
\dfrac{\varepsilon_k}{\left\|\widetilde{x}(0)-x_k(0)\right\|}\to 0
\quad\mbox{при}\ k\in\mathbb{N}.
$$
Поэтому, учитывая, что
\begin{equation}\label{uchII_3}
\nu_k(T)-\nu_k(t)=\varepsilon_k\int\limits_t^T\Upsilon_{\varepsilon_k}(\tau,\nu_k(\tau))d\tau,
\end{equation}
вправе перейти к пределу  при $k\to\infty$ и $t=0$ в
(\ref{10vv_3}), деленном на
$\left\|\widetilde{x}(0)-x_k(0)\right\|,$ и получить
$$
  \left(\Omega'_{(3)}(T,0,\widetilde{x}(0))-I\right)l=0.
$$
Таким образом, в случае (\ref{caI_3}) выполнено утверждение
(\ref{alt1_3}) теоремы~\ref{alter}, а в противном случае --
утверждение (\ref{caII_3}). Значит альтернатива
теоремы~\ref{alter} справедлива.

Предположим теперь, что последовательность $
\left\{\dfrac{\widetilde{x}(0)-x_{k}(0)}{\varepsilon_k}\right\}_{k\in\mathbb{N}}
$ сходится. Так как
\begin{eqnarray}
  \frac{\widetilde{x}(0)-\nu_k(t)}{\varepsilon_k}&=&\frac{\widetilde{x}(0)-x_k(0)+\nu_k(0)-\nu_k(t)}{\varepsilon_k}=\nonumber\\
  &=&\frac{\widetilde{x}(0)-x_k(0)}{\varepsilon_k}-\int\limits_0^t\Upsilon(\tau,\nu_k(\tau))d\tau,\label{uchI_3}
\end{eqnarray}
то последовательность $
\left\{\dfrac{\widetilde{x}(0)-\Omega(0,t,x_{k}(t))}{\varepsilon_{k}}
\right\}_{k\in\mathbb{N}}$ также сходится. Полагая
$h(t)=\lim\limits_{k\to\infty}\dfrac{\widetilde{x}(0)-\Omega(0,t,x_{k}(t))}{\varepsilon_{k}}$
и учитывая формулу (\ref{uchII_3}), перейдем к пределу в
(\ref{10vv_3}), деленном на $\varepsilon_k,$ в результате получим
$$
  -h(t)=-\Omega'_{(3)}(T,0,\widetilde{x}(0))h(t)+\Omega'_{(3)}(T,0,\widetilde{x}(0))
  \int\limits_t^T\Upsilon(\tau,\widetilde{x}(0))d\tau+$$
  $$+\int\limits_0^t\Upsilon(\tau,\widetilde{x}(0))d\tau,
$$
что в силу леммы~\ref{form_eta} означает (\ref{claa_3}).

Теорема доказана полностью.

В случае, когда порождающая система (\ref{np_3}) автономна, для
$T$-периодических решений возмущенной системы
\begin{equation}\label{psa_3}
   \dot x=f(x)+\varepsilon g(t,x,\varepsilon)
\end{equation}
получаем следующее следствие из теоремы~\ref{alter}.

\begin{cor}\label{cor_3}
Пусть выполнено условие (\ref{as1_3}) и цикл $\widetilde{x}$
является простым. Пусть выполнено условие (\ref{as2_3}) и $l$ --
вектор, о котором говорится в этом условии. Тогда, либо
$l=\lambda\dot{\widetilde{x}}(0),$ где $\lambda\not=0,$ либо
существует константа $c>0$ такая, что
\begin{equation}\label{cla_3}
  \left\|\widetilde{x}(t)-x_k(t)\right\|\le
  c\varepsilon_k,\qquad k\in\mathbb{N}.
\end{equation}
\end{cor}

В следующем пункте главы предлагается оценка для расстояния между
траекториями решений $x_k$ и $\widetilde{x},$ но при этом
дополнительно предполагается, что цикл $\widetilde{x}$ простой.

\section{Оценка скорости сходимости для случая, когда предельное
$T$-периодическое решение является простым циклом }

В настоящем пункте предполагается, что цикл $\widetilde{x}$
является простым. В сделанном предположении сопряженная система
\begin{equation}\label{ss}
  \dot z=-(f'(\widetilde{x}(t)))^*z
\end{equation}
допускает $n-1$ линейно-независимых не $T$-периодических решений
$z_1,z_2,...,z_{n-1},$ причем указанные решения всегда могут быть
выбраны так, что (см. \cite{mal}, формула~2.13)
\begin{equation}\label{orient0}
  \left<\dot{\widetilde{x}}(0),z_i(0)\right>=0,\qquad
  t\in\mathbb{R},\ i\in\overline{1,n-1}.
\end{equation}
Пусть  $Z_{n-1}(t)$ -- $n\times{n-1}$-матрица задаваемая формулой
$Z_{n-1}(t)=(z_1(t),z_2(t),...,z_{n-1}(t)).$
 Введем функцию $M^\bot:\mathbb{R}\to\mathbb{R}^{n-1}$ как
$$
  M^\bot(s)=\int\limits_{s-T}^s Z_{n-1}^*(\tau)g(\tau,\widetilde{x}(\tau),0)d\tau.
$$

Имеет место следующий результат.

\begin{thm}\label{conv_thm1}
Пусть выполнено условие (\ref{as1_3}). Тогда для любого
  $\theta\in[0,T]$
  имеем
\begin{equation}\label{mainprop}
  Z^*_{n-1}(\theta)\left(x_k(\theta-\Delta_{\varepsilon_k}(\theta))-\widetilde{x}(\theta)\right)
=\varepsilon_k D M^\bot(\theta)+o(\varepsilon_k),
\end{equation}
где $D$ -- невырожденная $(n-1)\times(n-1)$-матрица и
$\Delta_{\varepsilon_k}(\theta)\to 0$ при $k\to \infty$ равномерно
по отношению к $\theta\in[0,T].$ Более того, $\
x_k(\theta-\Delta_{\varepsilon_k}(\theta))\in
I(\theta,B_1^{n-1}(0)),$ $\theta\in[0,T],$ где
$I(\theta,\cdot):\mathbb{R}^{n-1}\to\mathbb{R}^n$ -- гладкая
поверхность, трансверсально пересекающая цикл $\widetilde{x}$ в
точке $\widetilde{x}(\theta).$
\end{thm}

Для доказательства теоремы~\ref{conv_thm1} определим, во-первых,
функцию $\Delta_{\varepsilon_k}(\theta),$ поверхность
$I(\theta,\cdot)$ и установим некоторые их свойства.

Обозначим через $Y_{n-1}(t)=(y_1(t),...,y_{n-1}(t))$ первые $n-1$
столбцов матрицы
$\left((Z_{n-1}(t),\widetilde{z}(t))^{-1}\right)^*,$ где
$\widetilde{z}$ -- $T$-периодическое решение системы (\ref{ss}),
удовлетворяющее
\begin{equation}\label{orient}
\left<\dot{\widetilde{x}}(t),\widetilde{z}(t)\right>=1,\quad
t\in\mathbb{R}.
\end{equation}
Последний выбор возможен (см. \cite{mal}, формула~2.13).

\begin{lem}\label{conv_lem1}
Предположим (\ref{orient0}) и (\ref{orient}). Тогда

1)
$\left((Z_{n-1}(t),\widetilde{z}(t))^{-1}\right)^*=\left(Y_{n-1}(t),\dot{\widetilde{x}}(t)\right),\quad
t\in\mathbb{R};$

2) любая функция из $\{y_1,...,y_{n-1}\}$ не $T$-периодична;

3) $Z^*_{n-1}(t)=\widetilde{D} Z^*_{n-1}(t+T),\quad
t\in\mathbb{R},$

где $\widetilde{D}$ -- постоянная $(n-1)\times (n-1)$-матрица,
собственные значения которой отличны от $+1.$

\end{lem}

\noindent {\bf Доказательство.} Так как
$\left((Z_{n-1}(t),\widetilde{z}(t))^{-1}\right)^*$ --
фундаментальная матрица $T$-периодической линейной системы $\dot
y=f'(\widetilde{x}(t))y$ (см. \cite{dem}, Гл.~III, лемма~\S~12),
то, пользуясь теорией Флоке, возможно записать следующую формулу
(см. \cite{dem}, Гл.~III, \S 15)
\begin{equation}\label{teflo}
\left((Z_{n-1}(t),\widetilde{z}(t))^{-1}\right)^*=\Phi(t)\left(\begin{array}{cc}
{\rm e}^{\Lambda t} & 0\\ 0 & 1\end{array}\right)S,\quad
t\in\mathbb{R},
\end{equation}
где $\Phi(t)$ -- $T$-периодическая $n\times n$ матрица Флоке,
$\Lambda$ -- $(n-1)\times (n-1)$ невырожденная матрица и $S$ --
подходящая невырожденная  $n\times n$-матрица. В силу
(\ref{orient0}) и (\ref{orient}) последний столбец матрицы
$\left((Z_{n-1}(t),\widetilde{z}(t))^{-1}\right)^*$ -- это
$\dot{\widetilde{x}}(t),$ $t\in\mathbb{R}.$ Из этого утверждения
следует, что $S$ имеет вид
$S=\left(S_{n-1},\left(\begin{array}{c} 0_{n-1\times 1} \\
s_0\end{array}\right)\right),$ где $S_{n-1}$ -- $n\times
(n-1)$-матрица и $s_0\not=0.$ Из этого мы можем заключить, что
матрица  $S_{n-1}$ не содержит линейно-зависимых с $\left(\begin{array}{c} 0_{n-1\times 1} \\
 s_0\end{array}\right)$ столбцов.
Следовательно, по крайней мере одна не $n$-я компонента любого
столбца  $S_{n-1}$ является ненулевой, что, очевидно, завершает
доказательство второго утверждения леммы.

Для доказательства третьего утверждения леммы заметим, что формула
(\ref{teflo}) позволяет утверждать, что
$$
\left(\begin{array}{cc} {\rm e}^{\Lambda T} & 0\\ 0 &
1\end{array}\right)S\left((Z_{n-1}(t+T),\widetilde{z}(t+T))^{-1}\right)^*=$$
$$=
S\left((Z_{n-1}(t),\widetilde{z}(t))^{-1}\right)^*,\quad
t\in\mathbb{R}.
$$
Но из последней формулы следует, что
$$
  \left.{\rm e}^{\Lambda
  T}S_{n-1}\right|_{\mathbb{R}^{n-1}}Z_{n-1}^*(t+T)=\left.S_{n-1}\right|_{\mathbb{R}^{n-1}}Z_{n-1}^*(t),
$$
где через $\left.S_{n-1}\right|_{\mathbb{R}^{n-1}}$ обозначена
матрица, составленная из первых $n-1$ строк матрицы $S_{n-1}.$
Таким образом, для завершения доказательства третьего утверждения
достаточно положить
$$
  \widetilde{D}:=\left(\left.S_{n-1}\right|_{\mathbb{R}^{n-1}}\right)^{-1}
  \left.{\rm e}^{\Lambda
  T}S_{n-1}\right|_{\mathbb{R}^{n-1}},
$$
при этом обращение матрицы
$\left.S_{n-1}\right|_{\mathbb{R}^{n-1}}$ возможно в силу
обратимости матрицы $S$ и установленного выше ее вида.

Лемма доказана.

Определим  непрерывно дифференцируемые функции
$h,I:\mathbb{R}\times\mathbb{R}^{n-1}\to\mathbb{R}^n$ как
$$
  h(\theta,r)=\widetilde{x}(\theta)+Y_{n-1}(\theta)r
$$
$$
  I(\theta,r)=\Omega(T,0,h(\theta,r)),
$$
где  $\Omega$ -- оператор сдвига по траекториям системы
(\ref{ps_3}).

\begin{lem}\label{conv_lem2}
Предположим (\ref{orient0}) и (\ref{orient}). Тогда
$\dot{\widetilde{x}}(\theta)\not\in
I'_r(\theta,0)\mathbb{R}^{n-1}$ при каждом $\theta\in[0,T].$
\end{lem}

\noindent {\bf Доказательство. }
 Предположим противное, тогда существует
  $r\in\mathbb{R}^{n-1},$ $r\not=0$ такое, что
$\dot{\widetilde{x}}(\theta)=I'_r(\theta,0)r.$
 Имеем
 $$
  \dot{\widetilde{x}}(\theta)=I'_r(\theta,0)r=\Omega'_{(3)}(T,0,h(\theta,0))Y_{n-1}(\theta)r=Y_{n-1}(\theta+T)r,
$$
$$
  \dot{\widetilde{x}}(\theta-T)=Y_{n-1}(\theta)r.
$$
Но $\dot{\widetilde{x}}(\theta)=\dot{\widetilde{x}}(\theta-T),$ и
получаем противоречие с утверждением леммы~\ref{conv_lem1}.

Лемма доказана.

\begin{cor}\label{conv_cor1} Предположим, что  $T>0$ -- наименьший период цикла
$\widetilde{x}.$ Тогда существует  $r_0>0$ такое, что
$$
  I\left(\theta,B_{r_0}(0)\right)\cap\widetilde{x}\left([0,T]\right)=\left\{\widetilde{x}(\theta)\right\}.
$$
\end{cor}

\begin{cor}\label{conv_cor2} Предположим, что  $T>0$ -- наименьший период цикла
$\widetilde{x}.$ Тогда существует $k_0>0$ такое, что
$$
  I\left(\theta,B_r(0)\right)\cap x_{\varepsilon_k}\left([0,T]\right)=
  \left\{x_{\varepsilon_k}(\theta-\Delta_{\varepsilon_k}(\theta))\right\},\qquad
  k>k_0,
$$
где $\Delta_{\varepsilon_k}:\mathbb{R}\to\mathbb{R}$ --
непрерывная функция такая, что $\Delta_{\varepsilon_k}(\theta)\to
0$ при $k\to \infty$ равномерно по отношению к
$\theta\in\mathbb{R}.$
\end{cor}

\noindent {\bf Доказательство теоремы \ref{conv_thm1}.}
 Сделаем в системе (\ref{ps_3}) замену переменных
$\nu_k(t,\theta)=\Omega(0,t,x_k(t-\Delta_{\varepsilon_k}(\theta)+\theta)),$
где $\Delta_{\varepsilon_k}(\theta)$ -- числа, о которых говорится
в следствии~\ref{conv_cor2}. Заметим, что
$x_k(t-\Delta_{\varepsilon_k}(\theta)+\theta)=\Omega(t,0,\nu_k(t,\theta))$
и, таким образом,
\begin{equation}\label{ob1}
 \dot x_k(t-\Delta_{\varepsilon_k}(\theta)+\theta)=
 f(\Omega(t,0,\nu_k(t,\theta))+\Omega'_{\xi}(t,0,
 \nu_k(t,\theta)) (\nu_k)'_t(t,\theta).
\end{equation}
С другой стороны, из (\ref{ps_3}) имеем
$$
  \dot
  x_k(t-\Delta_{\varepsilon_k}(\theta)+\theta)=
$$
\begin{equation}\label{ob2}
=  f(\Omega(t,0,\nu_k(t,\theta)))+{\varepsilon_k}
  g(t-\Delta_{\varepsilon_k}(\theta)+\theta,
  \Omega(t,0,\nu_k(t,\theta)),{\varepsilon_k}).
\end{equation}
Из (\ref{ob1}) и (\ref{ob2}) следует, что
$$
  (\nu_k)'_t(t,\theta)={\varepsilon_k}\left(\Omega'_{\xi}
  (t,0,\nu_k(t,\theta))\right)^{-1}g(t-\Delta_{\varepsilon_k}
  (\theta)+\theta,\Omega(t,0,\nu_k(t,\theta)),{\varepsilon_k}),
$$
и так как
$$
  \nu_k(0,\theta)=x_k(-\Delta_{\varepsilon_k}(\theta)+
  \theta)=x_k(T-\Delta_{\varepsilon_k}(\theta)+\theta)=
  \Omega(T,0,\nu_k(T,\theta)),
$$
окончательно получаем
$$
  \nu_k(t,\theta)=\Omega(T,0,\nu_k(T,\theta))+{\varepsilon_k}\int\limits_0^t
  \left(\Omega'_{\xi}(\tau,0,\nu_k(\tau,\theta))\right)^{-1}\circ
$$
\begin{equation}\label{ob3}
\circ g(\tau-
  \Delta_{\varepsilon_k}(\theta)+\theta,\Omega(\tau,0,\nu_k(\tau,\theta)),{\varepsilon_k})
  d\tau.
\end{equation}
Так как $\nu_k(t,\theta)\to \widetilde{x}(\theta)$ при $k\to
\infty$ можем записать $\nu_k(t,\theta)$ в виде
\begin{equation}\label{rep}
  \nu_k(t,\theta)=\widetilde{x}(\theta)+{\varepsilon_k}
  \mu_k(t,\theta).
\end{equation}
Докажем теперь, что функции $\mu_k$ равномерно ограниченны по
отношению к $k\in\mathbb{N}.$ Для этого, во-первых, вычтем
$\widetilde{x}(\theta)$ из обеих частей (\ref{ob3}), получив
\begin{eqnarray}
{\varepsilon_k}\mu_k(t,\theta)&=&{\varepsilon_k}\,\Omega'_{\xi}(T,0,\widetilde{x}(\theta))
\mu_k(T,\theta)+o({\varepsilon_k}\mu_k(T,\theta))+\nonumber\\
& &
+{\varepsilon_k}\int\limits_0^t\left(\Omega'_{\xi}(\tau,0,\nu_k(\tau,\theta))
\right)^{-1}g(\tau-\Delta_{\varepsilon_k}(\theta)+\nonumber\\
& &+\theta, \Omega(\tau,0,\nu_k(\tau,\theta)),{\varepsilon_k})
  d\tau.\label{ob4}
\end{eqnarray}
Так как $x_k(-\Delta_{\varepsilon_k}(\theta)+\theta)\in
I\left(\theta,B_{r_0}^{n-1}(0)\right),$ то по определению $I$
существует $r_k(\theta)\in B^{n-1}_{r_0}(0)$ такое, что
\begin{equation}\label{Q}
x_k(-\Delta_{\varepsilon_k}(\theta)+\theta)=
\Omega(T,0,h(\theta,r_k(\theta))),
\end{equation} причем на основании
(\ref{as1_3}) имеем
\begin{equation}\label{convr}
r_k(\theta)\to 0\ \ \mbox{\rm при}\ \ k\to \infty.
\end{equation}
Пользуясь равенством (\ref{Q}), получаем следующее представление
для $\varepsilon_k\mu_k$
\begin{eqnarray}{\varepsilon_k}\mu_k(T,\theta) & =
& \nu_k(T,\theta)-\widetilde{x}(\theta)=\Omega(0,T,x_k
(-\Delta_{\varepsilon_k}(\theta)+\theta))-\widetilde{x}(\theta)=\nonumber\\
& & =
\Omega(0,T,\Omega(T,0,h(\theta,r_k(\theta))))-\widetilde{x}(\theta)=\nonumber\\
& & =
h(\theta,r_k(\theta))-\widetilde{x}(\theta)=Y_{n-1}(\theta)r_k(\theta).
\nonumber
\end{eqnarray}
Следовательно, формула (\ref{ob4}) при $t=T$ может быть переписана
как
\begin{eqnarray}
Y_{n-1}(\theta)r_k(\theta)&=&Y_{n-1}(T+\theta)r_k(\theta)+
o(Y_{n-1}(\theta)r_k(\theta))+\nonumber\\
& &+
{\varepsilon_k}\int\limits_0^T\left(\Omega'_{\xi}(\tau,0,\nu_k(\tau,\theta))
\right)^{-1}g(\tau-\Delta_{\varepsilon_k}(\theta)+\nonumber\\
& &+\theta, \Omega(\tau,0,\nu_k(\tau,\theta)),{\varepsilon_k})
  d\tau.\label{bisp}
\end{eqnarray}
На основании (\ref{bisp})  сейчас будет установлено существование
такого $c>0,$ что
\begin{equation}\label{estt}
  \|r_k(\theta)\|\le {\varepsilon_k}
  c,\quad k\in\mathbb{N},\ \theta\in[0,T].
\end{equation}
Предположим противное, тогда можем считать, что
$\{\theta_k\}_{k\in\mathbb{N}}\subset[0,T],$ $\theta_k\to\theta_0$
при $k\to\infty,$ -- такая последовательность, что
$\|r_k(\theta_k)\|=\varepsilon_k c_k,$ где $c_k\to\infty$ при
$k\to\infty.$ Положим
$q_k=\dfrac{r_k(\theta_k)}{\|r_k(\theta_k)\|},$ тогда из
(\ref{bisp}) имеем
\begin{eqnarray}\label{ff}
 Y_{n-1}(\theta_k)q_k&=&Y_{n-1}(T+\theta_k) q_k+\frac{o(Y_{n-1}(\theta_k)r_{\varepsilon_k}(\theta_k))}
 {\|r_{\varepsilon_k}(\theta_k)\|}+\nonumber\\
 &  &+\frac{1}{c_k}\int\limits_0^t
 \left(\Omega'_{\xi}(\tau,0,\nu_k(\tau,\theta_k))
 \right)^{-1}g(\tau-\theta_0+\Delta_{\varepsilon_k}(\theta_k)+\nonumber\\
 & &+
 \theta_k,\Omega(\tau,0,\nu_k(\tau,\theta_k)),\varepsilon_k)
  d\tau.
\end{eqnarray}
Без ограничения общности можем считать, что последовательность
$\{q_k\}_{k\in\mathbb{N}}$ сходится, положим
$q_0=\lim_{k\to\infty}q_k,$ тогда $\|q_0\|=1.$ С другой стороны,
из (\ref{ff}) имеем $Y_{n-1}(\theta_0)q_0=Y_{n-1}(T+\theta_0)
q_0,$ то есть приходим к противоречию с утверждением~2
леммы~\ref{conv_lem1}. Таким образом, (\ref{estt}) выполнено для
некоторого $c>0$ и функции $\mu_k$ равномерно ограниченны по
отношению к $k\in\mathbb{N}.$ Из (\ref{rep}) также заключаем
\begin{equation}\label{wehav}
  x_k(t-\Delta_{\varepsilon_k}(\theta)+\theta)-\widetilde{x}(t+\theta)={\varepsilon_k}\mu_k(t,\theta).
\end{equation}

\vskip3mm

\noindent Следовательно,
$x_k(\theta-\Delta_{\varepsilon_k}(\theta))\to
\widetilde{x}(\theta)$ со скоростью ${\varepsilon_k}>0.$

\vskip3mm \noindent  Для завершения доказательства
теоремы~\ref{conv_thm1} остается установить  (\ref{mainprop}). Для
этого введем новые функции
$a_k:\mathbb{R}\times\mathbb{R}\to\mathbb{R}$ и
$b_k:\mathbb{R}\times\mathbb{R}\to\mathbb{R}^{n-1}$ согласно
следующим формулам
$$
  a_k(t,\theta)=\widetilde{z}^*(t+\theta)\mu_k(t,\theta),
$$
\begin{equation}\label{b}
  b_k(t,\theta)=Z_{n-1}^*(t+\theta)\mu_k(t,\theta).
\end{equation}
Пользуясь утверждением~1 леммы~\ref{conv_lem1}, можем представить
$x_k(t-\Delta_{\varepsilon_k}(\theta)+\theta)-\widetilde{x}(t+\theta)$
в виде
$$
 x_k(t-\Delta_{\varepsilon_k}(\theta)+\theta)-\widetilde{x}(t+\theta)=
$$
\begin{equation}\label{fw}
 ={\varepsilon_k}
  \dot {\widetilde{x}}(t+\theta) a_k(t,\theta)+{\varepsilon_k} Y_{n-1}(t+\theta)b_k(t,\theta).
\end{equation}
Вычитая (\ref{np_3}), где $x(t)$ заменено функцией
$\widetilde{x}(t+\theta),$ из  (\ref{ps_3}), где $x(t)$ заменено
функцией $x_k(t-\Delta_{\varepsilon_k}(\theta)+\theta),$ получаем
$$
  \dot x_k(t-\Delta_{\varepsilon_k}(\theta)+\theta)-\dot
  {\widetilde{x}}(t+\theta)=
$$
$$
  =f'(\widetilde{x}(t+\theta))(x_k(t-\Delta_{\varepsilon_k}(\theta)+\theta)-\widetilde{x}(t+\theta))+
$$
\begin{eqnarray}\label{bis1}
   +{\varepsilon_k} g(t-\Delta_{\varepsilon_k}(\theta)+\theta,x_k(t-\Delta_{\varepsilon_k}(\theta)+\theta),{\varepsilon_k})+\widetilde{
  o}_t({\varepsilon_k}),
\end{eqnarray}
где функция $\widetilde {o}_t(\cdot)$ имеет те же свойства, что и
функция $o(\cdot),$ введенная ранее, более того $\widetilde{
o}_{t+T}(\cdot)=\widetilde {o}_t(\cdot)$ для любого
$t\in\mathbb{R}.$ Подставляя (\ref{fw}) в (\ref{bis1}) и учитывая,
что
$$
f'(\widetilde{x}(t+\theta)){\varepsilon_k} \dot
{\widetilde{x}}(t+\theta) a_k(t,\theta)={\varepsilon_k} \ddot
{\widetilde{x}}(t+\theta) a_k(t,\theta)
$$
и
$$
f'(\widetilde{x}(t+\theta)){\varepsilon_k} Y_{n-1}(t+\theta)
b_k(t,\theta)={\varepsilon_k}
\dot{Y}_{n-1}(t+\theta)b_k(t,\theta),
$$
получаем
$$
  {\varepsilon_k}
\dot {\widetilde{x}}(t+\theta) (a_k)'_t(t,\theta)
  +{\varepsilon_k} \dot{Y}_{n-1}(t+\theta)({b_k})'_t(t,\theta)=
$$
$$=g(t-\Delta_{\varepsilon_k}(\theta)+\theta,x_k(t-\Delta_{\varepsilon_k}(\theta)+\theta),\varepsilon_k)+\widetilde
 { o}_t({\varepsilon_k}).
$$
Из предыдущего равенства имеем
$$
  {\varepsilon_k}\left(
  b_k\right)'_t(t,\theta)=$$
$$={\varepsilon_k}
Z_{n-1}(t+\theta)g(t-\Delta_{\varepsilon_k}(\theta)+\theta,x_k(t-\Delta_{\varepsilon_k}(\theta)+\theta),\varepsilon_k)+Z_{n-1}(t+\theta)\widetilde
{o}_t({\varepsilon_k})
$$
и, следовательно,
$$
  b_{k}(0,\theta)=b_{k}(-T,\theta)+\int\limits_{-T}^0
  Z_{n-1}^*(\tau+\theta)g(\tau-\Delta_{\varepsilon_k}(\theta)+
$$
\begin{equation}\label{prob}
  +\theta,
  x_k(\tau-\Delta_{\varepsilon_k}(\theta)+\theta),{\varepsilon_k})d\tau+\frac{o(\varepsilon_k)}{\varepsilon_k}.
\end{equation}
Из (\ref{prob}) и утверждения~3 леммы~\ref{conv_lem1} получаем
\begin{eqnarray}
  b_{k}(0,\theta)& =&(I-\widetilde{D})^{-1}\int\limits_{-T}^0 Z_{n-1}^*(\tau+\theta)g(\tau-\Delta_{\varepsilon_k}(\theta)+\nonumber\\
  & &+\theta,
  x_k(\tau-\Delta_{\varepsilon_k}(\theta)+\theta),{\varepsilon_k})d\tau+
  \frac{o(\varepsilon_k)}{\varepsilon_k}\nonumber
\end{eqnarray}
или, вводя замену переменных $\tau+\theta=s$ в интеграле,
$$
  b_{k}(0,\theta) =(I-\widetilde{D})^{-1}\int\limits_{\theta-T}^\theta
Z_{n-1}^*g(s,
  x_k(s),{\varepsilon_k})d\tau+
  \frac{o(\varepsilon_k)}{\varepsilon_k}.
$$
 С другой стороны, из (\ref{fw}) имеем
$$
\left<z_i(t+\theta),x_k(t-\Delta_{\varepsilon_k}(\theta)+\theta)-\widetilde{x}(t+\theta)\right>={\varepsilon_k}
b_k(t,\theta)
$$
и, таким образом, для завершения доказательства достаточно
положить $D:=(I-\widetilde{D})^{-1}.$

Теорема доказана.

Рассмотрим теперь некоторые приложения формулы  (\ref{mainprop}) к
описанию поведения решений $x_k,$ когда $k\to\infty.$ Ниже
$\angle(a,b)$ -- угол между векторами $a,b\in\mathbb{R}^n,$
принадлежащий отрезку $[0,\pi].$

\begin{cor}\label{conv_cor3} Пусть выполнены все предположения теоремы~\ref{conv_thm1}.
Тогда для любых  $i\in\overline{1,n-1}$ и $\theta\in[0,T]$ таких,
что $M^\bot_i(\theta)\not=0,$ существует $j\in\overline{1,n-1}$,
при котором
$$
\cos\angle\left(z_j(\theta),
x_k(\theta-\Delta_{\varepsilon_k}(\theta))-
  \widetilde{x}(\theta)
  \right)\not=0
$$
для достаточно больших $k\in\mathbb{N}.$
\end{cor}

Доказательство следствия вытекает из формулы
$$
  \left\|z_i(\theta)\right\|\cdot \left\|x_k(\theta-\Delta_{\varepsilon_k}(\theta))-
  \widetilde{x}(\theta)\right\|\cdot \cos\angle\left(z_i(\theta),
x_k(\theta-\Delta_{\varepsilon_k}(\theta))-
  \widetilde{x}(\theta)
  \right)=
$$
\begin{equation}\label{mainprop1}
=\left[{\varepsilon_k} D
M^\bot(\theta)\right]^i+o({\varepsilon_k}),
\end{equation}
получаемой подстановкой выражения
$$
  \left<z_i(\theta),x_k(\theta-\Delta_{\varepsilon_k}(\theta))-\widetilde{x}(\theta)\right>=
$$
$$=  \left\|z_i(\theta)\right\|\cdot
\left\|x_k(\theta-\Delta_{\varepsilon_k}(\theta))-
  \widetilde{x}(\theta)\right\|\cdot \cos\angle\left(z_i(\theta),
x_k(\theta-\Delta_{\varepsilon_k}(\theta))-
  \widetilde{x}(\theta)
  \right)
$$
в главную формулу (\ref{mainprop}). Действительно, если
$M^\bot_i(\theta)\not=0,$ то существует $j\in\overline{1,n-1}$
такое, что $\left[DM^\bot(\theta)\right]^j\not=0.$

Следующий результат является непосредственным следствием формулы
(\ref{mainprop1}).

\begin{cor}\label{conv_cor4} Пусть выполнены все предположения теоремы~\ref{conv_thm1}.
Если существует по крайней мере одно $i_*\in\overline{1,n-1}$
такое, что  $M^\bot_{i_*}(\theta)\not=0,$ то
$$
  c_1{\varepsilon_k}\le\left\|x_k(\theta-\Delta_{\varepsilon_k}(\theta))-
  \widetilde{x}(\theta)\right\|\le
  c_2{\varepsilon_k}
$$
для некоторых  $0<c_1\le c_2,$ любых  $\theta\in[0,T]$ и $k\ge
k_0,$ где $k_0\in\mathbb{N}$ достаточно велико.
\end{cor}

Учитывая следствие~\ref{conv_cor2}, можем получить теперь
следующий факт.

\begin{cor}\label{conv_cor5} Пусть выполнены все предположения теоремы~\ref{conv_thm1}.
Если существует по крайней мере одно $i_*\in\overline{1,n-1}$
такое, что $M^\bot_{i_*}(\theta)\not=0,$ то
$$ x_k(t)\not=\widetilde{x}(\theta)\quad \mbox{ для
любых } t,\theta\in[0,T]$$ при условии, что $k\in\mathbb{N}$
достаточно велико.
\end{cor}

\noindent {\bf Доказательство.} Пусть $k_0>0$ -- то, о котором
говорится в следствии~\ref{conv_cor4}, и предположим, что
существуют $\varepsilon_{k_*},$ $k_*>k_0,$ $\theta_*,t_*\in[0,T]$
такие, что $x_{k_*}(t_*)=\widetilde{x}(\theta_*).$ Но согласно
следствию~\ref{conv_cor2} должно быть
$t_*=x_{k_*}(\theta_*-\Delta_{\varepsilon_{k_*}}(\theta_*)),$
противореча утверждению следствия~\ref{conv_cor4}.

Следствие доказано.

Для того, чтобы получить теперь достаточные условия,
обеспечивающие сходимость со скоростью большей, чем
${\varepsilon_k}>0,$ нам необходим следующий вспомогательный
результат.

\begin{lem}\label{conv_lem3} Пусть $k_0\in\mathbb{N}$ достаточно велико. Тогда для
любых $k>k_0$ и $\theta\in[0,T],$ удовлетворяющих условию
$$x_k(\theta-\Delta_{\varepsilon_k}(\theta))\not=\widetilde{x}(\theta),$$
существует $i_*\in\overline{1,n-1}$ такое, что
$$\left|\cos\angle\left(z_{i_*}(\theta),x_k(\theta-\Delta_{\varepsilon_k}(\theta))-\widetilde{x}(\theta)\right)
\right|\ge\alpha_*,$$ где  $\alpha_*>0$ не зависит от $k$ и
$\theta.$
\end{lem}

\noindent {\bf Доказательство.} Предположим противное, тогда можем
считать, что  последовательности
$\{\varepsilon_k\}_{k\in\mathbb{N}},$
$\{\theta_k\}_{k\in\mathbb{N}}\subset[0,T],$ $\theta_k\to
\theta_0$ при $k\to\infty,$
$\{\alpha_k\}_{k\in\mathbb{N}}\subset[-1,1],$ $\alpha_k\to 0$ при
$k\to\infty,$ $\{r_k\}_{k\in\mathbb{N}}\subset(0,1],$ $r_k\to 0$
при $k\to\infty$ таковы, что
$$
\cos\angle\left(z_{i}(\theta),I(\theta_k,r_k)-I(\theta_k,0)\right)=\alpha_k,\qquad
i\in\overline{1,n-1},
$$
или равносильно
$$
\frac{\left<z_i(\theta_k),I(\theta_k,r_k)-I(\theta_k,0)\right>}
{\left\|z_i(\theta_k)\right\|\cdot\left\|I(\theta_k,r_k)-I(\theta_k,0)\right\|}=\alpha_k,\qquad
i\in\overline{1,n-1},
$$
поэтому
\begin{equation}\label{ptl}
\frac{\left<z_i(\theta_n),I'_r(\theta_k,r_k)\frac{r_k}{\|r_k\|}+
\frac{o(r_k)}{\|r_k\|}\right>}
{\left\|z_i(\theta_k)\right\|\cdot\left\|I'_r(\theta_k,r_k)\frac{r_k}{\|r_k\|}+
\frac{o(r_k)}{\|r_k\|}\right\|}=\alpha_k,\qquad
i\in\overline{1,n-1}.
\end{equation}
Без ограничения общности можем считать, что  $\frac{r_k}{\|r_k\|}$
сходится при $k\to\infty.$ Положим
$q_0=\lim_{k\to\infty}\frac{r_k}{\|r_k\|},$ тогда $\|q_0\|=1.$
Переходя к пределу при $k\to\infty$ в (\ref{ptl}), получаем
\begin{equation}\label{zvio}
  \left<z_i(\theta_0),I'_r(\theta_0,0)q_0\right>=0\quad\mbox{для любых
}i\in\overline{1,n-1}.
\end{equation}
Разложим $I'_r(\theta_0,0)q_0$ следующим образом
$$
  I'_r(\theta_0,0)q_0=y_1(\theta_0)a^1+...+y_{n-1}(\theta_0)a^{n-1}+\dot{\widetilde{x}}(\theta)a^n.
$$
Из (\ref{zvio}) имеем, что $a^1=...=a^{n-1}=0$ и, таким образом,
$I'_r(\theta,0)q_0=a^n\dot{\widetilde{x}}(\theta),$ противореча
утверждению леммы~\ref{conv_lem2}.

Лемма доказана.

Суммируя теорему~\ref{conv_thm1} и лемму~\ref{conv_lem3}, получаем
утверждение.

\begin{cor}\label{conv_cor6} Пусть выполнены все предположения теоремы~\ref{conv_thm1}.
Предположим, что $M^\bot_i(\theta)=0$ для любого
$i\in\overline{0,n-1}.$ Тогда
$$
  \left\|x_k(\theta-\Delta_{\varepsilon_k}(\theta))-\widetilde{x}(\theta)\right\|=o({\varepsilon_k}).
$$
\end{cor}

\noindent {\bf Доказательство.} Предположим противное, тогда можем
считать, что последовательности
$\{\varepsilon_k\}_{k\in\mathbb{N}},$
$\{\theta_k\}_{k\in\mathbb{N}}\subset[0,T],$ $\theta_k\to
\theta_0$ при $k\to\infty$ и $c_*>0$ таковы, что
\begin{equation}\label{cor6}
  \frac{\left\|x_{\varepsilon_k}(\theta_k-\Delta_{\varepsilon_k}(\theta_k))-\widetilde{x}(\theta_k)\right\|}{\varepsilon_k}\ge c_*.
\end{equation}
Из  (\ref{cor6}) имеем, что предположения леммы~\ref{conv_lem3}
удовлетворены, пусть $i_*\in\overline{1,n-1}$ -- то число, о
котором говорится в этой лемме. Но  (\ref{cor6}) противоречит
соотношению (\ref{mainprop1}), когда $i:=i_*,$ что завершает
доказательство требуемого утверждения.

Следствие доказано.

Следствия \ref{conv_cor3} и \ref{conv_cor6} позволяют
сформулировать следующую альтернативу.

\begin{cor}\label{conv_cor7} Пусть выполнены все предположения
теоремы~\ref{conv_thm1}.  Положим $\theta_*\in[0,T].$ Тогда либо
существует $i_*\in\overline{1,n-1}$ такое, что
$$
\cos\angle\left(z_{i_*}(\theta),
x_k(\theta_*-\Delta_{\varepsilon_k}(\theta_*))-
  \widetilde{x}(\theta_*)
  \right)\not=0
$$
для достаточно больших $k\in\mathbb{N},$ либо
$$
  \left\|x_k(\theta_*-\Delta_{\varepsilon_k}(\theta_*))-\widetilde{x}(\theta_*)\right\|=o({\varepsilon_k}).
$$
\end{cor}

\section{Сопоставление полученных результатов с имеющимися в
литературе} В случае, когда дополнительно известно, что функция
$g$ непрерывно дифференцируема и линеаризованная система
(\ref{ls_3}) не имеет единичных мультипликаторов, сходимость в
(\ref{as1_3}) со скоростью $\varepsilon_k>0$ следует из формулы
разложения решений $x_k$ в ряд по степеням $\varepsilon_k,$
даваемой методом малого параметра Пуанкаре (см. Б.~П.~Демидович
\cite{dem}, Гл.~III, \S~24, М.~Розо \cite{rozo}, Гл.~9, \S~1).
Если относительно системы (\ref{ls_3}) известно существование
мультипликатора $+1$ алгебраической кратности $1,$ то сходимость в
(\ref{as1_3}) со скоростью $\varepsilon_k>0$ для случая
аналитических и дважды непрерывно дифференцируемых правых частей
системы (\ref{ps_3}) доказана соответственно И.~Г.~Малкиным (см.
\cite{mal}, формула~4.1) и В.~С.~Лудом (см. \cite{loud},
формула~1.3 теоремы~1). Однако в работах указанных авторов не
отмечался установленный в главе факт (следствие~3.8) о том, что
скорость сходимости может иметь и больший, чем $\varepsilon_k>0$
порядок. Также в классических работах не указаны свойство
(\ref{claa_3}) и следствия~3.4-3.6, связанные с качественными
свойствами поведения $T$-периодических решений системы
(\ref{ps_3}) при $\varepsilon\to 0.$ Если же возмущение всего лишь
непрерывно, что является общим предположением теорем о
существовании, доказанных в \cite{kam}, \cite{alex3},
\cite{alex2}, \cite{kraper}, \cite{krastr}, \cite{mit},
\cite{muh3}-\cite{muh1}, \cite{sam}, \cite{camaza}-\cite{dancer},
\cite{zan}-\cite{laz}, \cite{mawphd1}-\cite{ortega}, а также
утверждений глав~1 и 2, то результаты о скорости сходимости
$T$-периодических решений системы (\ref{ps_3}) при $\varepsilon\to
0$  в литературе отсутствуют. Предложенные в настоящей главе
теоремы частично заполняют этот пробел.

Теорема~\ref{alter} опубликована автором в \cite{asim}, где также
проведено ее доказательство для $t=T.$

\def\bibname{Список литературы}


\end{document}